\documentclass[11pt]{amsart}
\usepackage{amsmath}
\usepackage{amssymb}
\usepackage{mathrsfs}
\usepackage{esint}
\usepackage{tikz}
\usepackage{comment}
\usepackage{ulem}
\usepackage{xcolor}
\usepackage{hyperref}

\theoremstyle{plain}
\newtheorem{theorem}{Theorem}[section]
\newtheorem{lemma}[theorem]{Lemma}
\newtheorem{prop}[theorem]{Proposition}
\newtheorem{cor}[theorem]{Corollary}

\theoremstyle{definition}
\newtheorem{remark}[theorem]{Remark}
\newtheorem{definition}[theorem]{Definition}

\numberwithin{equation}{section}

\makeatletter
\@namedef{subjclassname@2020}{%
  \textup{2020} Mathematics Subject Classification}
\makeatother

\def\be{\begin{equation}}
\def\ee{\end{equation}}

\def\R{\mathbb R}

\def\Z{\mathbb Z}
\def\C{\mathbb C}

\def\epsilon{\varepsilon}

\def\({\left(}
\def\){\right)}

\begin{document}

\title[Geometric Differential Equations near Infinity]
{Solution Spaces of the Minimal Surface and Monge-Amp\`{e}re Equations near Infinity}

\author[Han]{Qing Han}
\address{Department of Mathematics\\
University of Notre Dame\\
Notre Dame, IN 46556, USA} 
\email{qhan@nd.edu}
\author[Hung]{Jui-Yun Hung}
\address{Mathematical Sciences Institute\\
The Australian National University\\
Canberra, ACT 2600, Australia} 
\email{JUI-YUN.HUNG@anu.edu.au}


\begin{abstract}
In this paper, we study the constructions of all solutions to the minimal surface equation 
and the Monge-Amp\`{e}re equation near infinity
and the parametrization of these solutions by harmonic functions. 
It is well-known that Dirichlet boundary value problems in exterior domains
do not always admit solutions.
We demonstrate that prescribing asymptotic behaviors forms a new type of problem
leading to all solutions near infinity. The harmonic functions determining the asymptotic behaviors
play the role of ``free data" as the boundary values in the boundary value problems.
\end{abstract}

\thanks{Q. Han acknowledges the support of NSF Grant DMS-2305038.}


\maketitle

\section{Introduction}\label{sec-Intro}

Classification of solutions in the entire space and in exterior domains forms an important topic in differential equations. 
For differential equations of significant geometric background, 
global solutions are linear functions or quadratic polynomials 
and solutions in exterior domains are asymptotic to those global solutions near infinity, 
with asymptotic behaviors described precisely. 
In this paper, we study whether there always exists a solution with a prescribed asymptotic behavior. 

We start with the minimal surface equation. 
A classical result asserts that any minimal graph over $\R^n$
is a hyperplane, due to Bernstein \cite{Bern1015}  for $n=2$,
De Giorgi \cite{DeGiorgi1965} for $n=3$,  Almgren \cite{Alm1966} for $n=4$, 
and Simons \cite{Simons1968} for $n=5, 6, 7$. 
Bombieri, De Giorgi, and Giusti \cite{BGG1969} demonstrated that 
there are minimal graphs but not hyperplanes over $\R^n$ for $n\ge 8$. 
On the other hand, Moser \cite{Mo} proved that 
any minimal graphs over $\mathbb R^n$ with bounded slopes are hyperplanes, for any $n\ge 2$.

We now consider minimal graphs defined near infinity. 
Let $u$ be a solution of the minimal surface equation outside a bounded subset in the Euclidean space $\mathbb R^n$. 
Bers \cite{BL} proved that $\nabla u(x)$ has a limit 
as $x\rightarrow\infty$ for $n=2$, 
and Simon \cite{Sim} extended this result to dimension $3\leq n\leq 7$.
In arbitrary dimensions,  Moser \cite{Mo} proved that, if $|\nabla u|$ is bounded near infinity, then 
$\nabla u(x)$ has a limit 
as $x\rightarrow\infty$. 
Concerning asymptotic expansions near infinity, Bers \cite{BL} proved that, for $n=2$, 
there exist constants $a, c_1, c_2$ such that 
\begin{equation}\label{eq-MSE-expansion-dim-2-Bers-intro} 
u(x)=c_1x_1+c_2x_2+a\log|x|+O(1)\quad\text{as } x\to\infty.\end{equation} 
Here and hereafter, the expression $O(|x|^{-k})$ is used to indicate 
a term that is bounded in absolute value by a constant multiple of $|x|^{-k}$ for $|x|$ large. 
For higher dimensions $n\ge 3$,  under the additional assumption $\nabla u(x)\to0$ as $x\to\infty$,  
Schoen \cite{Sc} proved that
there exist constants $a, b, c_j\in\mathbb{R}$, such that
\begin{equation}\label{eq-MSE-expansion-dim-higher-Schoen-intro} 
u(x)=a+b|x|^{2-n}+\sum\limits_{j=1}^nc_jx_j|x|^{-n}+O(|x|^{-n}) \quad\text{as } x\to\infty.\end{equation} 
The estimates \eqref{eq-MSE-expansion-dim-2-Bers-intro} and \eqref{eq-MSE-expansion-dim-higher-Schoen-intro} 
demonstrate that the fundamental solution of the Laplace operator provides the correct approximation of the next order 
after linear functions.  
Refer to \cite{Simon1994} for more related results. 

It is possible to follow Bers and Schoen to approximate the solutions of the minimal surface equation near infinity up 
to arbitrary orders. 
Han and Wang \cite{Han-Wang2024} took a different approach and used a single function to 
identify behaviors of solutions near infinity precisely. 
Their result can be summarized as follows. 
Let $u\in C^2(\R^n\setminus \overline{B}_1)$ be a solution of 
\begin{equation}\label{eq-MSE-intro} 
\Delta u-\frac{u_iu_j}{1+|\nabla u|^2}u_{ij}=0\quad\text{in }\mathbb R^n\setminus\overline{B}_1.\end{equation}
Assume $\nabla u(x)\to 0$ as $x\to\infty$. Then,

$\mathrm{(i)}$ for $n =2$, there exist a constant $a$ and a function $v\in C^\infty(B_1)$ such that 
\begin{equation}\label{eq-MSE-expansion-identity-dim-2-intro}
u(x)=a\log|x|+v\Big(\frac{x}{|x|^2}\Big)\quad\text{for any }x\in \R^2\setminus\overline{B}_1;\end{equation}

$\mathrm{(ii)}$ for $n \geq 3$, there exist a constant $a$ and a function $v\in C^\infty(B_1)$ such that 
\begin{equation}\label{eq-MSE-expansion-identity-dim-higher-intro}
u(x)=a+|x|^{2-n}v\Big(\frac{x}{|x|^2}\Big)\quad\text{for any }x\in \mathbb R^n\setminus\overline{B}_1.\end{equation}

Without any effort, we can always introduce a function $v\in C^\infty(B_1\setminus\{0\})$ such that 
\eqref{eq-MSE-expansion-identity-dim-2-intro} or \eqref{eq-MSE-expansion-identity-dim-higher-intro} holds. 
The main feature here is the assertion 
that $v$ is a smooth function in $B_1$. 
The decompositions in \eqref{eq-MSE-expansion-identity-dim-2-intro} and 
\eqref{eq-MSE-expansion-identity-dim-higher-intro} illustrate that $u$ near infinity is determined 
completely by the function $v$ near the origin.

In this paper, we discuss the converse. 
Can we find all solutions $u$ to \eqref{eq-MSE-intro} near infinity through $v$?
To answer this question, we need to understand the function $v$ more precisely. 
A further study reveals the following property. 
Let $v$ be the function as in \eqref{eq-MSE-expansion-identity-dim-2-intro} 
or \eqref{eq-MSE-expansion-identity-dim-higher-intro}. 
We will prove that there exists a unique pair of a harmonic function $h$
and an analytic function $w$ near $0\in \R^n$ such that 
\begin{align}\label{eq-decomposition-v-MSE-intro}
v(y)=h(y)+|y|^2w(y)\quad\text{near }y=0\in \R^n.\end{align}
Here, we write $y=x/|x|^2$, through the Kelvin transform. 
The functions $h$ and $w$ are determined by $u$ but have different roles. 
In fact, there is no known restriction on the harmonic function $h$. 
However, once $h$ is known, the function $w$ can be determined from the equation \eqref{eq-MSE-intro}. 
In this sense, we say that the asymptotic behavior of $u$ is determined by $h$. 
Now, our existence question can be reformulated as follows. 
Does any harmonic function $h$ near the origin determine a solution of the equation \eqref{eq-MSE-intro} near infinity 
whose asymptotic behavior is determined by $h$?
We can interpret harmonic functions here as ``free data'', the same role as boundary values in the boundary value problems. 

We prove the following result. 

\begin{theorem}\label{thrm-harmonic-fcn-correspondence-MSE-intro}
Let $n \ge 2$, $a$ be a constant, and $h$ be a harmonic function in $B_1$.

$\mathrm{(i)}$ For $n =2$, there exists a unique analytic function $w$ in $B_R$ for some $R \in (0,1)$ such that 
$$u(x)=a\log|x|+h\Big(\frac{x}{|x|^2}\Big) +|x|^{-2}w\Big(\frac{x}{|x|^2}\Big)$$
is a solution to \eqref{eq-MSE-intro} in $\R^n \setminus \overline{B}_{1/R}$.

$\mathrm{(ii)}$ For $n \geq 3$, there exists a unique analytic function $w$ in $B_R$ for some $R \in (0,1)$ such that   
$$u(x)=a+|x|^{2-n}h\Big(\frac{x}{|x|^2}\Big) +|x|^{-n}w\Big(\frac{x}{|x|^2}\Big)$$
is a solution to \eqref{eq-MSE-intro} in $\R^n \setminus \overline{B}_{1/R}$.
\end{theorem}

Note that a bounded harmonic function near infinity 
in dimension $n \ge 3$ can be written as
\begin{align*}
a + |x|^{2-n} h\big( x/|x|^2 \big), 
\end{align*}
for some constant $a \in \R$ and harmonic function $h$ near $0$.
Hence, Theorem \ref{thrm-harmonic-fcn-correspondence-MSE-intro} establishes a
one-to-one and onto correspondence between the space of bounded harmonic functions
near infinity and the space of minimal graphs with horizontal tangent planes
at infinity.

We now turn our attention to the Monge-Amp\`{e}re equation. 
A classical result asserts that any classical convex solution of 
$\det(\nabla^2u)=1$ in $\R^n$
is a quadratic polynomial, due to J\"orgens \cite{Jorgens1954} for $n = 2$, 
Calabi \cite{Calabi1958} for $n \le 5$, and Pogorelov \cite{Pogorelov} for any $n \ge 2$. 
Caffarelli \cite{Caffarelli1995} extended this result
from classical solutions to viscosity solutions. 
Ferrer, Mart\'{i}nez, and Mil\'{a}n \cite{FeMaMi1999} studied convex classical solutions of $\det(\nabla^2u) = 1$ 
outside a bounded subset of $\R^2$. 
Caffarelli and Li \cite{CaLi2003} studied convex viscosity solutions of $\det(\nabla^2u) = 1$ outside a bounded subset of $\R^n$
and proved that solutions are well-approximated by quadratic polynomials at infinity. 
For 2-dimensional case, quadratic polynomials should be modified by an additional logarithmic term. 
To present their results, we introduce 
$$\mathcal{A}=\{A : A \text{ is a real } n \times n \text{ positive definite matrix with }\det(A) = 1\}.$$
Let $u\in C(\mathbb{R}^n\setminus \overline{B}_1)$ be a locally convex viscosity solution of 
$\det(\nabla^2 u)=1$ in $\R^n\setminus \overline{B}_1.$
Caffarelli and Li \cite{CaLi2003} proved that $u \in C^\infty(\mathbb{R}^n\setminus \overline{B}_1)$, and 

$\mathrm{(i)}$ for $n = 2$, there exist some $c, d \in \mathbb{R}$, $b\in \mathbb{R}^2$, and $A \in \mathcal{A}$, such that
\begin{align*}
u(x)=\frac12x^TAx+\langle b,x\rangle+d\log\sqrt {x^TAx}+c+O(|x|^{-1})\quad\text{as } x\to\infty;
\end{align*}

$\mathrm{(ii)}$ for $n \geq 3$, there exist some $c\in\mathbb{R}$, $b \in \mathbb{R}^n$, and $A\in \mathcal{A}$, such that  
\begin{align*}
u(x)=\frac12x^TAx+\langle b,x\rangle+c+O(|x|^{2-n})\quad\text{as } x\to\infty.
\end{align*}

It is possible to follow Caffarelli and Li to approximate convex solutions of the Monge-Amp\`{e}re equation near infinity up 
to arbitrary orders. 
Han and Wang \cite{Han-Wang2024} took a different approach and used one to two functions to 
identify behaviors of solutions near infinity precisely. 
Their result can be summarized as follows. 
We first introduce a notation. For a positive definite matrix $A$ above, we denote by $\sqrt{A}$ a positive square root of $A$; 
namely, $\sqrt{A}$ is a positive definite matrix such that $(\sqrt{A})^2=A$. 
Let $u\in C(\mathbb{R}^n\setminus \overline{B}_1)$ be a locally convex viscosity solution of 
\begin{equation}\label{eq-MA-intro}\det(\nabla^2 u)=1 \quad\text{in } \R^n\setminus \overline{B}_1.\end{equation}
Then, 

$\mathrm{(i)}$ for $n = 2$, there exist some $d\in\mathbb{R}$, $b \in \mathbb{R}^2$, $A\in \mathcal{A}$, 
and a function $v\in C^\infty(B_r)$ for some $r\in (0,1)$, 
such that, for any $x\in  \R^2\setminus \overline{B}_R$  for some $R>0$, 
\begin{equation}\label{eq-MA-expansion-identity-dim-2-intro}u(x)=\frac{1}{2}x^TAx+ \langle b, x\rangle+d\log\sqrt {x^TAx}+
v\Big(\frac{\sqrt{A}x}{|\sqrt{A}x|^2}\Big);\end{equation}

$\mathrm{(ii)}$ for $n \geq 3$, there exist some $c\in\mathbb{R}$, $b \in \mathbb{R}^n$, $A\in \mathcal{A}$, 
and a function $v\in C^\infty(B_r)$ for $n$ even and $v\in C^{n-1,1}(B_r)$ for $n$ odd, for some $r\in (0,1)$,  
such that, for any $x\in  \R^n\setminus \overline{B}_R$ for some $R>0$, 
\begin{equation}\label{eq-MA-expansion-identity-dim-large-intro}u(x)=\frac{1}{2}x^TAx+ \langle b, x\rangle +c+
|\sqrt{A}x|^{2-n}v\Big(\frac{\sqrt{A}x}{|\sqrt{A}x|^2}\Big).\end{equation}
Moreover, for $n$ odd, 
there exist functions $\widetilde{v}, \widehat{v}\in C^\infty(B_1)$ such that  
\begin{equation}\label{eq-expansion-MA-odd-dim-intro} 
v(y)=\widetilde{v}(y)+|y|^n\widehat{v}(y).\end{equation}

We can decompose $v$ in \eqref{eq-MA-expansion-identity-dim-2-intro} 
and \eqref{eq-MA-expansion-identity-dim-large-intro} further, similarly
as for $v$ in \eqref{eq-MSE-expansion-identity-dim-2-intro} 
and \eqref{eq-MSE-expansion-identity-dim-higher-intro}. 
Now, there is a difference between even dimensions and odd dimensions. 
Let $v$ be as in \eqref{eq-MA-expansion-identity-dim-2-intro} 
or \eqref{eq-MA-expansion-identity-dim-large-intro}. If $n$ is even, 
we will prove that there exists a unique pair of a harmonic function $h$ 
and an analytic function $w$ near $0\in \R^n$ such that 
\eqref{eq-decomposition-v-MSE-intro} holds. 
If $n$ is odd, we will prove that
there exist a unique harmonic function $h$ and unique analytic functions $\widetilde{w}$ 
and $\widehat{w}$ near $0\in \R^n$ such that 
\begin{align}\label{eq-decomposition-v-MA-odd-intro}
v(y)=h(y)+|y|^2\widetilde{w}(y)+|y|^{n}\widehat{w}(y)\quad\text{near }y=0\in \R^n.\end{align}
As for the minimal surface equation, $h$ is determined by $u$. 
Once $h$ is known, $w$ or $\widetilde{w}, \widehat{w}$
can be determined from the equation \eqref{eq-MA-intro}. 
In this sense, we say that the asymptotic behavior of $u$ is determined by $h$. 

We prove the following result. 

\begin{theorem}
\label{thrm-harmonic-fcn-correspondence-MA-intro}
Let $n \ge 2$, $b \in \R^n$, $c,d \in \R$, $A$ be an $n \times n$ positive definite matrix with $\det A = 1$, 
and $h$ be a harmonic function in $B_1$.

$\mathrm{(i)}$ For $n =2$, there exists a unique analytic function $w$ in $B_r$ for some $r \in (0,1)$ such that 
$$u(x)=\frac{1}{2}x^T A x+\langle b, x\rangle + d \log\sqrt{x^T A x} 
+ h\Big(\frac{\sqrt{A}x}{|\sqrt{A}x|^2}\Big) +|\sqrt A x|^{-2}w\Big(\frac{\sqrt{A}x}{|\sqrt{A}x|^2}\Big)$$
is a solution to \eqref{eq-MA-intro} in $\R^n \setminus \overline{B}_{R}$ for some $R>0$.

$\mathrm{(ii)}$ For $n \geq 4$ and $n$ even, there exists a unique analytic function $w$ in $B_r$ for some $r \in (0,1)$ such that   
$$u(x)=\frac{1}{2}x^T A x+\langle b, x\rangle + c+|\sqrt Ax|^{2-n}h\Big(\frac{\sqrt A x}{|\sqrt Ax|^2}\Big) 
+|\sqrt A x|^{-n}w\Big(\frac{\sqrt Ax}{|\sqrt Ax|^2}\Big)$$
is a solution to \eqref{eq-MA-intro} in $\R^n \setminus \overline{B}_{R}$ for some $R>0$.

$\mathrm{(iii)}$ For $n \geq 3$ and $n$ odd, there exists a unique pair of analytic functions $\widetilde{w}$ 
and $\widehat{w}$ in $B_r$ for some $r \in (0,1)$ such that
\begin{align*}
u(x)&=\frac{1}{2}x^T A x+\langle b, x\rangle + c+|\sqrt Ax|^{2-n}h\Big(\frac{\sqrt A x}{|\sqrt Ax|^2}\Big)\\
&\qquad\qquad+|\sqrt A x|^{-n}\widetilde{w}\Big(\frac{\sqrt Ax}{|\sqrt Ax|^2}\Big)
+|\sqrt A x|^{2-2n}\widehat{w}\Big(\frac{\sqrt Ax}{|\sqrt Ax|^2}\Big)
\end{align*}
is a solution to \eqref{eq-MA-intro} in $\R^n \setminus \overline{B}_{R}$ for some $R>0$.
\end{theorem}

We briefly describe the proof of Theorem \ref{thrm-harmonic-fcn-correspondence-MSE-intro} and 
Theorem \ref{thrm-harmonic-fcn-correspondence-MA-intro} for $A=I$, the identity matrix. 
We start with the minimal surface equation. 
The Kelvin transform changes the minimal surface equation for $u$ as in \eqref{eq-MSE-intro}
to a quasilinear elliptic equation for $v$ in $B_1$, as in 
\eqref{eq-nonlinear-Poission-MSE-fully} with \eqref{eq-MSE-dim2-fully} and \eqref{eq-MSE-dim3-above-fully}. 
For any given harmonic function $h$ near the origin, we will prove the existence of a solution $v$ 
in the form \eqref{eq-decomposition-v-MSE-intro}. 
To this end, we expand $w$ in \eqref{eq-decomposition-v-MSE-intro} by its Taylor series and 
derive precise estimates of its homogeneous Taylor polynomials. 
The derivation is delicate due to the complicated nonlinear structure of the equation. 
The study of the Monge-Amp\`{e}re equation follows a similar process, only more complicated. 
For even dimensions, $v$ satisfies a fully nonlinear equation. 
For odd dimensions, $\widetilde{v}$ and $\widehat{v}$ in \eqref{eq-expansion-MA-odd-dim-intro} 
satisfy a fully nonlinear elliptic system. 
Fully nonlinear structures add extra layers of difficulties.

On a technical point, the degree of the factor $|y|^2$ in the decomposition of $v$ plays a crucial role in our analysis. 
The decompositions \eqref{eq-decomposition-v-MSE-intro} and \eqref{eq-decomposition-v-MA-odd-intro}
demonstrate that $|y|^2$ distinguishes the harmonic part and non-harmonic part. 
In fact, any harmonic function does not allow a factor of $|y|^2$ unless it is identically zero. 

It is well-known that Dirichlet 
problems in exterior domains
do not always admit solutions.
Refer to \cite{Aiolfi&Bustos2022, BaoLi2012, Osserman1986, Krust1989, KutevTomi1998, Kuwert1993, WangBao2013}
for details.
Theorems \ref{thrm-harmonic-fcn-correspondence-MSE-intro} and
\ref{thrm-harmonic-fcn-correspondence-MA-intro}
demonstrate that prescribing asymptotic behaviors forms a new type of problem
leading to all solutions near infinity. The harmonic functions determining the asymptotic behaviors
play the role of ``free data" as the boundary values in the boundary value problems.

Recently, asymptotic expansions of solutions in exterior domains 
similar to those in \cite{BL, CaLi2003, FeMaMi1999, Han-Wang2024, Sc} 
are established for the special Lagrangian equation in \cite{Han-Marchenko202?} and \cite{LLY2020}, 
for the maximal surface equation in \cite{HongYuan2020}, 
and for the Lagrangian mean curvature equation in \cite{BaoLiu2020}. 
We expect results in Theorem \ref{thrm-harmonic-fcn-correspondence-MSE-intro} and 
Theorem \ref{thrm-harmonic-fcn-correspondence-MA-intro} hold for these equations.

The method in this paper can be generalized to other geometric equations 
such as the special Lagrangian equation. 
Due to the more complicated structure of the special Lagrangian equation, 
several results in this paper for the Monge-Ampere equation need to be 
improved. We will not pursue along this line in this paper.

The paper is organized as follows. 
In Section \ref{sec-polynomials}, we derive some basic estimates of homogeneous polynomials. 
In Section \ref{sec-decompositions}, we discuss decompositions of analytic functions and 
introduce the notion of purely unharmonic functions.  
In Section \ref{sec-linear-eqs}, we study linear equations and prove the existence of 
purely unharmonic solutions under additional structure conditions. 
In section \ref{sec-model-fully-nonlinear}, we study two simplified model equations.
In Section \ref{sec-fully-nonlinear}, we study a general class of nonlinear equations 
and prove the existence of purely unharmonic solutions. 
In Sections \ref{sec-MSE} and \ref{sec-MA}, we study the minimal surface equation and 
the Monge-Amp\`{e}re equation and prove 
Theorem \ref{thrm-harmonic-fcn-correspondence-MSE-intro} and 
Theorem \ref{thrm-harmonic-fcn-correspondence-MA-intro}, respectively.
In section \ref{sec-optimal-decomposition-MA-odd},
we prove the decomposition \eqref{eq-decomposition-v-MA-odd-intro} for solutions
of the Monge-Amp\`{e}re equation in odd dimensions.
In section \ref{sec-conti-depend}, we study the dependence of the solutions on harmonic functions.

\section{Estimates of Polynomials}\label{sec-polynomials}

Throughout the paper, we fix $n\geq 2$ as the dimension of $\mathbb R^n$. 
For simplicity, all functions are assumed to be smooth in a neighborhood of the origin, 
and our focus will be on analytic functions. 
In this section, we study polynomials. 

We first introduce two notations. 
For each nonnegative integer $k$, define 
$$\mathscr P_k=\{\text{homogeneous polynomials of degree $k$ in $\mathbb R^n$}\},$$
and 
$$\mathscr H_k=\{\text{homogeneous harmonic polynomials of degree $k$ in $\mathbb R^n$}\}.$$
We always identify a homogeneous polynomial $P_k \in \mathscr{P}_k$ 
with its restriction $P_k|_{\mathbb{S}^{n-1}}$ to the unit sphere $\mathbb S^{n-1}$. 
It is obvious that $P_k$, as a homogeneous polynomial of degree $k$, 
is uniquely determined by its restriction $P_k|_{\mathbb{S}^{n-1}}$ to the unit sphere.

It is well-known that 
$$\mathscr P_k=\mathscr H_k\oplus \mathscr H_k^\perp,$$
where $\perp$ is in the $L^2$ sense. 
We emphasize that $\mathscr{H}_k^\perp$ is the orthogonal complement of $\mathscr{H}_k$ in $\mathscr{P}_k$.
Here, we can integrate either the involved functions over the unit ball 
or their restrictions over the unit sphere. 
For any $P_k\in \mathscr P_k$, there exists a unique decomposition 
$P_k=\widetilde{P}_k+\widehat{P}_k$ for some $\widetilde{P}_k\in \mathscr H_k$ and $\widehat{P}_k\in \mathscr H_k^\perp$.

Note that 
$$\mathscr P_0=\mathscr H_0=\{\text{constants}\},\quad \mathscr P_1=\mathscr H_1=\{\text{linear functions}\}.$$
Then,  $\dim \mathscr P_0=\dim \mathscr H_0=1$ 
and $\dim \mathscr P_1=\dim \mathscr H_1=n$. 
It is well-known that 
\begin{align*}
\dim \mathscr P_k=\left(\begin{matrix} n+k-1\\k\end{matrix}\right)\quad\text{for }k\ge 1,
\end{align*}
and 
\begin{align*}
\dim \mathscr H_k=\left(\begin{matrix} n+k-1\\k\end{matrix}\right)
-\left(\begin{matrix} n+k-3\\k-2\end{matrix}\right)\quad\text{for }k\ge 2. 
\end{align*}
Hence, 
\begin{align*}
\dim \mathscr H_k^\perp=\dim \mathscr P_{k-2}\quad\text{for }k\ge 2.\end{align*}
It is straightforward to verify that, for any $k\ge 2$, 
\begin{align*}
c_1k^{n-1}\le\dim \mathscr P_k\le c_2k^{n-1},\quad
c_1k^{n-2}\le\dim \mathscr H_k\le c_2k^{n-2},
\end{align*}
for some positive constants $c_1$ and $c_2$ depending only on $n$. 

We now discuss further decompositions of homogeneous polynomials. 
For each $k\ge 2$, the space $\mathscr P_k$ has the following orthogonal decomposition: 
$$\mathscr P_k=\mathscr H_k\oplus |x|^2\mathscr P_{k-2}.$$
This gives a characterization $\mathscr H_k^\perp=|x|^2\mathscr P_{k-2}$. 
In other words, for any $P_k\in \mathscr P_k$, 
there exists a unique pair $p_k\in \mathscr H_k$ and $P_{k-2}\in  \mathscr P_{k-2}$ such that 
\begin{align}\label{eq-decomposition-poly}  
P_k=p_k+|x|^2P_{k-2}.\end{align}
Upon an iteration, we obtain 
$$\mathscr P_k = \bigoplus_{j=0}^{[k/2]} |x|^{2j} \mathscr H_{k-2j}.$$
Specifically, we have
$$\mathscr P_k=\mathscr H_k\oplus |x|^2\mathscr H_{k-2}\oplus |x|^4\mathscr H_{k-4}
\oplus\cdots\oplus |x|^k\mathscr H_{0}\quad\text{for $k$ even},$$
and 
$$\mathscr P_k=\mathscr H_k\oplus |x|^2\mathscr H_{k-2}\oplus |x|^4\mathscr H_{k-4}
\oplus\cdots\oplus |x|^{k-1}\mathscr H_{1}\quad\text{for $k$ odd}.$$
Then, for any $P_k\in \mathscr P_k$, there exists a unique decomposition given by 
$$P_k=\sum_{j=0}^{[k/2]} |x|^{2j} p_{k-2j},$$
for some $p_{k-2j}\in \mathscr H_{k-2j}$ for $j = 0,\cdots,[k/2]$. 
We also have the following decomposition:  
$$\mathscr H_k^\perp=\bigoplus_{j=1}^{[k/2]} |x|^{2j}\mathscr H_{k-2j}.$$
Refer to Chapter 5 in \cite{HarmonicFcnTheory} for the decomposition of homogeneous polynomials.

The decomposition in \eqref{eq-decomposition-poly} demonstrates that 
any homogeneous polynomial  can be decomposed into a sum of 
a homogeneous harmonic polynomial of the same degree and the $|x|^2$-multiple of 
a homogeneous polynomial of degree 2 less. 
We will discuss to what extend such a decomposition holds for arbitrary functions, not necessarily 
homogeneous polynomials. We first prove a preliminary result. 

\begin{lemma}\label{lemma-Laplace-inj-sur} 
For $k\ge 2$, 
$\Delta: \mathscr H_{k}^\perp\to \mathscr P_{k-2}$ is injective and surjective.
\end{lemma} 

\begin{proof} 
Since $\mathscr H_{k}^\perp$ and $\mathscr P_{k-2}$ 
have the same dimension, $\Delta: \mathscr H_{k}^\perp\to \mathscr P_{k-2}$ is injective 
if and only if it is surjective.  
For the injectiveness, we only need to prove that 
if $\Delta P_k=0$ for some $P_k\in \mathscr H_{k}^\perp$, then $P_k=0$. 
This is clear from the definitions of $\mathscr H_{k}$ and $\mathscr H_{k}^\perp$. 
\end{proof} 

We now provide another point of view of Lemma \ref{lemma-Laplace-inj-sur}. 
We first record several useful formulas, which will be used repeatedly. 
For any $m\ge 0$ and $u\in C^2(B_1)$, a simple computation yields 
$$\Delta(|x|^{m+2} u)=|x|^{m+2}\Delta u+2(m+2)|x|^{m}\langle x, \nabla u\rangle+(m+2)(n+m)|x|^{m}u,$$
where $\langle \cdot, \cdot\rangle$ is the inner product in $\mathbb R^n$. 
In particular, for $P_k\in \mathscr P_k$, we have 
\begin{align}\label{eq-Laplace-|x|-power}
\Delta(|x|^{m+2} P_k)=|x|^{m+2}\Delta P_k+(m+2)(2k+n+m)|x|^{m}P_k.
\end{align}
The following result plays a fundamental role in our subsequent discussions. 

\begin{lemma}\label{lemma-Laplace-estimate-poly-L2} 
For any $m\ge 0$, $k\ge 0$, and $Q_{k}\in \mathscr P_{k}$, there exists a unique $P_k\in \mathscr P_k$ 
satisfying 
\begin{align*}
\Delta (|x|^{m+2}P_k)=|x|^mQ_{k}.\end{align*}  
Moreover, $P_k$ satisfies 
\begin{align}\label{eq-estimate-L2-shpere}
\|P_k\|_{L^2(\mathbb{S}^{n-1})}
\le \big[(m+2)(m+2k+n)\big]^{-1} \|Q_{k}\|_{L^2(\mathbb{S}^{n-1})}.\end{align}
\end{lemma} 

\begin{proof}
The existence and uniqueness of $P_k$ 
is equivalent to the statement that  
$$\Delta: |x|^{m+2}\mathscr P_{k}\to |x|^{m}\mathscr P_{k}\text{ is injective and surjective}.$$ 
The case $m=0$ was already discussed in Lemma \ref{lemma-Laplace-inj-sur}. 
For general $m \ge 1$, it is straightforward to prove the injectiveness. 
Indeed, assume $\Delta(|x|^{m+2} P_k) = 0$ for some $P_k \in \mathscr P_k$, 
if $m$ is even, then we can apply Lemma \ref{lemma-Laplace-inj-sur} to $|x|^m P_k$; 
if $m$ is odd, then $|x|^{m+2}P_k$ is smooth near $0$ by the regularity of harmonic functions.
In either case, we have $P_k = 0$.
We now present a proof of the surjectiveness by deriving an explicit expression of the solution $P_k$ 
in terms of $Q_{k}$. Such an expression is important in the derivation of the desired estimate.

Let $Q_{k}\in \mathscr P_{k}$ be given. Then, there exist unique $q_{k-2j}\in \mathscr H_{k-2j}$ for $0 \le j \le [k/2]$ such that 
$$Q_{k}=\sum_{j=0}^{[k/2]} |x|^{2j} q_{k-2j}.$$
For $P_k\in \mathscr P_k$, we will find $p_{k-2j}\in \mathscr H_{k-2j}$ for $0 \le j \le [k/2]$ such that $P_k$ defined by   
$$P_k=\sum_{j=0}^{[k/2]} |x|^{2j} p_{k-2j}$$
satisfies $\Delta (|x|^{m+2}P_k)=|x|^mQ_{k}$. 
To this end, we require, for $0 \le j \le [k/2]$,  
\begin{align}\label{eq-Poisson-special-p-q}\Delta(|x|^{m+2j+2}p_{k-2j})=|x|^{m+2j}q_{k-2j}.\end{align}
By \eqref{eq-Laplace-|x|-power} and $\Delta p_{k-2j}=0$, we have 
$$\Delta(|x|^{m+2j+2}p_{k-2j})=(m+2j+2)(m+2k-2j+n)|x|^{m+2j}p_{k-2j}.$$
To have \eqref{eq-Poisson-special-p-q}, we simply take
\begin{align}\label{eq-p_i-q_i}
p_{k-2j}=\big[(m+2j+2)(m+2k-2j+n)\big]^{-1}q_{k-2j}.\end{align}
This finishes the proof of the surjectiveness and, in fact, 
establishes an explicit expression of $P_k$ in terms of $Q_{k}$. 

For $n\ge 2$ and $0\le j \le [k/2]$, we have
\begin{align}\label{eq-estimate-coeff-min}\begin{split}
&(m+2j+2)(m+2k-2j+n) \\
&\quad= (m+2)(m+2k+n) + 2j\big[ (m+2k + n) - (m+2) - 2j \big]\\
&\quad = (m+2)(m+ 2k +n) + 2j(n-2 + 2k-2j).\\
&\quad\geq (m+2)(m+2k+n).
\end{split}\end{align}
Hence, for $0 \le j \le [k/2]$,  
$$|p_{k-2j}|\le \big[(m+2)(m+2k+n)\big]^{-1}|q_{k-2j}|\quad\text{on }\R^n.$$
Now, we have
\begin{align*}P_k=\sum_{j=0}^{[k/2]} p_{k-2j}
\quad\text{on }\mathbb{S}^{n-1},\end{align*}
and 
\begin{align*}Q_{k}=\sum_{j=0}^{[k/2]} q_{k-2j}
\quad\text{on }\mathbb{S}^{n-1}.\end{align*}
It is clear that 
\begin{align*}
\int_{\mathbb S^{n-1}}P_k^2d\sigma
=\int_{\mathbb S^{n-1}}\sum_{j=0}^{[k/2]} p_{k-2j}^2 d\sigma, \end{align*}
and 
\begin{align*}\int_{\mathbb S^{n-1}}Q_{k}^2d\sigma
=\int_{\mathbb S^{n-1}}\sum_{j=0}^{[k/2]} q_{k-2j}^2d\sigma. 
\end{align*}
By \eqref{eq-p_i-q_i} and \eqref{eq-estimate-coeff-min}, we obtain 
\begin{align*}
\int_{\mathbb S^{n-1}}P_k^2d\sigma
\le  \big[(m+2)(m+2k+n)\big]^{-2}\int_{\mathbb S^{n-1}}Q_{k}^2d\sigma. 
\end{align*}
This is \eqref{eq-estimate-L2-shpere}. 
\end{proof} 

We next reformulate Lemma \ref{lemma-Laplace-estimate-poly-L2}. 
For any $Q \in \mathscr{P}_{k}$, we examine 
\begin{align}\label{eq-Poisson-degree-k-case2}\Delta P=Q.\end{align}
By Lemma \ref{lemma-Laplace-estimate-poly-L2}, 
there exists a unique solution $P\in \mathscr H_{k+2}^\perp$ of \eqref{eq-Poisson-degree-k-case2}. 
In the following, we write 
\begin{align*}
P = \Delta_{k+2}^{-1} Q.\end{align*} 
Lemma \ref{lemma-Laplace-estimate-poly-L2} also yields
\begin{align}\label{eq-estimate-poly-special-norm-0-deg-|x|^2-case2}
\|\Delta_{k+2}^{-1} Q\|_{L^2(\mathbb S^{n-1})}
\le \big[2(2k+n)\big]^{-1} \|Q\|_{L^2(\mathbb S^{n-1})}.\end{align}

\begin{lemma}\label{lemma-estimate-degree-|x|2-case2}
Assume 
$Q\in |x|^{2\ell}\mathscr{P}_{k-2\ell}$ with $k\ge 2\ell$. 
Then, $\Delta_{k+2}^{-1} Q\in |x|^{2\ell+2}\mathscr{P}_{k-2\ell}$ and
\begin{align}
\label{eq-estimate-poly-special-norm-deg-|x|^2-case2}
\|\Delta_{k+2}^{-1} Q\|_{L^2(\mathbb S^{n-1})} \leq \big[ (2\ell + 2) (2k-2\ell  + n) \big]^{-1} \|Q\|_{L^2(\mathbb S^{n-1})}.
\end{align}
\end{lemma}

\begin{proof} Write $Q=|x|^{2\ell}\widetilde Q$ for some $\widetilde Q \in \mathscr{P}_{k-2\ell}$.
By Lemma \ref{lemma-Laplace-estimate-poly-L2}, 
there exists a unique $\widetilde P\in \mathscr{P}_{k-2\ell}$ such that 
$\Delta(|x|^{2\ell+2}\widetilde P)=|x|^{2\ell}\widetilde Q$. The uniqueness of $\Delta_{k+2}^{-1} Q$ yields 
$\Delta_{k+2}^{-1} Q=|x|^{2\ell+2}\widetilde P\in |x|^{2\ell+2}\mathscr{P}_{k-2\ell}$.  
With $2 \ell$ and $k-2\ell$ replacing $m$ and $k$ in 
Lemma \ref{lemma-Laplace-estimate-poly-L2}, respectively, we have 
\eqref{eq-estimate-poly-special-norm-deg-|x|^2-case2}. 
\end{proof} 

We point out that \eqref{eq-estimate-poly-special-norm-deg-|x|^2-case2} 
improves \eqref{eq-estimate-poly-special-norm-0-deg-|x|^2-case2}. 
In fact, consider 
$$h(s)=(2s+2)(2k-2s+n).$$
It is easy to verify that $h$ is increasing on $[0, k/2]$ and in particular $h(s)\ge h(0)$ 
for any $s$ with $0\le 2s\le k$. 

Lemma \ref{lemma-estimate-degree-|x|2-case2} plays an essential role in the proof of 
many results later on. It demonstrates that an additional factor of $|x|^2$ in $Q \in \mathscr{P}_{k}$ 
improves the estimate of $\Delta_{k+2}^{-1}$. 

We next prove a well-known estimate for derivatives of homogeneous polynomials. 

\begin{lemma}\label{lemma-homo-poly-first-derivatives-L2} 
Assume $P_k\in\mathscr P_k$ for some $k \ge 0$. Then, for $i = 1,\cdots,n$, 
\begin{align}\label{eq-homo-poly-first-derivatives-L2}
\|\partial_iP_k\|_{L^2(\mathbb{S}^{n-1})} \leq C k \|P_k\|_{L^2(\mathbb{S}^{n-1})},
\end{align}
where $C$ is a positive constant depending only on $n$.
\end{lemma}

\begin{proof} If $k = 0$, then both sides of \eqref{eq-homo-poly-first-derivatives-L2} are $0$. 
Assume that $k \geq 1$.
The proof consists of two steps.

{\it Step 1.} We first consider $P_k\in\mathscr H_k$. 
A simple integration by parts yields
$$\int_{B_1} |\nabla P_k|^2 dx = \int_{\mathbb{S}^{n-1}} P_k \langle \nabla P_k, \theta \rangle d\sigma
= k \|P_k\|_{L^2 (\mathbb{S}^{n-1})}^2.$$
We point out that $\nabla P_k$ here is the gradient of $P_k$ with respect to the Euclidean coordinates. 
Now, $|\nabla P_k|^2$ is a homogeneous polynomial of degree $2k - 2$. Then,
$$\int_{B_1} |\nabla P_k|^2 dx = \frac{1}{n + 2k - 2} \int_{\mathbb{S}^{n-1}} |\nabla P_k|^2 d\sigma,$$
and hence 
$$\|\nabla P_k\|_{L^2(\mathbb S^{n-1})} = \sqrt{k (n + 2k - 2)}  \|P_k\|_{L^2(\mathbb S^{n-1})}.$$
In particular, we have 
\begin{align}
\label{eq-polynomial-derivative-estimate-step1}
\|\nabla P_k\|_{L^2(\mathbb S^{n-1})} \leq C k \|P_k\|_{L^2(\mathbb S^{n-1})},    
\end{align}
where $C$ is a positive constant depending only on $n$. 

{\it Step 2.} We now consider $P_k\in\mathscr P_k$. 
Write
\begin{align}
\label{eq-harmonic-decomposition}
P_k(x) = \sum_{j=0}^{[ k/2]} |x|^{2j} p_{k-2j},
\end{align}
where $p_{k-2j}\in\mathscr{H}_{k-2j}$ for $j=0, \cdots, [ k/2]$. 
Note that \eqref{eq-harmonic-decomposition} is an orthogonal decomposition in $L^2 (\mathbb{S}^{n-1}).$ 
It is obvious that 
$$\partial_i P_k=\sum_{j=0}^{[ k/2 ]} \big(2j |x|^{2j-2} x_i p_{k-2j} + |x|^{2j}\partial_i p_{k-2j}\big).$$
Hence,
\begin{align*}
\|\partial_i P_k\|_{L^2 (\mathbb{S}^{n-1})} 
&\leq 2\Big\| \sum_{j=0}^{[ k/2 ]} j  p_{k-2j} \Big\|_{L^2 (\mathbb{S}^{n-1})} 
+\Big\| \sum_{j=0}^{[ k/2 ]} \partial_i p_{k-2j}\Big\|_{L^2 (\mathbb{S}^{n-1})}\\
& = 2\Big(\sum_{j=0}^{[ k/2 ]} j^2\|  p_{k-2j} \|_{L^2 (\mathbb{S}^{n-1})}^2 \Big)^{1/2} 
+ \Big(\sum_{j=0}^{[ k/2 ]}\| \partial_ip_{k-2j} \|_{L^2 (\mathbb{S}^{n-1})}^2 \Big)^{1/2}\\
& \leq C k \|P_k\|_{L^2 (\mathbb{S}^{n-1})},
\end{align*}
where we used \eqref{eq-polynomial-derivative-estimate-step1} in the last inequality.
This yields the desired estimate. 
\end{proof}

\begin{lemma}\label{lemma-homo-poly-higher-derivatives} 
Assume $P_k\in\mathscr P_k$ for some $k \ge 1$. 
Then, for any  $\gamma \in \Z_+^n$ with $|\gamma| \leq k$, 
$$\frac{1}{\gamma!}|\partial^\gamma P_k|_{L^\infty(B_1)} \leq C^k \|P_k\|_{L^2(\mathbb{S}^{n-1})},$$
where $C$ is a positive constant depending only on $n$.
\end{lemma}

\begin{proof} We first claim, for any $\gamma \in \Z_+^n$ with $|\gamma|\leq k$, 
\begin{align}
\label{eq-polynomial-derivative-estimate-step3-claim}
\|\partial^\gamma P_k\|_{L^2 (\mathbb{S}^{n-1})} \leq \frac{k!}{(k - |\gamma|)!} 
C^{|\gamma|} \|P_k\|_{L^2 (\mathbb{S}^{n-1})},
\end{align}
where $C$ is the constant in Lemma \ref{lemma-homo-poly-first-derivatives-L2}.
We will prove \eqref{eq-polynomial-derivative-estimate-step3-claim} by induction. 
By Lemma \ref{lemma-homo-poly-first-derivatives-L2}, 
we have \eqref{eq-polynomial-derivative-estimate-step3-claim} for $|\gamma| = 1$. 
Assume that \eqref{eq-polynomial-derivative-estimate-step3-claim} holds for 
some $\gamma \in \Z_+^n$ with $|\gamma| \leq k-1$. 
Applying Lemma \ref{lemma-homo-poly-first-derivatives-L2} to $\partial^\gamma P_k$, we get
$$\|\partial_i \partial^\gamma P_k\|_{L^2 (\mathbb{S}^{n-1})} 
\leq C (k - |\gamma|) \|\partial^{\gamma} P_k\|_{L^2 (\mathbb{S}^{n-1})} 
\leq C^{|\gamma|+1} \frac{k!}{(k - |\gamma| - 1)!} \|P_k\|_{L^2 (\mathbb{S}^{n-1})}.$$
This proves \eqref{eq-polynomial-derivative-estimate-step3-claim}. 

Next, we write
$$P_k(x) = \sum_{|\alpha| = k} c_\alpha x^{\alpha}.$$
Note that $\partial^\alpha P_k$ is constant if $|\alpha| = k$. 
Hence, with $\gamma = \alpha$ in \eqref{eq-polynomial-derivative-estimate-step3-claim}, we get
\begin{align*}
\mathrm{vol}(\mathbb{S}^{n-1})^{1/2}\sum_{|\alpha| = k} |c_\alpha|
&= \sum_{|\alpha| = k} \frac{1}{\alpha!} \|\partial^\alpha P_k\|_{L^2 (\mathbb{S}^{n-1})}\\
&\leq \|P_k\|_{L^2 (\mathbb{S}^{n-1})} C^k\sum_{|\alpha| = k} \frac{k!}{\alpha!} 
= \|P_k\|_{L^2 (\mathbb{S}^{n-1})} (nC)^k,
\end{align*}
where we used the identity
$$\sum_{|\alpha| = k} \frac{k!}{\alpha!} = n^k.$$
By renaming $C$, we obtain 
$$\sum_{|\alpha| = k} |c_\alpha| \leq C^k \|P_k\|_{L^2 (\mathbb{S}^{n-1})},$$
and, in particular, 
$$|c_\alpha| \leq C^k \|P_k\|_{L^2 (\mathbb{S}^{n-1})}.$$
Hence, for any $\gamma \in \Z_+^n$ with $|\gamma| \leq k$,
\begin{align*}
\frac{1}{\gamma!}|\partial^\gamma P_k|_{L^\infty(B_1)} 
&\leq \sum_{|\alpha| = k, \alpha \geq \gamma} \binom{\alpha}{\gamma} |c_\alpha| 
\leq C^k \|P_k\|_{L^2 (\mathbb{S}^{n-1})} \binom{k}{|\gamma|} \sum_{|\alpha| = k} 1\\
&\leq C^k \|P_k\|_{L^2 (\mathbb{S}^{n-1})}\bigg( \sum_{i=0}^k \binom{k}{i}\bigg) 
\bigg( \sum_{|\alpha| = k} \frac{k!}{\alpha!}\bigg)\\
&= (2nC)^k \|P_k\|_{L^2 (\mathbb{S}^{n-1})}.
\end{align*}
This is the desired estimate.
\end{proof}

As a consequence, we obtain the following characterization of analytic functions near $0$.

\begin{cor}\label{cor-constructing-analytic-functions}
Let $\{P_k\}_{k\geq 1}$ be a sequence of homogeneous polynomials in $\R^n$ 
with $ P_k \in \mathscr{P}_k$. Assume, for any $k \ge 1$, 
\begin{align}\label{eq-analyticity-estimate-general}\|P_k\|_{L^2(\mathbb{S}^{n-1})} \leq A^k,\end{align}
for some positive constant $A$.
Then, the function given by
$$u = \sum_{k=1}^\infty P_k$$
is analytic in $B_R \subseteq \R^n$, 
for some positive constant $R$ depending only on $n$ and $A$, and, for any $k \geq 1$, 
$P_k$ is the $k$-th homogeneous Taylor polynomial of $u$ at $0$.
\end{cor}

\begin{proof}
For each $k \geq 1$, write 
$$P_k(x) = \sum_{|\alpha |= k} a_\alpha x^\alpha,$$
for some constants $a_\alpha$ and $\alpha \in \Z_+^n$ with $|\alpha| = k$.
Then, for any $\alpha \in \Z_+^n$ with $|\alpha| = k$,
$$a_\alpha = \frac{1}{\alpha!} \partial^\alpha P_k(0).$$
By Lemma \ref{lemma-homo-poly-higher-derivatives}, we get, for any $\alpha \in \Z_+^n$ with $|\alpha| = k$,
$$|a_\alpha| = \frac{1}{\alpha!} |\partial^\alpha P_k(0)| \leq C^k\|P_k\|_{L^2(\mathbb{S}^{n-1})} \leq (AC)^k,$$
for some positive constant $C$ depending only on $n$.
Therefore, the series
$$u(x) = \sum_{\alpha\neq 0} a_\alpha x^\alpha$$
converges uniformly and absolutely in $B_R$,
for some positive constant $R$ depending only on $n$ and $A.$
This finishes the proof.
\end{proof}

The estimate \eqref{eq-analyticity-estimate-general} is referred to as an {\it analyticity estimate}. 

\smallskip 

We have derived various estimates of polynomials in the $L^2$-norms. 
These estimates are sufficient for linear equations with constant coefficients but 
not adequate for linear equations with variable coefficients or nonlinear equations. 
To estimate products or, more generally, compositions involving solutions, 
we need to have Banach algebra type estimates. To this end, we need to have an estimate 
of derivatives of orders related to the Sobolev embedding.

Fix $m \ge 1$. 
Consider an $H^{m}$-inner product 
on $\mathbb S^{n-1}$ given by, for any $u,v \in H^{m}(\mathbb{S}^{n-1})$, 
\begin{align}
\label{eq-H^m-inner-product-with-constant}
\langle u,v \rangle_{H^{m}(\mathbb{S}^{n-1})} 
=  C(m,n) \sum_{i=0}^{m} \langle \nabla_\theta^i u, \nabla_\theta^i v\rangle_{L^2(\mathbb{S}^{n-1})},
\end{align}
for some positive constant $C(m,n)$ to be chosen. 
It is known that $\Delta_{\theta}$ is self-adjoint in $H^{m}(\mathbb{S}^{n-1})$. 
More specifically, for any $u,v \in C^\infty(\mathbb{S}^{n-1})$, 
$$\langle \Delta_{\theta} u, v\rangle_{H^{m}(\mathbb{S}^{n-1})} = \langle  u,  \Delta_{\theta} v\rangle_{H^{m}(\mathbb{S}^{n-1})}.$$
In fact, for any $u,v \in C^\infty(\mathbb{S}^{n-1})$,
\begin{align}
\label{eq-formula-H^m-2-point-homogeneous}
\langle u,v\rangle_{H^m(\mathbb S^{n-1})} = C(m,n) \langle p(\Delta_\theta) u, v\rangle_{L^2(\mathbb S^{n-1})},
\end{align}
for some polynomial $p$.
For example, $p(t) = 1-t$ for $m = 1$ and
$p(t) = t^2 + (n-3)t + 1$ for $m = 2$.
Refer to Theorem 11\cite{Helgason1959}.
We should point out that the Laplace-Beltrami operator may not be self-adjoint on
$C^\infty(M)$ with respect to $\langle\cdot,\cdot\rangle_{H^m(M)}$, $m \ge 2$,
for a general Riemannian manifold $(M,g)$.

By the Sobolev embedding, if $m > \frac{n-1}{2}$,
we have the Banach algebra property. 
We now choose the constant $C(m,n) > 1$ such that, for any $u,v \in H^{m}(\mathbb{S}^{n-1})$, 
\begin{align}
\label{eq-Banach-algebra-property-without-constant}
\|u v\|_{H^{m}(\mathbb{S}^{n-1})} \leq \|u\|_{H^{m}(\mathbb{S}^{n-1})}
\|v\|_{H^{m}(\mathbb{S}^{n-1})}.
\end{align}
By \eqref{eq-formula-H^m-2-point-homogeneous}, we have 
$\langle H_i,u\rangle_{H^{m}(\mathbb{S}^{n-1})} = 0$ if and only if $\langle H_i,u\rangle_{L^2(\mathbb S^{n-1})} = 0$, 
for any $H_i \in \mathscr H_i$
and $u \in C^\infty(\mathbb S^{n-1})$. 
Throughout this paper, we fix an integer
$$n_* > \frac{n-1}{2}.$$
For brevity, we write, for any $u,v \in H^{n_*}(\mathbb{S}^{n-1})$,  
\begin{align}
\label{eq-H^*1-inner-product}
\langle u, v\rangle_* = \langle u,v\rangle_{H^{n_*}(\mathbb S^{n-1})},\quad
\|u\|_{*}=\|u\|_{H^{n_*}(\mathbb{S}^{n-1})}.
\end{align}

We point out that $\|\cdot\|_*$ is a norm on $\mathbb{S}^{n-1}$. 
Hence, if $f$ and $g$ are smooth functions on $\R^n$ with $f|_{\mathbb{S}^{n-1}} = g|_{\mathbb{S}^{n-1}}$, 
then $\|f\|_* = \|g\|_*$. In particular, we have, for any $m \ge 0$ and any smooth function $u$ on $\R^n$,
$$\|u\|_* = \big\||x|^m u\big\|_*.$$

The following result is the counterpart of Lemma \ref{lemma-Laplace-estimate-poly-L2} 
for the $\|\cdot\|_*$-norm. 

\begin{lemma}\label{lemma-Laplace-estimate-poly-H*} 
For any $m\ge 0$, $k\ge 0$, and $Q_{k}\in \mathscr P_{k}$, there exists a unique $P_k\in \mathscr P_k$ 
satisfying 
\begin{align*}
\Delta (|x|^{m+2}P_k)=|x|^mQ_{k}.
\end{align*}  
Moreover, $P_k$ satisfies 
\begin{align*}
\|P_k\|_{*}
\le \big[(m+2)(m+2k+n)\big]^{-1} \|Q_{k}\|_{*}.\end{align*}
\end{lemma} 

The proof is similar as that of Lemma \ref{lemma-Laplace-estimate-poly-L2}, 
with the $\|\cdot\|_*$-norm replacing the $L^2$-norm. 
The following result is the counterpart of Lemma \ref{lemma-estimate-degree-|x|2-case2}
for the $\|\cdot\|_*$-norm. 

\begin{lemma}\label{lemma-estimate-degree-|x|2-case2-H*}
Assume 
$Q\in |x|^{2\ell}\mathscr{P}_{k-2\ell}$ with $k\ge 2\ell$. 
Then, $\Delta_{k+2}^{-1} Q\in |x|^{2\ell+2}\mathscr{P}_{k-2\ell}$ and
\begin{align*}
\|\Delta_{k+2}^{-1} Q\|_{*} \leq \big[ (2\ell + 2) (2k-2\ell  + n) \big]^{-1} \|Q\|_{*}.
\end{align*}
\end{lemma}

The following result is the counterpart of Lemma \ref{lemma-homo-poly-first-derivatives-L2}
for the $\|\cdot\|_*$-norm. 

\begin{lemma}\label{lemma-homo-poly-first-derivatives-H*} 
Assume $P_k\in\mathscr P_k$ for some $k \ge 0$. Then, for $i = 1,\cdots,n$, 
\begin{align}\label{eq-homo-poly-first-derivatives-H*}
\|\partial_iP_k\|_{*} \leq C k \|P_k\|_{*},
\end{align}
where $C$ is a positive constant depending only on $n$.
\end{lemma}

\begin{proof} If $k = 0$, then both sides of \eqref{eq-homo-poly-first-derivatives-H*} are $0$. Assume that $k \geq 1$.

We first consider $P_k\in\mathscr H_k$. 
By \eqref{eq-homo-poly-first-derivatives-L2}, we have \begin{align*}
\|\nabla P_k\|_{L^2(\mathbb{S}^{n-1})} \leq C k \|P_k\|_{L^2(\mathbb{S}^{n-1})},
\end{align*}
where $C$ is a positive constant depending only on $n$.
Let $\lambda_k = k(k+n-2)$ be the $k$-th eigenvalue of $-\Delta_\theta$. 
Here, we arrange eigenvalues in an increasing order without counting multiplicity.
Then, $P_k|_{\mathbb{S}^{n-1}}\in \mathscr{E}_k$ and $\partial_iP_k|_{\mathbb{S}^{n-1}}\in \mathscr{E}_{k-1}$, 
the corresponding eigenspaces. 
Hence, for any $1 \leq i \leq n$,
\begin{equation}
\label{eq-estimate-first-derivative-of-homo-harmonic-poly-H*}
\begin{aligned}
\| \partial_i P_k\|_* 
&\le C(1 + \lambda_{k-1})^{\frac{n_*}{2}}\|\partial_i P_k\|_{L^2(\mathbb{S}^{n-1})}\\
& \leq C  k(1 + \lambda_{k})^{\frac{n_*}{2}}\|P_k\|_{L^2(\mathbb{S}^{n-1})} \le Ck \|P_k\|_*,
\end{aligned}
\end{equation}
where $C$ is a constant (different from line to line) 
depending only on $n$. 

Now, consider a general $P_k \in \mathscr{P}_k$. Write
$$P_k = \sum_{j=0}^{[ k/2]} |x|^{2j} p_{k-2j},$$
where $p_{k-2j} \in \mathscr{H}_{k - 2j}$ for $j = 0,\cdots, [ k/2]$.
As in the proof of Lemma \ref{lemma-homo-poly-first-derivatives-L2}, we have
\begin{align*}
\partial_i P_k 
= \sum_{j=0}^{[ k/2]} (2j |x|^{2j-2} x_i p_{k-2j} + |x|^{2j} \partial_i p_{k-2j}),
\end{align*}
and hence
\begin{align*}
\|\partial_i P_k\|_* 
&\leq  2 \|x_i\|_* \Big \| \sum_{j=0}^{[ k/2]}  jp_{k-2j}\Big\|_* +  \Big\|  \sum_{j=0}^{[ k/2]} \partial_i p_{k-2j}\Big\|_*\\
&= 2\|x_i\|_* \Big(  \sum_{j=0}^{[ k/2]}j^2\| p_{k-2j}\|_*^2\Big)^{\frac{1}{2}} 
+ \Big( \sum_{j=0}^{[ k/2]} \|  \partial_i p_{k-2j}\|_*^2 \Big)^{\frac{1}{2}}\\
& \leq (\|x_i\|_* + C) k \|P_k\|_*,
\end{align*}
where we used \eqref{eq-estimate-first-derivative-of-homo-harmonic-poly-H*} in the last inequality.
This proves the desired estimate.
\end{proof}

\section{Decompositions of Analytic Functions}\label{sec-decompositions} 

In this section, we discuss decompositions of analytic functions. 
The decomposition in \eqref{eq-decomposition-poly} demonstrates that 
any homogeneous polynomial  can be decomposed as a sum of 
a homogeneous harmonic polynomial of the same degree and the $|x|^2$-multiple of 
a homogeneous polynomial of degree 2 less. 
We now discuss to what extent such a decomposition holds for arbitrary functions, not necessarily
homogeneous polynomials. 

The main result in this section is the following proposition. 

\begin{prop}\label{prop-decomposition-ana} 
Let $u$ be an analytic function in $B_1\subseteq \mathbb R^n$. 
Then, there exists a unique pair of a harmonic function $h$ and 
an analytic function $u_*$ in some ball $B_R$ such that 
\begin{align}\label{eq-decomposition-ana}u=h+|x|^2 u_*\quad\text{in }B_R,\end{align}
where $R\in (0,1)$ is a positive constant depending only on $n$ and $u$. 
\end{prop}

In the proof, we will demonstrate that $R$ depends only on $n$ 
and the constants in the analyticity estimate of $u$. 

\begin{proof} For the uniqueness, we only need to verify the following statement. 
If $h$ is a harmonic function and $u_*$ is
an analytic function in some ball $B_R$ such that 
$h+|x|^2u_*=0$ in $B_R$, then $h=u_*=0$ in $B_R$. First, by $|x|^2u_*=-h$, we know that $|x|^2u_*$ is harmonic. 
Now, by the analyticity of $u_*$, we can write 
$$u_*=\sum_{k=0}^\infty Q_k,$$
for some $Q_k\in \mathscr P_k$, for each $k\ge 0$. 
Here, the convergence is uniform and absolute in possibly some smaller ball $B_r\subseteq B_R$. 
By examining each homogeneous expression in the expansion of $|x|^2u_*$, 
we conclude that $\Delta (|x|^2Q_k)=0$ for any $k$. 
Of course, this holds in the entire $\mathbb R^n$. 
With $|x|^2Q_k\in \mathscr H_{k+2}^\perp$, 
the injectiveness in Lemma \ref{lemma-Laplace-inj-sur} implies $Q_k=0$, and hence $u_*=0$ in $B_r$.  
By the analyticity again, we obtain $u_*=0$ in $B_R$. Thus, $h=0$ in $B_R$. 
This finishes the proof of the uniqueness. 

We now prove the existence. First, we write
$$u=\sum_{k=0}^\infty P_k,$$
where $P_k$ is the $k$-th homogeneous Taylor polynomial of $u$ at $0$. 
For each $k \ge 0$, $P_k$ has a decomposition given by 
\begin{align}\label{eq-decompostion-P}P_k=\widetilde{P}_k+\widehat{P}_k,\end{align}
for some $\widetilde{P}_k\in \mathscr H_k$
and $\widehat{P}_k \in\mathscr H_k^\perp$. 
Since $u$ is analytic, we have, for any $k \ge 1$,
$$\|P_k\|_{L^2(\mathbb{S}^{n-1})} \leq C^k,$$
for some constant $C>0$ depending only on $n$ and $u$. 
By  
$$\|P_k\|^2_{L^2(\mathbb{S}^{n-1})}
=\|\widetilde{P}_k\|^2_{L^2(\mathbb{S}^{n-1})}+\|\widehat{P}_k\|^2_{L^2(\mathbb{S}^{n-1})},$$
we get 
$$\|\widetilde{P}_k\|_{L^2(\mathbb{S}^{n-1})}, \|\widehat{P}_k\|_{L^2(\mathbb{S}^{n-1})}\le C^k.$$ 
We point out that $\widehat{P}_0 = \widehat{P}_1 = 0$, and for $k \geq 2$,
$$\widehat{P}_k = |x|^2 \widehat{Q}_{k-2},$$
for some $\widehat{Q}_{k-2} \in \mathscr{P}_{k-2}$.
Corollary \ref{cor-constructing-analytic-functions} implies that 
$$h = \sum_{k=0}^\infty \widetilde{P}_k\quad\text{and}\quad |x|^2 u_* = \sum_{k=2}^\infty |x|^2\widehat{Q}_{k-2}$$
converge uniformly and absolutely in $B_R$, for some constant $R\in (0,1)$ depending only on $n$ and $u$. 
Since each $\widetilde{P}_k$ is harmonic, $h$ is harmonic in $B_R$.
Then, \eqref{eq-decompostion-P} and a simple summation yield 
the desired decomposition \eqref{eq-decomposition-ana}.
\end{proof}


The decomposition \eqref{eq-decomposition-ana} generalizes \eqref{eq-decomposition-poly} to arbitrary analytic functions. 
Proposition \ref{prop-decomposition-ana} is only conceptually important 
and will not be essentially used in this paper. 
We now point out an important property already used in the proof of Proposition \ref{prop-decomposition-ana}. 
{\it If $u$ is a harmonic function in $B_1$ with the form $u=|x|^2u_*$ 
for an analytic function $u_*$ in $B_1$, then $u=0$.}
In fact, by using the Taylor series of $u_*$, we can express $u$ by its Taylor series 
$$u=\sum_{k=0}^\infty |x|^2P_k,$$ 
for some $P_k\in \mathscr P_k$, for each $k\ge 0$.
We note that $|x|^2P_k$ cannot be harmonic unless $P_k$ is identically zero. 
It is natural to introduce the following terminology. 

\begin{definition}\label{def-unharmonic-functions} 
Let $u$ be an analytic function in $B_1$. It is called {\it purely unharmonic} at $0$ 
if there exists an analytic function $u_*$ in possibly a smaller ball $B_R\subseteq B_1$ 
such that $u=|x|^2u_*$ in $B_R$. 
\end{definition}

Purely unharmonic functions represent the extreme case that functions cannot be harmonic. 
It is natural to call $h$ and $|x|^2u_*$ in Proposition \ref{prop-decomposition-ana} the {\it harmonic part} 
and the {\it purely unharmonic part} of $u$, respectively. 
We point out that a part of the proof of Proposition \ref{prop-decomposition-ana} 
demonstrates that a nonzero function cannot be harmonic and purely unharmonic simultaneously.

\section{Linear Equations}
\label{sec-linear-eqs}

In this section, we study linear equations and prove the existence and uniqueness of 
purely unharmonic solutions. 

We studied the Poisson equation for homogeneous polynomials 
in Lemma \ref{lemma-Laplace-estimate-poly-L2}. 
We now discuss the Poisson equation for arbitrary functions. 
Let $f$ be a given function defined in $B_1\subseteq\mathbb R^n$ and consider
\begin{align}\label{eq-linear-Poission} 
\Delta u=f.
\end{align}
If $f$ is analytic in $B_1$, then any solutions of \eqref{eq-linear-Poission} in $B_1$ are analytic. 
We can take an arbitrary analytic solution $u$ of \eqref{eq-linear-Poission} 
and apply Proposition \ref{prop-decomposition-ana}. 
As a consequence, the purely unharmonic part of $u$ is a solution of \eqref{eq-linear-Poission}. 
It is easy to verify that it is the unique purely unharmonic solution. 
In fact, the uniqueness does not rely on the existence. 
This line of argument is the motivation of the introduction of purely unharmonic functions 
as in Definition \ref{def-unharmonic-functions}
and leads to the following result.

\begin{theorem}\label{thrm-existence-Poisson}
Let $f$ be an analytic function in $B_1\subseteq\mathbb R^n$. 
Then, there exists a unique analytic function $u_*$ in $B_R$ for some $R\in (0,1)$ 
such that $u=|x|^2u_*$ is a solution of \eqref{eq-linear-Poission} in $B_R$.
\end{theorem} 

\begin{proof}
We only consider the existence. 
The argument before the statement of the theorem 
cannot be generalized to general linear equations or nonlinear equations. 
We present a proof based on Lemma \ref{lemma-Laplace-estimate-poly-L2}. 
For brevity, we write $\|\cdot\|=\|\cdot\|_{L^2(\mathbb S^{n-1})}$, 
the $L^2$-norm on $\mathbb S^{n-1}$. 

We first compute formally. Write 
\begin{align}\label{eq-expansion-existence-Poisson}
u=\sum_{k=2}^\infty P_{k}, \,\, f=\sum_{k=0}^\infty Q_k,\end{align}
for some $P_k\in\mathscr{H}_k^\perp$, for any $k\ge 2$, 
and $Q_k\in\mathscr{P}_k$, for any $k\ge 0$. 
For unification, we set $P_0=P_1=0$. 
By substituting \eqref{eq-expansion-existence-Poisson} in \eqref{eq-linear-Poission}
and arranging all terms as degree-$k$ expressions, we have, for any $k\ge 2$,
\begin{align}\label{eq-linear-degree-k-existence-Poisson}\Delta P_{k}=Q_{k-2}.\end{align} 
By the analyticity of $f$, we assume, for any $i\ge 0$,  
$$\|Q_i\|\le A_0A^i,$$ 
for some positive constants $A_0$ and $A$. 

For any $k\ge 2$, by Lemma \ref{lemma-Laplace-estimate-poly-L2} with $m=0$ and $k$ replaced by $k-2$,  
there is a unique solution $P_{k}\in \mathscr H_{k}^\perp$ of \eqref{eq-linear-degree-k-existence-Poisson}
and 
\begin{align}\label{eq-analyticity-estimate-existence-Poisson}
\|P_{k}\|\le [2(2(k-2)+n)]^{-1}\|Q_{k-2}\|\le \frac{1}{2k}A_0A^k,
\end{align}
where we used $2(k-2)+n=k+(k-2)+(n-2)\ge k$. 
By Corollary \ref{cor-constructing-analytic-functions}, 
we conclude that the series defining $u$ in \eqref{eq-expansion-existence-Poisson}
converges uniformly and absolutely in $B_R$ for some $R\in (0,1)$ 
and that $u$ is analytic in $B_R$. 
A simple summation shows that $u$ is a solution to \eqref{eq-linear-Poission} in $B_R$.

Last, write $P_{k}=|x|^2R_{k-2}$ for some $R_{k-2}\in\mathscr P_{k-2}$, for any $k\ge 2$. 
Then, \eqref{eq-analyticity-estimate-existence-Poisson} yields a similar estimate of $R_k$. 
By Corollary \ref{cor-constructing-analytic-functions} again, we have 
$u=|x|^2u_*$ for some analytic function $u_*$ in $B_R$. In fact, 
$$u_* = \sum_{k=0}^\infty R_k,$$
where the convergence is uniform and absolute in $B_R$. 
\end{proof}

In the following, we study a class of second order linear elliptic equations satisfying a structure condition. Consider
\begin{align}
\label{eq-linear-equation-general}
\Delta u = a^{ijk\ell}x_i x_j \partial_{k\ell} u + b^{ij}x_i \partial_j u + cu + f,
\end{align}
where $a^{ijk\ell}$, $b^{ij}$, $c$, and $f$ are analytic functions in $B_1 \subseteq \R^n$. 
Here, we use the Einstein summation convention.

First, we study the structure of the expression in the right-hand side of \eqref{eq-linear-equation-general}. 
Recall that $\langle \cdot,\cdot \rangle_{*}$ is the inner product defined in \eqref{eq-H^m-inner-product-with-constant}
and $\|\cdot\|_*$ is the corresponding norm. Let $a^{pqrs}_i$, $b^{pq}_i$, $c_i$, 
and $Q_i$ be the $i$-th homogeneous Taylor polynomials of $a^{pqrs}$, $b^{pq}$, $c$, and $f$ at $0$, respectively. 
Then, we can write
\begin{align}
\label{eq-Taylor-series-expression-linear}
a^{ijk\ell}x_i x_j \partial_{k\ell} u + b^{ij}x_i \partial_j u + cu + f = \sum_{k=0}^\infty Q_k + \sum_{i=0}^\infty T_i (u)\quad\text{near }0,
\end{align}
where 
\begin{align}
\label{eq-T_i-operator-linear-equation}
T_i = a^{pqrs}_i x_p x_q \partial_{rs} + b^{pq}_i x_p \partial_q + c_i.
\end{align}
Note that $T_{i}$ is a linear operator from $C^\infty(\R^n)$ into $C^\infty (\R^n)$ 
and maps $\mathscr{P}_{k} $ into $\mathscr{P}_{k + i}$, for any $k \ge 0$. 

In the following, we prove the existence of purely unharmonic solutions 
to \eqref{eq-linear-equation-general} under an additional condition. 
It is motivated by our study of the minimal surface equation and the Monge-Ampère equation in exterior domains.

First, it is convenient to introduce the following terminology. 

\begin{definition}\label{def-degree-|x|^2}
Let $P \in \mathscr{P}_k$ be a homogeneous polynomial of degree $k$ with $P \in |x|^{2\ell} \mathscr{P}_{k - 2\ell}$ 
for some $0\le \ell \le k/2$.
The largest such $\ell$ is called the {\it degree of $|x|^2$} of $P$ and is denoted by
\begin{align*}
\ell = \deg_{|x|^2} P.
\end{align*}
For $P = 0$, we define $\deg_{|x|^2} P = \infty$.
\end{definition}

If $P \in \mathscr{H}_k^\perp$ satisfies $\Delta P = Q$ for some $Q \in \mathscr{P}_{k-2}$, then 
$$\deg_{|x|^2}P = \deg_{|x|^2}Q +1.$$

Next, we introduce the following terminology.
\begin{definition}
\label{def-operators-of-type-i}
Fix an integer $i \ge 0$. 
A linear operator $T:C^\infty(\R^n) \rightarrow C^\infty(\R^n)$ 
is called to be {\it of type $i$ that preserves the degree of $|x|^2$} if
\begin{itemize}
    \item[(i)] $T:\mathscr{P}_{k} \rightarrow \mathscr{P}_{k+ i}$
    for any $k\ge 0$; 
    \item[(ii)] for any $u \in C^\infty(\R^n)$ and $l \ge 0$,
    $$T\big(|x|^{2l}u \big) = |x|^{2l} S_{l}(u),$$
    for some operator $S_l:C^\infty(\R^n) \rightarrow C^\infty(\R^n).$
\end{itemize}
\end{definition}

For convenience, we say that $T_i$ {\it preserves the degree of $|x|^2$} for homogeneous polynomials, 
although, strictly speaking, $T_i$ does not decrease the degree of $|x|^2$.

Now, we are ready to prove the main result of this section.

\begin{theorem}\label{thrm-existence-linear-equation}
Let $a^{ijk\ell}$, $b^{ij}$, $c$, and $f$ be analytic functions in $B_1$, for $i,j,k,\ell =1, \cdots, n$, 
such that \eqref{eq-Taylor-series-expression-linear} holds near $0$, 
for some sequence of homogeneous polynomials $\{Q_k\}_{k\ge 0}$, with $Q_k \in \mathscr{P}_k$, 
and operators $T_i$ of type $i$ that preserve the degree of $|x|^2$, for $i \ge 0$. 
Then, for any harmonic function $h$ in $B_1$, there exists a unique analytic function $u_*$ in $B_R$ for some $R \in (0,1)$ 
such that $u =h+ |x|^2 u_*$ is a solution of \eqref{eq-linear-equation-general} in $B_R$.
\end{theorem}

We point out that the key assumption in Theorem \ref{thrm-existence-linear-equation} 
is the preservation of the degree of $|x|^2$ by $T_{i}$.

\begin{proof} 
By introducing $u\mapsto u-h$, we consider the case $h=0$ and prove 
that there exists a unique analytic function $u_*$ in $B_R$ 
such that $u=|x|^2u_*$ is a solution of \eqref{eq-linear-equation-general} in $B_R$, 
for some $R\in (0,1)$. 
Recall that $\langle \cdot,\cdot \rangle_{*}$ is the inner product defined in \eqref{eq-H^m-inner-product-with-constant}
and $\|\cdot\|_*$ is the corresponding norm. 
We will use the estimate \eqref{eq-Banach-algebra-property-without-constant} 
repeatedly.

Assume that $u = |x|^2 u_*$ is a solution to \eqref{eq-linear-equation-general} near $0$ 
with the right-hand side given by \eqref{eq-Taylor-series-expression-linear}. 
Let $\{P_k\}_{k\ge 2}$ be the sequence of the homogeneous Taylor polynomials of $u$ at $0$, 
with $P_k \in \mathscr{H}_k^\perp$. It is straightforward to check that, for any $k \ge 2$,
\begin{align}
\label{eq-recursive-equation-linear}
\Delta P_k = \sum_{i+j = k-2} T_iP_j + Q_{k-2}.
\end{align}
By $P_k \in \mathscr{H}_k^\perp$ for any $k \ge 2$, $P_k$ is uniquely determined by $Q_{k-2}$ and $P_2,\cdots, P_{k-2}$. 
Hence, any solution $u = |x|^2 u_*$ of \eqref{eq-linear-equation-general} for some analytic function $u_*$ near $0$ 
with the right-hand side of \eqref{eq-linear-equation-general} given by \eqref{eq-Taylor-series-expression-linear} 
has the same Taylor series at $0$. This proves the uniqueness. 
We point out that the uniqueness does not require the preservation of the degree of $|x|^2$ by $T_i$.

Set $P_0 = P_1 = 0$. 
By using Lemma \ref{lemma-Laplace-estimate-poly-H*} successively for $k=2, 3, \cdots$, we conclude 
that there is a unique solution $P_{k}\in \mathscr H_{k}^\perp$ of 
\eqref{eq-recursive-equation-linear}. 
It suffices to prove that, for any $k \geq 2$,
\begin{align}
\label{eq-exp-growth-P_k-linear}
\|P_k\|_* \le C_*^k,
\end{align}
for some positive constant $C_*$.
Then, the existence follows from Corollary \ref{cor-constructing-analytic-functions}. 

Recall that $T_i$ is given by \eqref{eq-T_i-operator-linear-equation}. 
Since $a^{pqrs}$, $b^{pq}$, and $c$ are analytic, Lemma \ref{lemma-homo-poly-first-derivatives-H*} yields, 
for any $i \ge 0$, $k\ge 2$, 
and $P \in \mathscr{P}_{k}$, 
\begin{align}
\label{eq-estimate-T_i P-linear}
\|T_i P\|_* \leq A_0A^i k (k-1) \|P\|_*, 
\end{align}
for some constants $A_0>0$ and $A>1$, depending only on $n$, $a^{pqrs}$, $b^{pq}$, and $c$. 
Similarly, since $f$ is analytic, then
\begin{align}
\label{eq-estimate-Q_j-linear}
\|Q_i\|_* \le B_0B^i,
\end{align}
for some constants $B_0>0$ and $B>1$, depending only on $n$ and $f$. 
By \eqref{eq-recursive-equation-linear},  we get, for any $k\ge 2$, 
\begin{align}
\label{eq-recursive-zeroth-iteration-linear}
P_k = \Delta_k^{-1} Q_{k-2} +  \sum_{i+j = k-2} \Delta_k^{-1} T_i P_{j}.
\end{align}
Here, $P_k$ is expressed in terms of $Q_{k-2}$ and $P_{j}$ for $j\le k-2$. 
We can replace $P_k$ by $P_{j}$ in \eqref{eq-recursive-zeroth-iteration-linear}  
and substitute in \eqref{eq-recursive-zeroth-iteration-linear}. 
By continuing this process, we get an expression of $P_k$ in terms of $Q_i$ for $i=0, 1, \cdots, k-2$. 
We introduce some notations to simplify such an expression. 

We define, for any $k\ge 0$, 
\begin{align}\label{eq-def-R_k,0-linear}R_{k, 0}=Q_{k},\end{align}
and inductively, for any $1\le \ell\le k/2$, 
\begin{align*}
R_{k, \ell}=\sum_{\substack{i+j=k \\ j\ge 2\ell}}T_i\Delta_j^{-1}R_{j-2, \ell-1}.\end{align*}
For convenience, we also set, for any $\ell>k/2$, 
$$R_{k, \ell}=0.$$
Hence, for any $1\le \ell\le k/2$, 
\begin{align}\label{eq-def-R_k,l-linear}
R_{k, \ell}=\sum_{i+j=k}T_i\Delta_j^{-1}R_{j-2, \ell-1}.\end{align}
We claim, for any $k \ge 2$, 
\begin{align}
\label{eq-further-decomp-of-P_k-linear}
P_k=\sum_{\ell=0}^{[k/2]- 1}\Delta_k^{-1}R_{k-2, \ell}. 
\end{align}
First, \eqref{eq-further-decomp-of-P_k-linear} holds for $k=2, 3$ 
by \eqref{eq-recursive-zeroth-iteration-linear} and $P_0=P_1=0$. 
We prove \eqref{eq-further-decomp-of-P_k-linear} for $k\ge 4$ by induction. 
Assume, for any $2\le j<k$, 
\begin{align*}
P_j=\sum_{\ell=0}^{[j/2]- 1}\Delta_j^{-1}R_{j-2, \ell}. 
\end{align*}
By a simple substitution of $P_j$ in \eqref{eq-recursive-zeroth-iteration-linear}, we have 
\begin{align*}
\Delta P_k& = Q_{k-2} +  \sum_{i+j = k-2}\sum_{ \ell=0}^{[j/2]- 1}T_i  \Delta_j^{-1}R_{j-2, \ell}\\
& = Q_{k-2} +  \sum_{i+j = k-2}\sum_{ \ell=1}^{[j/2]}T_i  \Delta_j^{-1}R_{j-2, \ell-1}. 
\end{align*}
We point out that $R_{j-2,\ell-1} = 0$ if $\ell >[ j/2 ]$. 
By switching the order of summations, we get 
\begin{align*}
\Delta P_k  = Q_{k-2} +\sum_{ \ell=1}^{[k/2]-1}\sum_{i+j = k-2}T_i  \Delta_j^{-1}R_{j-2, \ell-1}.
\end{align*}
The second sum in the right-hand side is $R_{k-2,\ell}$ by \eqref{eq-def-R_k,l-linear}. 
Hence, 
\begin{align*}
\Delta P_k  = R_{k-2,0} +\sum_{ \ell=1}^{[k/2]-1}R_{k-2,\ell} = \sum_{ \ell=0}^{[k/2]-1} R_{k-2,\ell}.
\end{align*}
This proves \eqref{eq-further-decomp-of-P_k-linear}.

We now proceed to estimate $R_{k, \ell}$ for each pair $k, \ell$ with $0\le \ell\le k/2$. 
First, it is easy to verify $R_{k, \ell}\in |x|^{2\ell}\mathscr{P}_{k-2\ell}\subseteq\mathscr{P}_{k}$ for any such $k, \ell$. 
Take any $1\le \ell\le j/2\le k/2$. 
By $R_{j-2, \ell-1}\in |x|^{2(\ell-1)}\mathscr{P}_{j -2-2(\ell-1)}\subseteq \mathscr P_{j-2}$
and Lemma \ref{lemma-estimate-degree-|x|2-case2-H*}, 
with $k = j-2$ and $\ell$ replaced by $\ell-1$, we get 
$$\|\Delta_{j}^{-1}R_{j-2, \ell-1}\|_*\leq \frac{1}{[2(\ell-1)+2\big]\big[ 2(j-2) -2(\ell-1) + n \big]} \|R_{j-2, \ell-1}\|_*.$$
It is obvious that 
$$2(j-2) -2(\ell-1) + n=j+(j-2\ell)+ (n-2)\ge j. $$
Then, 
$$\|\Delta_{j}^{-1}R_{j-2, \ell-1}\|_*\leq \frac{1}{2\ell j} \|R_{j-2, \ell-1}\|_*.$$
By \eqref{eq-estimate-T_i P-linear} and 
$\Delta_{j}^{-1}R_{j-2, \ell-1}\in \mathscr{P}_{j}$, we have
\begin{align*}\|T_i \Delta_{j}^{-1}R_{j-2, \ell-1}\|_*
&\leq A_0A^ij(j-1) \|\Delta_{j}^{-1}R_{j-2, \ell-1}\|_*\\
&\le \frac{A_0A^i(j-1)}{2\ell}\|R_{j-2, \ell-1}\|_*
\le \frac{k-1}{2\ell}A_0A^i\|R_{j-2, \ell-1}\|_*.\end{align*}
By \eqref{eq-def-R_k,l-linear}, we obtain 
\begin{align*}\|R_{k, \ell}\|_*\le  \frac{k-1}{2\ell}\sum_{i+j=k}A_0A^i\|R_{j-2, \ell-1}\|_* 
= \frac{k-1}{2\ell}\sum_{i+j=k-2}A_0A^i\|R_{j, \ell-1}\|_*.\end{align*}
For ease of iterations, we relabel indices and write 
\begin{align*}\|R_{k, \ell}\|_*\le  \frac{k-1}{2\ell}\sum_{i_1+j_2=k-2}A_0A^{i_1}\|R_{j_2, \ell-1}\|_*.\end{align*}
Then, for any $j_2$ allowed in the sum, we get 
\begin{align*}\|R_{j_2, \ell-1}\|_*\le  \frac{j_2-1}{2(\ell-1)}\sum_{i_2+j_3=j_2-2}A_0A^{i_2}\|R_{j_3, \ell-2}\|_*.\end{align*}
Note that $j_2\le k-2$. A simple substitution yields 
\begin{align*}\|R_{k, \ell}\|_*\le  \frac{(k-1)(k-3)}{2^2\ell(\ell-1)}\sum_{i_1+i_2+j_3=k-4}A_0^2A^{i_1+i_2}\|R_{j_3, \ell-2}\|_*.\end{align*}
After $\ell$ times, we have 
\begin{align*}\|R_{k, \ell}\|_*\le  \frac{(k-1)(k-3)\cdots(k-2\ell+1)}{2^\ell\ell!}
\sum_{i_1+\cdots+i_\ell+j_{\ell+1}=k-2\ell}A_0^\ell A^{i_1+\cdots+i_\ell}\|R_{j_{\ell+1}, 0}\|_*.\end{align*}
It is obvious that 
$$(k-1)(k-3)\cdots(k-2\ell+1)\le(2k)(2k-2)\cdots(2k-2\ell+2)= \frac{2^\ell k!}{(k-\ell)!}.$$
By \eqref{eq-estimate-Q_j-linear} and \eqref{eq-def-R_k,0-linear}, 
we obtain
\begin{align*}\|R_{k, \ell}\|_*
\leq \frac{k!}{\ell!(k-\ell)!} A_0^\ell B_0\sum_{m+j=k-2\ell} A^{m}B^{j}\sum_{i_1+\cdots+i_\ell=m}1.\end{align*} 
This holds for any $1\le \ell\le k/2$.  
By the binomial theorem, we have 
\begin{align*}\frac{k!}{\ell!(k-\ell)!} \leq  2^k.
\end{align*} 
In addition, the number of nonnegative integers $i_1,\cdots,i_{\ell}$ with 
$i_1 + \cdots + i_{\ell}=m$ is given by 
$$\frac{(m+\ell-1)!}{m!(\ell-1)!},$$
and, by the binomial theorem again, is bounded by $2^{m+\ell}$.
Hence, 
for any $1\le \ell\le k/2$, 
\begin{align*}\|R_{k, \ell}\|_*
\leq B_0 2^k (2A_0)^\ell \sum_{m+j=k-2\ell} (2A)^{m}B^{j} 
\le B_0 2^k (2A_0)^\ell (2A + B)^{k-2\ell}.\end{align*} 
Thus, by choosing a constant $D>2A+B$, we have,
for any $1 \leq \ell \leq k/2$,  
\begin{align}\label{eq-estimate-R_k,ell-linear-0}
\|R_{k, \ell}\|_*\leq  B_0 2^k (2A_0)^\ell D^{k-2\ell}.\end{align} 
This obviously holds for $\ell=0$. 

Next, in \eqref{eq-further-decomp-of-P_k-linear}, we have 
$R_{k-2, \ell}\in  |x|^{2\ell}\mathscr P_{k-2-2\ell}\subseteq \mathscr P_{k-2}$ 
for any $k\ge 2$ and $0\le \ell\le k/2- 1$. 
By Lemma \ref{lemma-estimate-degree-|x|2-case2-H*} again, we get, 
for any such $k$ and $\ell$,  
\begin{align*} 
\|\Delta_k^{-1}R_{k-2, \ell}\|_*
\le \big[(2\ell+2)\big(2(k-2)-2\ell+n\big)\big]^{-1}\|R_{k-2, \ell}\|_*.
\end{align*}
It is obvious that 
$$2(k-2)-2\ell+n=k+ (k-2\ell-2)+n-2\ge k.$$
With \eqref{eq-estimate-R_k,ell-linear-0}, we have 
\begin{align*} 
k\|\Delta_k^{-1}R_{k-2, \ell}\|_*
\le \frac{B_0}{2\ell+2}2^{k-2}(2A_0)^\ell D^{k-2\ell-2}.
\end{align*}
By a simple substitution in \eqref{eq-further-decomp-of-P_k-linear}, we obtain 
\begin{align}\label{eq-final-estimate-linear-k}
k\|P_k\|_* \leq   B_0(2D)^{k-2}\sum_{\ell=0}^{[k/2]- 1}\frac{1}{2\ell+2} \Big(\frac{2A_0}{D^2}\Big)^\ell 
\le  B_0E_0(2D)^{k-2},
\end{align}
for some positive constant $E_0$, if $D$ is large such that $D^2>2A_0$. 
This finishes the proof of \eqref{eq-exp-growth-P_k-linear}.
\end{proof}


To end this section, we discuss a special class of linear equations. 
By choosing $a^{ijk\ell}$ and $b^{ij}$ in \eqref{eq-linear-equation-general} appropriately, we have the following more 
specific form:  
\begin{align}\label{eq-linear-equation-general-special-factors}\begin{split} 
\Delta u &= a|x|^2\Delta u+\widetilde{a} x_kx_\ell\partial_{k\ell}u+a^k|x|^2 x_\ell \partial_{k\ell}u
+a^{k\ell}|x|^4 \partial_{k\ell}u\\
&\qquad+ b  \langle x, \nabla u\rangle+ b^k|x|^2\partial_ku + cu + f,
\end{split}\end{align}
for some analytic functions $a, \widetilde{a}, a^k, a^{k\ell}$, $b, b^k$, $c$, and $f$ in $B_1$, for $k, \ell=1, \cdots,n$. 
The equation \eqref{eq-linear-equation-general-special-factors} is of particular interest to us, 
due to its close connections with the minimal surface equation and the Monge-Amp\`{e}re equation.

\begin{lemma}\label{lemma-T-preserving} 
Let $a, \widetilde{a}, a^k, a^{k\ell}$, $b, b^k$, and $c$ be analytic functions in $B_1$, 
for $k,\ell=1, \cdots,n$, 
and $T_i$ be given by \eqref{eq-T_i-operator-linear-equation} 
for the equation \eqref{eq-linear-equation-general-special-factors}, for any $i\ge 0$,  
Then, $T_i$ is an operator of type $i$ that preserves the degree of $|x|^2$. 
\end{lemma} 

\begin{proof}
Take any $u \in C^\infty(\R^n)$ and positive integer $\ell$. Then, 
\begin{align}
\label{eq-partial_i(|x|^2l u)}
\partial_i \big(|x|^{2\ell} u\big) = 2\ell|x|^{2\ell-2} x_i u + |x|^{2\ell} \partial_i u,
\end{align}
and
\begin{equation}
\label{eq-partial_ij(|x|^2l u)}
\begin{aligned}
\partial_{ij}\big(|x|^{2\ell} u\big) &= 2\ell(2\ell-2)|x|^{2\ell-4} x_i x_j u +  2\ell \delta_{ij}|x|^{2\ell-2} u\\
&\qquad+ 2\ell|x|^{2\ell-2}  (x_i \partial_j u + x_j \partial_i u) + |x|^{2\ell} \partial_{ij}u.
\end{aligned}
\end{equation}
The desired conclusion follows from a straightforward computation.
\end{proof}

As a consequence, the existence as in Theorem \ref{thrm-existence-linear-equation} 
holds for \eqref{eq-linear-equation-general-special-factors}. 

\section{Model Equations}
\label{sec-model-fully-nonlinear}

In this section, we discuss several model equations. 
Our main purpose is to present ideas of proving the existence of purely unharmonic solutions
and discuss encountered difficulties. 

Recall that $\langle \cdot,\cdot \rangle_{*}$ is the inner product defined in \eqref{eq-H^m-inner-product-with-constant}
and $\|\cdot\|_*$ is the corresponding norm. 
We will use the estimate \eqref{eq-Banach-algebra-property-without-constant} repeatedly.

First, we consider 
\begin{align}\label{eq-nonlinear-special-case3} 
\Delta u=|x|^2u\Delta u+f.\end{align}

\begin{prop}
\label{prop-nonlinear-special-case3}
Let $f$ be an analytic function in $B_1\subseteq \R^n$. 
Then, there exists a unique analytic function $u_*$ in $B_R$ 
such that $u=|x|^2u_*$ is a solution of \eqref{eq-nonlinear-special-case3} in $B_R$, 
where $R \in (0,1)$ is a positive constant depending only on $n$ and the analyticity of $f$. 
\end{prop}

\begin{proof} 
We first consider the uniqueness. 
Assume that $u=|x|^2u_*$ is a solution of \eqref{eq-nonlinear-special-case3} 
for some analytic function $u_*$ near 0, 
and $\{P_k\}_{k\geq 0}$ is the sequence of the homogeneous Taylor polynomials of $u$ at $0$. 
Then, $P_0=P_1=0$ and $P_k\in  \mathscr H_k^\perp$ for any $k\ge 2$. 
Let $Q_i$ be the $i$-th homogeneous Taylor polynomial of $f$ 
at $0$.  
Replacing each term in \eqref{eq-nonlinear-special-case3} by its Taylor series 
and collecting terms of the same degree, we have, for any $k \geq 2$,
\begin{align}
\label{eq-recursive-eq-nonlinear-special-case3}
\Delta P_k = |x|^2 \sum_{i+j = k-2} P_i \Delta P_j + Q_{k-2}.
\end{align}
By $P_k\in  \mathscr H_k^\perp$ for any $k\ge 2$, 
$P_k$ is uniquely determined by $Q_{k-2}$ and $P_{i}$ for $i\le k-2$. 
Since $P_0 = P_1 = 0$, the sequence $\{P_k\}_{k\geq 0}$ is uniquely determined by $f$. 
In summary, every solution $u=|x|^2u_*$ to \eqref{eq-nonlinear-special-case3}  
for some analytic function $u_*$ near 0
has the same Taylor series at $0$. This proves the uniqueness.

Next, we study the existence. 
We take $P_0=P_1=0$.
By using Lemma \ref{lemma-Laplace-estimate-poly-H*} successively for $k=2, 3, \cdots$, we conclude 
that there is a unique solution $P_{k}\in \mathscr H_{k}^\perp$ of 
\eqref{eq-recursive-eq-nonlinear-special-case3}. 
It suffices to prove, for any $k \ge 2$,
\begin{align}
\label{eq-exp-growth-P_k-nonlinear-special-case3}
\|P_k\|_* \le C_*^k,
\end{align}
for some positive constant $C_*$.
Then, the existence follows from Corollary \ref{cor-constructing-analytic-functions}. 
By the analyticity of $f$, we assume, for any $i\ge 0$,  
$$\|Q_i\|_* \leq A_0A^i,$$
for some positive constants $A_0$ and $A$. 

By successively eliminating $\Delta P_j$ in the right-hand side 
of \eqref{eq-recursive-eq-nonlinear-special-case3} for $j=k-2, \cdots, 2$, we have,  
for any $k\ge 2$, 
\begin{equation*}
\Delta P_k = \sum_{2\ell + r = k-2} |x|^{2\ell} \sum_{i_1 + \cdots + i_\ell + j = r} P_{i_1}\cdots P_{i_\ell} Q_{j}+Q_{k-2},
\end{equation*}
and hence
\begin{equation}
\label{eq-recursive-equation-nonlinear-special-case3}
P_k = \sum_{2\ell + r = k-2}  \sum_{i_1 + \cdots + i_\ell + j = r}\Delta_k^{-1}(|x|^{2\ell} P_{i_1}\cdots P_{i_\ell} Q_{j})
+\Delta_k^{-1}Q_{k-2}.
\end{equation}
We note that the summation is over $1 \le \ell \le k/2 -1$ and appears only for $k\ge 6$.
By applying Lemma \ref{lemma-estimate-degree-|x|2-case2-H*} to each term 
in \eqref{eq-recursive-equation-nonlinear-special-case3}, with $k$ replaced by $k-2$, 
we have, for any $k\ge 2$, 
\begin{align*}
\|P_k\|_* & \leq \sum_{2\ell + r = k-2} \frac{1}{(2\ell + 2)(2\ell + 2r + n)} 
\sum_{i_1 + \cdots + i_\ell + j = r}  \|P_{i_1}\|_*\cdots  \|P_{i_\ell}\|_*  \|Q_j\|_*\\
&\qquad+\frac{1}{2(2(k-2) + n)} \|Q_{k-2}\|_*,
\end{align*}
where we used $2(k-2)-2\ell + n = 2\ell + 2r + n$.
It is obvious that, for $2\ell+r = k-2$,
$$2\ell + 2r + n = k+(n-2) +r \geq k.$$
Hence,  for any $k \geq 2$, 
\begin{align*} 
k \|P_k\|_* 
\leq A_0\sum_{2\ell + r = k-2} \frac{1}{2\ell + 2} \sum_{i_1 + \cdots + i_\ell + j = r}  \|P_{i_1}\|_*\cdots  \|P_{i_\ell}\|_*  A^{j}
+\frac12A_0A^{k-2}.
\end{align*}
By renaming indices, we obtain, for any $k \geq 2$, 
\begin{align}\label{eq-recursive-estimate-last-nonlinear-special-case3}
k \|P_k\|_* \leq A_0\sum_{2\ell + r +j = k-2} \frac{A^j}{2\ell + 2} \sum_{i_1 + \cdots + i_\ell = r}  \|P_{i_1}\|_*\cdots  \|P_{i_\ell}\|_*
+\frac12A_0A^{k-2}.
\end{align}

In view of \eqref{eq-recursive-estimate-last-nonlinear-special-case3}, 
we set $a_0=a_1=0$ and inductively, for any $k\ge 2$, 
\begin{align}\label{eq-recursive-equation-ode-majorant-nonlinear-special-case3} 
k a_k = A_0\sum_{2\ell + r +j = k-2} \frac{A^j}{2\ell + 2} \sum_{i_1 + \cdots + i_\ell = r}  a_{i_1}\cdots  a_{i_\ell}+\frac12A_0A^{k-2}.
\end{align} 
Consider formally 
\begin{align}\label{eq-formal-expansion-nonlinear-special-case3} 
\varphi(t) = \sum_{k=2}^\infty a_k t^k\quad\text{near }t=0.\end{align}
First, it is obvious that the left-hand side of 
\eqref{eq-recursive-equation-ode-majorant-nonlinear-special-case3}  is the $(k-1)$-th coefficient of the power series of 
$\varphi'(t)$, for any $k \ge 2$. 
Next, note that the summation over $i_1 + \cdots + i_\ell = r$ is the $r$-th coefficient of the power series of $\varphi^\ell$, 
for $\ell \ge 1$. 
Set 
$$\beta(t)= A_0\sum_{i=0}^\infty A^it^i=\frac{A_0}{1-At}\quad\text{near }t=0.$$
Then, $\beta$ is analytic near $t=0$. 
Hence, the right-hand side of \eqref{eq-recursive-equation-ode-majorant-nonlinear-special-case3} is the $(k-2)$-th coefficient of 
$$\beta(t)\sum_{\ell = 0}^\infty \frac{1}{2\ell + 2} t^{2\ell} \big(\varphi(t)\big)^{\ell}.$$
We next define 
$$G(t,s) = \sum_{\ell=0}^\infty \frac{1}{2\ell+2} t^{2\ell} s^{\ell}\quad\text{near }(t,s)=(0,0).$$
Then, $G$ is an analytic function near $0 \in \R^2$. 
Therefore, for any $k \ge 2$, the right-hand side of 
\eqref{eq-recursive-equation-ode-majorant-nonlinear-special-case3} is the $(k-2)$-th coefficient of 
the power series of $\beta(t)G(t, \varphi(t)).$
We are led to an identity  
\begin{align}\label{eq-formal-diff-ineq-nonlinear-special-case3} 
t\varphi' = t^2 \beta(t)G(t, \varphi).\end{align}
We point out that the derivation of 
\eqref{eq-formal-diff-ineq-nonlinear-special-case3} is formal but the functions $\beta$ and $G$ 
are defined rigorously as analytic functions. 

We now derive \eqref{eq-exp-growth-P_k-nonlinear-special-case3}
from \eqref{eq-recursive-estimate-last-nonlinear-special-case3}. 
Inspired by \eqref{eq-formal-diff-ineq-nonlinear-special-case3}, we consider the initial value problem
\begin{align}
\label{eq-majorant-ode-nonlinear-special-case3}
\varphi' = t \beta(t)G(t, \varphi),\quad \varphi(0) = 0.
\end{align}
It is clear that \eqref{eq-majorant-ode-nonlinear-special-case3} 
admits a unique analytic solution $\varphi = \varphi(t)$ near $t = 0$, with $\varphi'(0) = 0$. 
Write $\varphi$ as in \eqref{eq-formal-expansion-nonlinear-special-case3}. 
Then, for any $k \ge 2$, \eqref{eq-recursive-equation-ode-majorant-nonlinear-special-case3} holds 
and 
\begin{align}
\label{eq-exp-growth-majorant-ode-sol-nonlinear-special-case3}
|a_k| \leq C_*^k,
\end{align}
for some positive constant $C_*$.  
Inductively, we can prove that $a_k\ge 0$ for any $k\ge 2$. 

We next prove that, for any $k \ge 2$,
\begin{align}
\label{eq-P_k<=a_k-nonlinear-special-case3}
\|P_k\|_* \le a_k.
\end{align}
First, by taking $k=2$ in \eqref{eq-recursive-estimate-last-nonlinear-special-case3}
and \eqref{eq-recursive-equation-ode-majorant-nonlinear-special-case3},  
we have \eqref{eq-P_k<=a_k-nonlinear-special-case3} for $k=2$. 
Next, a simple induction  argument based on 
\eqref{eq-recursive-estimate-last-nonlinear-special-case3}
and \eqref{eq-recursive-equation-ode-majorant-nonlinear-special-case3} yields 
\eqref{eq-P_k<=a_k-nonlinear-special-case3} for any $k\ge 2$. 
Finally, by \eqref{eq-exp-growth-majorant-ode-sol-nonlinear-special-case3}
and \eqref{eq-P_k<=a_k-nonlinear-special-case3}, we obtain \eqref{eq-exp-growth-P_k-nonlinear-special-case3}.
\end{proof}

Next, consider 
\begin{align}\label{eq-nonlinear-special-case4} 
\Delta u=ux_ix_j\partial_{ij} u+f.\end{align}
Assume that $u=|x|^2u_*$ is a solution of \eqref{eq-nonlinear-special-case4} 
for some analytic function $u_*$ near 0, 
and $\{P_k\}_{k\geq 0}$ is the sequence of the homogeneous Taylor polynomials of $u$ at $0$. 
Then, $P_0=P_1=0$ and $P_k\in  \mathscr H_k^\perp$ for any $k\ge 2$. 
Let $Q_i$ be the $i$-th homogeneous Taylor polynomial of $f$ 
at $0$.  
Replacing each term in \eqref{eq-nonlinear-special-case4} by its Taylor series 
and collecting terms of the same degree, we have, for any $k \geq 2$,
\begin{align}
\label{eq-recursive-eq-nonlinear-special-case4}
\Delta P_k = \sum_{i+j = k-2} j(j-1)P_i P_j + Q_{k-2}.
\end{align}
This is the recursive equation for \eqref{eq-nonlinear-special-case4}. 
As in the proof of Proposition \ref{prop-nonlinear-special-case3}, we can get, for any $k\ge 2$, 
\begin{align*}
k\|P_k\|_*\le \frac12 \sum_{i+j = k-2} j(j-1)\|P_i\|_* \|P_j\|_* + \frac12\|Q_{k-2}\|_*.
\end{align*}
Then, a formal argument leads to the recursive relation 
\begin{align*}
ka_k= \frac12 \sum_{i+j = k-2} j(j-1)a_ia_j + \frac12A_0A^{k-2},
\end{align*}
for some constants $A_0,A>0$, and then to the equation 
\begin{align}\label{eq-formal-diff-ineq-nonlinear-not-working-special-case4} 
\varphi' = \frac12t^3 \varphi\varphi''+t\beta,\end{align}
for some analytic function $\beta$. 
The equation \eqref{eq-formal-diff-ineq-nonlinear-not-working-special-case4} can be interpreted 
as a perturbation of $\varphi'$ by $\varphi''$. 
It is unlikely that \eqref{eq-formal-diff-ineq-nonlinear-not-working-special-case4} admits an analytic solution. 
Since $\varphi''$ appears in the right-hand side of the equation, although as a perturbative term, 
we need a $\varphi''$ also in the left-hand side. 
A more sophisticated argument is needed. 

\begin{prop}
\label{prop-nonlinear-special-case4}
Let $f$ be an analytic function. 
Then, there exists a unique analytic function $u_*$ in $B_R$ for some $R \in (0,1)$ 
such that $u=|x|^2u_*$ is a solution of \eqref{eq-nonlinear-special-case4} in $B_R$, 
where $R \in (0,1)$ is a constant depending only on $n$ and the analyticity of $f$. 
\end{prop}

\begin{proof} To get the uniqueness, we can proceed similarly as in the proof of Proposition \ref{prop-nonlinear-special-case3}. 
We only record the recursive equations. 
Assume that $u=|x|^2u_*$ is a solution of \eqref{eq-nonlinear-special-case4} 
for some analytic function $u_*$ near 0, 
and $\{P_k\}_{k\geq 0}$ is the sequence of the homogeneous Taylor polynomials of $u$ at $0$. 
Then, $P_0=P_1=0$ and $P_k\in  \mathscr H_k^\perp$ for any $k\ge 2$. 
Let $Q_i$ be the $i$-th homogeneous Taylor polynomial of $f$ 
at $0$.  The recursive equation is given by \eqref{eq-recursive-eq-nonlinear-special-case4}, i.e., 
\begin{align*}
\Delta P_k = \sum_{i+j = k-2} j(j-1)P_i P_j + Q_{k-2}.
\end{align*} 

Next, we study the existence. 
We take $P_0=P_1=0$.
By using Lemma \ref{lemma-Laplace-estimate-poly-H*} successively for $k=2, 3, \cdots$, we conclude 
that there is a unique solution $P_{k}\in \mathscr H_{k}^\perp$ of 
\eqref{eq-recursive-eq-nonlinear-special-case4}. 
We write 
\begin{align}\label{eq-recursive-equiv-nonlinear-special-case4}
P_k = \sum_{i+j = k-2} j(j-1)\Delta_k^{-1}(P_i P_j) + \Delta_k^{-1}Q_{k-2}.
\end{align}
It suffices to prove, for any $k \ge 2$,
\begin{align}
\label{eq-exp-growth-P_k-nonlinear-special-case4}
\|P_k\|_* \le C_*^k,
\end{align}
for some positive constant $C_*$.
Then, the existence follows from Corollary \ref{cor-constructing-analytic-functions}. 
By the analyticity of $f$, we assume, for any $i\ge 0$,  
$$\|Q_i\|_* \leq A_0A^i,$$
for some positive constants $A_0$ and $A$. 

For any $k\ge 2$, by $P_k\in |x|^2\mathscr P_{k-2}$, we assume 
\begin{align}\label{eq-decomposition-P-nonlinear-special-case4}
P_k=\sum_{1\le \ell\le k/2}R_{k,\ell},\end{align}
for some $R_{k,\ell}\in |x|^{2\ell}\mathscr P_{k-2\ell}$. 
For unification, we set 
$$R_{k,\ell}=0\quad\text{for $\ell=0$ or $\ell>[k/2]$}.$$ 
A simple substitution in \eqref{eq-recursive-equiv-nonlinear-special-case4} yields 
\begin{align*}
P_k = \sum_{i+j = k-2}\sum_{\substack{p\le i/2\\ q\le j/2}} j(j-1)\Delta_k^{-1}(R_{i,p}R_{j,q}) + \Delta_k^{-1}Q_{k-2}.
\end{align*}
Note that $R_{i,p}R_{j,q}\in |x|^{2(p+q)}\mathscr P_{i+j-2(p+q)}$ and hence 
$$\Delta_k^{-1}(R_{i,p}R_{j,q})\in |x|^{2(p+q+1)}\mathscr P_{i+j-2(p+q)}.$$
For any $2\le \ell\le k/2$, we consider all pairs $(p,q)$ with $p+q=\ell-1$ 
and write  
\begin{align*}P_k = \sum_{i+j = k-2}\sum_{2\le \ell\le k/2}\sum_{p+q=\ell-1} j(j-1)\Delta_k^{-1}(R_{i,p}R_{j,q}) 
+ \Delta_k^{-1}Q_{k-2}.
\end{align*}
By comparing with \eqref{eq-decomposition-P-nonlinear-special-case4}, we have, for any $k\ge 2$, 
\begin{align}\label{eq-expression-P-k1-nonlinear-special-case4}R_{k,1}= \Delta_k^{-1}Q_{k-2},\end{align}
and, for any $k\ge 2$ and $2\le\ell\le k/2$, 
\begin{align}\label{eq-expression-P-kl-nonlinear-special-case4}
R_{k,\ell}=\sum_{i+j = k-2}\sum_{p+q=\ell-1} j(j-1)\Delta_k^{-1}(R_{i,p}R_{j,q}).\end{align}
In fact, we should define $R_{k,\ell}$ by \eqref{eq-expression-P-k1-nonlinear-special-case4}
and \eqref{eq-expression-P-kl-nonlinear-special-case4} successively and then prove 
\eqref{eq-decomposition-P-nonlinear-special-case4} by induction. 

By applying Lemma \ref{lemma-estimate-degree-|x|2-case2-H*} with $\ell = 0$ and $k$ replaced by $k-2$, we get 
\begin{align}\label{eq-estimate-P-k1-nonlinear-special-case4}
\|R_{k,1}\|_*\le \frac{1}{2(2k-4+n)}\|Q_{k-2}\|_* \le \frac{1}{2k}A_0 A^{k-2}.\end{align}
By applying Lemma \ref{lemma-estimate-degree-|x|2-case2-H*} again with $\ell = p+q$ and $k=i+j$, 
we have
$$\|\Delta_k^{-1}(R_{i,p}R_{j,q})\|_*\le
\big[\big(2(p+q)+2\big)\big(2(i+j)-2(p+q)+n\big)\big]^{-1}\|R_{i,p}R_{j,q}\|_*.$$
It is obvious that 
$$2(i+j) - 2(p+q) + n=(i+j)+(i+j-2(p+q))+n\ge i+j+2.$$
Hence, 
$$\|\Delta_k^{-1}(R_{i,p}R_{j,q})\|_*\le \big[2(p+q+1)(i+j+2)\big]^{-1}\|R_{i,p}\|_* \|R_{j,q}\|_*.$$
By a simple substitution in \eqref{eq-expression-P-kl-nonlinear-special-case4}, we obtain, 
for any $k\ge 2$ and $2\le\ell\le k/2$, 
\begin{align}\label{eq-estimate-P-kl-nonlinear-special-case4}
\|R_{k,\ell}\|_*\le\frac{1}{2k\ell}\sum_{i+j = k-2}\sum_{p+q=\ell-1} j(j-1)\|R_{i,p}\|_* \|R_{j,q}\|_*.\end{align} 
The first summation in the right-hand side of \eqref{eq-estimate-P-kl-nonlinear-special-case4} 
is over $(i,j)$ with $i+j=k-2$ and $i, j\ge 2$. 
We now insert a factor $i(i-1)$ and write, 
for any $k\ge 2$ and $2\le\ell\le k/2$, 
\begin{align}\label{eq-estimate-P-kl-stronger-nonlinear-special-case4}
k\ell\|R_{k,\ell}\|_*\le \sum_{i+j = k-2}\sum_{p+q=\ell-1} i(i-1)j(j-1)\|R_{i,p}\|_* \|R_{j,q}\|_*.\end{align}
The factor $k\ell$ suggests that we need to introduce a majorant analytic function of two variables. 

We first introduce a sequence $\{a_{k,\ell}\}$ of two indices. 
Set 
\begin{align}\label{eq-identity-a-kl-0-stronger-nonlinear-special-case4}
a_{k,\ell}=\begin{cases}
0&\text{for $k=0,1$ and $\ell\ge 0$},\\
0&\text{for $k\ge 0$ and $\ell= 0$},\\
\|\Delta_k^{-1} Q_{k-2}\|_*&\text{for $k\ge 2$ and $\ell= 1$},
\end{cases}
\end{align}
and, for any $k, \ell\ge 2$,
\begin{align}\label{eq-identity-a-kl-stronger-nonlinear-special-case4}
k\ell a_{k,\ell}=\sum_{i+j = k-2}\sum_{p+q=\ell-1} i(i-1)j(j-1)a_{i,p}a_{j,q}.\end{align}

We now introduce a function $\varphi=\varphi(s,t)$ formally defined by 
\begin{align}\label{eq-formal-expansion-nonlinear-special-case4} \varphi(s,t)=\sum_{k,\ell\ge 0}a_{k,\ell}s^kt^\ell.
\end{align}
Then, 
\begin{align}\label{eq-formal-boundary-condition-special-case4}
\varphi(0,t)=\partial_s\varphi(0,t)=0,\quad \varphi(s,0)=0,\end{align}
and 
\begin{align*}st\partial_{st}\varphi-\sum_{k\ge 2} k a_{k, 1}s^k t=s^2t(s^2\partial_{ss}\varphi)^2.\end{align*}
Set 
$$\beta(s)=\sum_{k=1}^\infty (k+1)a_{k+1, 1}s^k.$$ 
This is an analytic function near $s=0$ and $\beta(0)=0$ 
by \eqref{eq-estimate-P-k1-nonlinear-special-case4} and \eqref{eq-identity-a-kl-0-stronger-nonlinear-special-case4}. 
Hence, 
\begin{align}
\label{eq-majorant-ode-nonlinear-special-case4}
\partial_{st}\varphi= \beta(s) +s^5(\partial_{ss}\varphi)^2.
\end{align}
We point out that the derivation of 
\eqref{eq-majorant-ode-nonlinear-special-case4} is formal but the function $\beta$ 
is defined rigorously as an analytic function. 

We now derive \eqref{eq-exp-growth-P_k-nonlinear-special-case4}
from \eqref{eq-estimate-P-kl-stronger-nonlinear-special-case4}. 
We study the problem \eqref{eq-majorant-ode-nonlinear-special-case4}
and \eqref{eq-formal-boundary-condition-special-case4}, and prove there exists an analytic solution 
$\varphi=\varphi(s,t)$ near $(s,t)=(0,0)$. 
To this end, we consider 
\begin{align}\label{eq-majorant-pde-nonlinear-special-case4}
\partial_{t}\psi= \beta(s) + s^5(\partial_{s}\psi)^2, \quad \psi(s, 0)=0.
\end{align}
We treat \eqref{eq-majorant-pde-nonlinear-special-case4} as 
an initial value problem for $\psi=\psi(s,t)$, with $t$ interpreted as the time variable. 
By the Cauchy-Kovalevskaya theorem, \eqref{eq-majorant-pde-nonlinear-special-case4}
admits an analytic solution $\psi=\psi(s,t)$ near $(s,t)=(0,0)$. 
To recover an analytic solution $\varphi$ of \eqref{eq-majorant-ode-nonlinear-special-case4}
and \eqref{eq-formal-boundary-condition-special-case4}, we set 
\begin{align}\label{eq-definition-phi-nonlinear-special-case4}\varphi(s,t)=\int_0^s\psi(\tau, t)d\tau.\end{align}
Then, 
\begin{align}\label{eq-definition-psi-nonlinear-special-case4}\psi(s,t)=\partial_s\varphi(s,t).\end{align}
We should point out that \eqref{eq-majorant-pde-nonlinear-special-case4} is actually obtained from 
\eqref{eq-majorant-ode-nonlinear-special-case4} via \eqref{eq-definition-psi-nonlinear-special-case4}. 
It is obvious that $\varphi$ given by \eqref{eq-definition-phi-nonlinear-special-case4} satisfies 
\eqref{eq-majorant-ode-nonlinear-special-case4} and 
$$\varphi(0,t)=0,\quad \varphi(s,0)=0.$$
Next, by \eqref{eq-majorant-pde-nonlinear-special-case4} and $\beta(0)=0$, we get 
$\partial_{t}\psi(0,t)=0$. Then, $\psi(0,t)$ is a constant, independent of $t$, and hence 
$\psi(0,t)=\psi(0,0)=0$. Thus, 
$$\partial_s\varphi(0,t)=0.$$
Therefore, $\varphi$ is an analytic solution of \eqref{eq-majorant-ode-nonlinear-special-case4}
and \eqref{eq-formal-boundary-condition-special-case4}. 
Write $\varphi$ as in \eqref{eq-formal-expansion-nonlinear-special-case4}. 
Then, for any $k \ge 2$ and $\ell\ge 2$, \eqref{eq-identity-a-kl-stronger-nonlinear-special-case4} holds 
and
\begin{align}\label{eq-exp-growth-majorant-ode-sol-nonlinear-special-case4}
|a_{k,\ell}| \leq C_*^{k+\ell},
\end{align}
for some positive constant $C_*$.  
Inductively, we can prove that $a_{k,\ell}\ge 0$ for any $k\ge 2$ and $\ell\ge 1$.

We next prove that, for any $k \ge 2$ and $1\le\ell\le k/2$,
\begin{align}
\label{eq-P_k<=a_k-nonlinear-special-case4}
\|R_{k,\ell}\|_* \le a_{k,\ell}.
\end{align}
First, \eqref{eq-P_k<=a_k-nonlinear-special-case4} holds for $k\ge 2$ and $\ell=1$ 
by \eqref{eq-expression-P-k1-nonlinear-special-case4}
and \eqref{eq-identity-a-kl-0-stronger-nonlinear-special-case4}. 
Next, a simple induction  argument based on 
\eqref{eq-estimate-P-kl-stronger-nonlinear-special-case4}
and \eqref{eq-identity-a-kl-stronger-nonlinear-special-case4} yields 
\eqref{eq-P_k<=a_k-nonlinear-special-case4} for any $k \ge 2$ and $1\le\ell\le k/2$. 
Finally, by \eqref{eq-exp-growth-majorant-ode-sol-nonlinear-special-case4}
and \eqref{eq-P_k<=a_k-nonlinear-special-case4}, we obtain \eqref{eq-exp-growth-P_k-nonlinear-special-case4}.
\end{proof}

We now make several remarks. 
In the proof of Proposition \ref{prop-nonlinear-special-case4}, we insert a factor $i(i-1)$
in \eqref{eq-estimate-P-kl-nonlinear-special-case4}
to get \eqref{eq-estimate-P-kl-stronger-nonlinear-special-case4}. 
In this way, we obtain the equation of the majorant analytic function in the form 
\eqref{eq-majorant-ode-nonlinear-special-case4}. 
If we use the original \eqref{eq-estimate-P-kl-nonlinear-special-case4}, 
the equation of a majorant analytic function is given by 
\begin{align}\label{eq-majorant-ode-nonlinear-weaker-special-case4}
\partial_{st}\varphi=s^3\varphi\partial_{ss}\varphi+\beta(s).
\end{align}
It is not clear we can obtain an analytic solution of 
\eqref{eq-majorant-ode-nonlinear-weaker-special-case4}
and \eqref{eq-formal-boundary-condition-special-case4} 
by the Cauchy-Kovalevskaya theorem as easily as in the proof. 

The methods to prove Proposition \ref{prop-nonlinear-special-case3} 
and Proposition \ref{prop-nonlinear-special-case4} are different. 
In the proof of Proposition \ref{prop-nonlinear-special-case3}, 
we estimate each $P_k$ directly in terms of the homogeneous Taylor polynomials of $f$ of appropriate degrees. 
However, in the proof of Proposition \ref{prop-nonlinear-special-case4}, we need to decompose $P_k$ 
according to the degree of $|x|^2$ and estimate each term in the decomposition separately. 
The discussion before the statement of Proposition \ref{prop-nonlinear-special-case4}
demonstrates that such a decomposition is necessary.

The nonlinear elliptic equation leading to the estimate 
\eqref{eq-estimate-P-kl-stronger-nonlinear-special-case4}
is given by 
\begin{align}\label{eq-nonlinear-special-case4-squared} 
\Delta u=(x_ix_j\partial_{ij} u)^2+f.\end{align}
In any aspect, the nonlinear structure in \eqref{eq-nonlinear-special-case4-squared} 
is worse than that in \eqref{eq-nonlinear-special-case4}.
The proof of Proposition \ref{prop-nonlinear-special-case4} simply says that 
our method works for this \lq\lq worse\rq\rq\, equation
and suggests that the result should hold for more general nonlinear equations.

\section{Nonlinear Equations}
\label{sec-fully-nonlinear}

In this section, we study a class of nonlinear equations and systems. 
The results in this section play an important role in our study later.

Recall that $\langle \cdot,\cdot \rangle_{*}$ is the inner product defined in \eqref{eq-H^m-inner-product-with-constant}
and $\|\cdot\|_*$ is the corresponding norm. 
We will use the estimate \eqref{eq-Banach-algebra-property-without-constant} repeatedly.

Consider
\begin{align}
\label{eq-fully-nonlinear-equation}
\Delta u = F\big(x,u, x \nabla u, x^2 \nabla^2 u \big),
\end{align}
where $F$ is an analytic function in $B_1 \times \R \times \R^{n^2} \times \R^{n^4}$, 
$$x \nabla u:\,\text{the collection of all }x_i \partial_j u\text{ for }1 \leq i ,j\leq n,$$
and
$$x^2 \nabla^2 u:\,\text{the collection of all }x_i x_j \partial_{k\ell} u\text{ for }1 \leq i ,j,k,\ell\leq n.$$ 
We will study the existence and uniqueness of purely unharmonic solutions to \eqref{eq-fully-nonlinear-equation}.

First, we study the structure of the nonlinear functional $F(x,u,x\nabla u, x^2 \nabla^2 u)$ near $0$.
 Assume that $u$ is a $C^2$-function near $0$ with $u(0) = 0$. 
 Since $F$ is analytic near $0$, for $|x|$ sufficiently small, we can write
\begin{align}
\label{eq-Taylor-series-nonlinear-functional}
F(x,u,x\nabla u, x^2 \nabla^2 u) = \sum_{\alpha,l,\beta,\gamma} 
a_{\alpha,l,\beta,\gamma} x^\alpha u^l (x\nabla u)^\beta (x^2 \nabla^2 u)^\gamma,
\end{align}
for some constants $a_{\alpha,l,\beta,\gamma}$, 
where $\alpha \in \Z_+^n$, $l \ge 0$, $\beta \in \Z_+^{n^2}$, and $\gamma \in \Z_+^{n^4}$. 
Here, we write, for $\beta = (\beta_{ij}) \in \Z_+^{n^2}$,
$$(x\nabla u)^\beta=  \prod_{i,j=1}^n (x_i\partial_j u)^{\beta_{ij}},$$
and, for $\gamma = (\gamma_{pqrs}) \in \Z_+^{n^4}$,
$$(x^2 \nabla^2 u)^\gamma = \prod_{p,q,r,s = 1}^n (x_px_q\partial_{rs} u)^{\gamma_{pqrs}}.$$
By \eqref{eq-Taylor-series-nonlinear-functional}, we write
\begin{align}
\label{eq-Taylor-series-nonlinear-functional-1}
F(x,u,x\nabla u, x^2 \nabla^2 u) = f_0 + \sum_{i=1}^\infty \sum_{j=0}^\infty T_{ij}(u,\cdots, u)\quad\text{near }0,
\end{align}
where
$$f_0 = \sum_{\alpha \in \Z_+^n} a_{\alpha,0,0,0} x^\alpha,$$
and, for $i \ge 1$ and $j \ge 0$,
\begin{align}\label{eq-expression-T-ij}
T_{ij}(u,\cdots, u) = \sum_{l + |\beta| + |\gamma| = i} \sum_{|\alpha| = j} a_{\alpha,l,\beta,\gamma}
 x^\alpha u^l (x\nabla u)^\beta (x^2\nabla^2 u)^\gamma.\end{align}
Note that $T_{ij}$ is a multi-linear operator from $C^\infty(\R^n)^i$ into $C^\infty (\R^n)$ 
and maps $\mathscr{P}_{k_1}\times \cdots \times \mathscr{P}_{k_i}$
into $\mathscr{P}_{k_1 +\cdots + k_i + j}$, for any nonnegative integers $k_1,\cdots,k_i$. 

In the following, we prove the existence of purely unharmonic solutions to \eqref{eq-fully-nonlinear-equation} 
under an additional condition. It is motivated by our study of the minimal surface equation 
and the Monge-Ampère equation in exterior domains.

We first introduce the following terminology.

\begin{definition}
\label{def-operators-of-type-(i,j)}
Fix integers $i \geq 1$ and $j \ge 0$. 
A multi-linear operator $T:C^\infty(\R^n)^i \rightarrow C^\infty(\R^n)$ from $i$-copies of $C^\infty(\R^n)$ into $C^\infty(\R^n)$ 
is called to be {\it of type $(i,j)$ that preserves the degree of $|x|^2$} if
\begin{itemize}
    \item[(i)] $T:\mathscr{P}_{k_1}\times\cdots\times\mathscr{P}_{k_i} \rightarrow \mathscr{P}_{k_1 +\cdots + k_i + j}$
    for any $k_1, \cdots, k_i \ge 0$; 
    \item[(ii)] for any $u_1,\cdots,u_i \in C^\infty(\R^n)$ and $l_1, \cdots, l_i \ge 0$,
    $$T\big(|x|^{2l_1}u_1,\cdots,|x|^{2l_i} u_i \big) = |x|^{2(l_1 + \cdots + l_i)} B_{(l_1, \cdots, l_i)}(u_1, \cdots, u_i),$$
    for some operator $B_{(l_1, \cdots, l_i)}:C^\infty(\R^n)^i \rightarrow C^\infty(\R^n).$
\end{itemize}
\end{definition}

Note that, for fixed $i \ge 1$ and $j \ge 0$, the collection of all operators of type $(i,j)$ 
that preserve the degree of $|x|^2$ forms a vector space over $\R$. 
It is clear that Definition \ref{def-operators-of-type-(i,j)} generalizes Definition \ref{def-operators-of-type-i} 
from operators acting on one function to operators acting on multiple functions. 
An operator of type $(1,j)$ that preserves the degree of $|x|^2$ is an operator of type $j$ 
that preserves the degree of $|x|^2$ in Section \ref{sec-linear-eqs}.

In the next two sections, we will discuss examples of operators satisfying Definition \ref{def-operators-of-type-(i,j)}. 
Refer to Lemmas \ref{lemma-operators-in-MSE}, \ref{lemma-Monge-Ampere-operators-n>=3}, 
and \ref{lemma-Monge-Ampere-operators-n=2}. 
We next present methods to construct new operators of type $(i,j)$ that preserve 
the degree of $|x|^2$ from existing ones. 

\begin{lemma}
\label{lemma-constructing-new-operators}
Let $T$ be an operator of type $(i,j)$ and $T'$ be an operator of type $(k,\ell)$ 
that preserve the degree of $|x|^2$, and $Q \in \mathscr{P}_m$ be a homogeneous polynomial of degree $m$.
\begin{itemize}
    \item[(i)] For any fixed $1 \leq s \leq i$, the operator given by
    $$(u_1,\cdots, u_{s-1}, u_{s+1},\cdots,u_{i}) \mapsto T(u_1,\cdots,u_{s-1}, Q, u_{s+1}, \cdots,u_i)$$
    is an operator of type $(i-1,j+m)$ that preserves the degree of $|x|^2$.
    \item[(ii)] The operator given by
    $$(u_1,\cdots,u_i) \mapsto Q\cdot T(u_1,\cdots,u_i)$$
    is an operator of type $(i,j+m)$ that preserves the degree of $|x|^2$.
    \item[(iii)] The operator given by
    $$(u_1,\cdots,u_i,v_1,\cdots,v_k)\mapsto T(u_1,\cdots,u_i)T'(v_1,\cdots,v_k)$$
    is an operator of type $(i+k,j+\ell)$ that preserves the degree of $|x|^2$.
    \item[(iv)] For any permutation $\sigma$ of $\{1,\cdots,i\}$, the operator given by
    $$(u_1,\cdots,u_i) \mapsto T(u_{\sigma_1},\cdots, u_{\sigma_i})$$
    is an operator of type $(i,j)$ that preserves the degree of $|x|^2$.
\end{itemize}
\end{lemma}

In the next theorem, we prove our main existence result if $F$ in \eqref{eq-fully-nonlinear-equation} has a special structure.

\begin{theorem}
\label{thrm-existence-fully-nonlinear}
Let $F$ be an analytic function in $B_1 \times\R \times\R^{n^2} \times \R^{n^4}$ such that 
\eqref{eq-Taylor-series-nonlinear-functional-1} holds near $0$,
for some analytic function $f_0$ near $0 \in \R^n$ and operators $T_{ij}$ of type $(i,j)$ 
that preserve the degree of $|x|^2$, for $i \ge 1$ and $j \ge 0$. 
Then, there exists a unique analytic function $u_*$ in $B_R$ for some $R \in (0,1)$ such that $u = |x|^2 u_*$ 
is a solution to \eqref{eq-fully-nonlinear-equation} in $B_R$. 
\end{theorem}

We point out that the key assumption in Theorem \ref{thrm-existence-fully-nonlinear} 
is the preservation of the degree of $|x|^2$ by $T_{ij}$.

\begin{proof}
Let $\{Q_k\}_{k\ge 0}$ be the sequence of the homogeneous Taylor polynomials of $f_0$ at $0$, 
with $Q_k \in \mathscr{P}_k$. Assume that $u = |x|^2 u_*$ is a solution to \eqref{eq-fully-nonlinear-equation} 
near $0$ with $F$ given by \eqref{eq-Taylor-series-nonlinear-functional-1}. 
Let $\{P_k\}_{k\ge 2}$ be the sequence of the homogeneous Taylor polynomials of $u$ at $0$, with $P_k \in \mathscr{H}_k^\perp$. 
It is straightforward to check that, for any $k \ge 2$,
\begin{align}
\label{eq-recursive-equation-P_k-fully-nonlinear}
\Delta P_k = Q_{k-2} + \sum_{m+j = k-2} \sum_{i=1}^{[m/2]} \sum_{k_1 + \cdots + k_i = m} T_{ij}(P_{k_1},\cdots, P_{k_i}).
\end{align}
By $P_k \in \mathscr{H}_k^\perp$ for any $k \ge 2$, $P_k$ is uniquely determined by $Q_{k-2}$ and $P_2,\cdots, P_{k-2}$. 
Hence, any solution $u = |x|^2 u_*$ of \eqref{eq-fully-nonlinear-equation} for some analytic function $u_*$ 
near $0$ with $F$ given by \eqref{eq-Taylor-series-nonlinear-functional-1} has the same Taylor series at $0$. 
This proves the uniqueness. We point out that the uniqueness does not require the preservation of the degree of $|x|^2$ by $ T_{ij}$.

Set $P_0 = P_1 = 0$. 
By using Lemma \ref{lemma-Laplace-estimate-poly-H*} successively for $k=2, 3, \cdots$, we conclude 
that there is a unique solution $P_{k}\in \mathscr H_{k}^\perp$ of 
\eqref{eq-recursive-equation-P_k-fully-nonlinear}. 
It suffices to prove that, for any $k \geq 2$,
\begin{align}
\label{eq-exp-growth-P_k-fully-nonlinear}
\|P_k\|_* \le C_*^k,
\end{align}
for some positive constant $C_*$.
Then, the existence follows from Corollary \ref{cor-constructing-analytic-functions}.

Recall that $T_{ij}$ is as in \eqref{eq-Taylor-series-nonlinear-functional-1} and \eqref{eq-expression-T-ij}. 
Since $F$ is analytic, we have, for any $i \ge 1$, $j \ge 0$, $k_1,\cdots, k_i \ge 2$, 
and $P_{k_1} \in \mathscr{P}_{k_1},\cdots, P_{k_i} \in \mathscr{P}_{k_i}$,
\begin{align}
\label{eq-estimate-for-T_ij-fully-nonlinear}
\|T_{ij}(P_{k_1},\cdots, P_{k_i})\|_* \leq A^i B_0 B^j k_1(k_1-1)\|P_{k_1}\|_*\cdots k_i(k_i-1)\|P_{k_i}\|_*,
\end{align}
for some constants $A, B_0, B>0$, depending only on $n$ and $F$.
By \eqref{eq-recursive-equation-P_k-fully-nonlinear}, we have
\begin{align}
\label{eq-recursive-zeroth-iteration-fully-nonlinear}
P_k = \Delta_k^{-1} Q_{k-2} +  \Delta_k^{-1}\sum_{m+j = k-2} \sum_{i=1}^{[m/2]} \sum_{k_1 + \cdots + k_i = m}  T_{ij}(P_{k_1},\cdots, P_{k_i}).
\end{align}
We now further decompose $P_k$ as
\begin{align}
\label{eq-further-decomp-of-P_k-fully-nonlinear}
P_k = \sum_{1 \leq \ell \leq k/2} R_{k,\ell},
\end{align}
for some $R_{k,\ell} \in |x|^{2\ell} \mathscr{P}_{k-2\ell}$, for $1 \leq \ell \leq k/2$.
To this end, define 
$$R_{0,\ell} = R_{1,\ell} = 0\text{ for any }\ell \geq 0,$$
$$R_{k,0} = 0\text{ for any }k \ge 0,$$
and $R_{k,\ell}$ inductively by, for $k \ge 2$ and $\ell = 1$,
\begin{align}
\label{eq-def-R_k,1-fully-nonlinear}
R_{k,1} = \Delta_k^{-1} Q_{k-2},
\end{align}
and, for $2 \leq \ell \leq k/2$,
\begin{align}
\label{eq-def-R_k,ell-fully-nonlinear}
R_{k,\ell} = \Delta_k^{-1}\sum_{m+j = k-2} \sum_{i=1}^\infty 
\sum_{\substack{k_1+\cdots + k_i = m\\
\ell_1 + \cdots + \ell_i = \ell-1}}  
 T_{ij}(R_{k_1,\ell_1},\cdots, R_{k_i,\ell_i}),
\end{align}
and $R_{k,\ell} = 0$ for any $\ell > k/2$. We point out that if $(k_1,\ell_1) +\cdots + (k_i,\ell_i) = (m,\ell-1)$ 
with $k_1,\cdots,k_i \ge 2$ and $\ell_1,\cdots,\ell_i \ge 1$, then $i \le \min \{[m/2], \ell-1\}$. 
Hence, the second sum in \eqref{eq-def-R_k,ell-fully-nonlinear} is over $1 \le i \le \min \{[m/2], \ell-1\}$. 
There are only finitely many nonzero terms in the right-hand side of \eqref{eq-def-R_k,ell-fully-nonlinear}.
We claim that \eqref{eq-further-decomp-of-P_k-fully-nonlinear} holds for such $R_{k,\ell}$. 
First, since $T_{ij}$ is an operator of type $(i,j)$ that preserves the degree of $|x|^2$, for each $i \ge 1$ 
and $j\ge 0$, it is straightforward to check that $R_{k,\ell} \in |x|^{2\ell} \mathscr{P}_{k-2\ell}$ 
for any $k \ge 2$ and $ \ell \ge 1$. By induction and \eqref{eq-recursive-zeroth-iteration-fully-nonlinear}, we have
\begin{align*}
P_k &= \Delta_k^{-1} Q_{k-2} +  \Delta_k^{-1}\sum_{m+j = k-2} \sum_{i=1}^{\infty} \sum_{k_1 + \cdots + k_i = m}  
T_{ij}(P_{k_1},\cdots, P_{k_i})\\
& = \Delta_k^{-1} Q_{k-2} +  \Delta_k^{-1}\sum_{m+j = k-2} \sum_{i=1}^{\infty} \sum_{k_1 + \cdots + k_i = m}  
\sum_{\ell_1,\cdots,\ell_i \ge 0}  T_{ij}(R_{k_1,\ell_1},\cdots, R_{k_i,\ell_i}).
\end{align*}
Now, we write
\begin{align}
\label{eq-summation-over-l_1,...,l_i}
\sum_{\ell_1,\cdots,\ell_i \ge 0}  T_{ij}(R_{k_1,\ell_1},\cdots, R_{k_i,\ell_i}) 
= \sum_{\ell=2}^\infty \sum_{\ell_1+\cdots + \ell_i = \ell - 1}  T_{ij}(R_{k_1,\ell_1},\cdots, R_{k_i,\ell_i}),
\end{align}
where we used $R_{k,0} = 0$ for any $k \ge 0$. Note that $R_{k_s,\ell_s} = 0$ if $\ell_s> k_s/2$. 
Hence, for any fixed $k\ge 2$, the sum in the right-hand side of \eqref{eq-summation-over-l_1,...,l_i} 
is over finitely many $\ell$ and a finite collection of $\ell_1,\cdots,\ell_i\ge 1$ with
$$2 \le \ell = \ell_1+\cdots + \ell_i+1 \le \frac{1}{2}(k_1 + \cdots + k_i)+1 \le \frac{1}{2}k.$$ 
By switching the order of summation, we get a sum over $2 \le \ell \le k/2$ of
\begin{align*}
\Delta_k^{-1}\sum_{m+j = k-2} \sum_{i=1}^{\infty} \sum_{\substack{k_1+\cdots + k_i = m\\
\ell_1 + \cdots + \ell_i = \ell-1}} T_{ij}(R_{k_1,\ell_1},\cdots, R_{k_i,\ell_i}).
\end{align*}
This is simply $R_{k,\ell}$ by \eqref{eq-def-R_k,ell-fully-nonlinear}, for $2 \le \ell \le k/2$. 
Together with \eqref{eq-def-R_k,1-fully-nonlinear}, we obtain \eqref{eq-further-decomp-of-P_k-fully-nonlinear}. 

Now, fix $k \ge 2$ and $\ell \ge 2$ with $\ell \leq k/2$. By Lemma \ref{lemma-estimate-degree-|x|2-case2-H*} 
and $R_{k,\ell} \in |x|^{2\ell}\mathscr{P}_{k-2\ell}$, we have
$$\|R_{k,\ell}\|_* \le \frac{1}{2\ell(2k-2\ell + n-2)}\|\Delta R_{k,\ell}\|_* \leq \frac{1}{2k\ell} \|\Delta R_{k,\ell}\|_*,$$
where we used
$$2k-2\ell+n-2 = k+(k-2\ell)+(n-2) \ge k.$$
By \eqref{eq-estimate-for-T_ij-fully-nonlinear} and \eqref{eq-def-R_k,ell-fully-nonlinear}, 
we get
\begin{equation}
\label{eq-estimate-for-R_k,l-fully-nonlinear}
\begin{aligned}
k\ell \|R_{k,\ell}\|_*
\leq \frac{1}{2} B_0 \sum_{m+j = k-2} \sum_{i=1}^{\infty} \sum_{\substack{k_1+\cdots + k_i = m\\
\ell_1 + \cdots + \ell_i = \ell-1}}
A^i B^j \prod_{s=1}^i k_s(k_s-1) \|R_{k_s,\ell_s}\|_*.
\end{aligned}
\end{equation}
We point out that \eqref{eq-estimate-for-T_ij-fully-nonlinear} holds for homogeneous polynomials of degree at least $2$. 
Since $R_{0,\ell} = R_{1,\ell} = 0$ for any $\ell \ge 0$, \eqref{eq-estimate-for-R_k,l-fully-nonlinear} is valid. 
In addition, \eqref{eq-estimate-for-R_k,l-fully-nonlinear} obviously holds for $k\ge 2$ and $\ell \ge 2$ with $\ell > k/2$, 
since the left-hand side is $0$.

We proceed formally. Define 
$$A_{0,\ell} = A_{1,\ell} = 0\text{ for any }\ell \ge 0,$$
$$A_{k,0} = 0\text{ for any }k \ge 0,$$
and, for any $k \ge 2$,
$$A_{k,1} = \|R_{k,1}\|_*,$$
and, for any $\ell \ge 2$,
\begin{align}
\label{eq-def-A_k,l}
k\ell A_{k,\ell}
= \frac{1}{2}B_0\sum_{m+j = k-2} \sum_{i=1}^{\infty} \sum_{\substack{k_1+\cdots + k_i = m\\
\ell_1 + \cdots + \ell_i = \ell-1}}
A^i B^j \prod_{s=1}^i k_s(k_s-1) A_{k_s,\ell_s}.
\end{align}
Set
$$\varphi(s,t) = \sum_{k\ge 2,\ell\ge 1} A_{k,\ell}s^k t^\ell\quad\text{near }0 \in \R^2.$$
Then, \eqref{eq-def-A_k,l} implies
\begin{align*}
&st\partial_{st} \varphi - \sum_{k\ge 2} k\|R_{k,1}\|_*s^k t 
= s^2 t \frac{1}{2} B_0  \Big( \sum_{j=0}^\infty B^j s^j \Big) \Big( \sum_{i=1}^{\infty} A^i \big(s^2\partial_{ss} \varphi \big)^i\Big),\\
&\varphi(s,0) = \varphi(0,t) = \partial_s \varphi(0,t) = 0,
\end{align*}
or equivalently,
\begin{equation}
\begin{aligned}
\label{eq-majorant-pde-second-order-fully-nonlinear-1}
&\partial_{st} \varphi = \beta(s) +   \frac{1}{2} B_0 b(s)  a \big(s^2\partial_{ss} \varphi \big),\\
&\varphi(s,0) = \varphi(0,t) = \partial_s \varphi(0,t) = 0,
\end{aligned}
\end{equation}
where
$$\beta(s) = \sum_{k=1}^\infty(k+1) \|R_{k+1,1}\|_*s^k,\quad b(s) = \sum_{j=0}^\infty B^j s^{j+1},
\quad\text{and}\quad a(\tau) = \sum_{i=1}^\infty A^i \tau^i.$$
Note that,
by Lemma \ref{lemma-estimate-degree-|x|2-case2-H*} with $\ell = 0$ and $k$ replaced by $k-2$,
for any $k \ge 2$,
$$\|R_{k,1}\|_*=  \|\Delta_k^{-1} Q_{k-2}\|_* \le \frac{1}{2(2k-4+n)} \|Q_{k-2}\|_* \le \frac{1}{2k}C^k,$$
for some constant $C>0$, depending only on $f_0$. 
Hence, $\beta$, $a$, and $b$ are analytic functions near $0$. We point out that the derivation 
of \eqref{eq-majorant-pde-second-order-fully-nonlinear-1} is formal, but $\beta$, $a$, 
and $b$ are defined rigorously as analytic functions near $0$.

We now prove that \eqref{eq-majorant-pde-second-order-fully-nonlinear-1} 
admits an analytic solution $\varphi = \varphi(s,t)$ near $0 \in \R^2$. Consider
\begin{equation}
\label{eq-majorant-pde-first-order-fully-nonlinear}
\begin{aligned}
& \partial_t \psi = \beta(s) + \frac{1}{2} B_0 b(s) a\big(s^2 \partial_s \psi \big),\\
& \psi(s,0) = 0.
\end{aligned}
\end{equation}
By the Cauchy-Kovalevskaya theorem, \eqref{eq-majorant-pde-first-order-fully-nonlinear} 
has a unique analytic solution $\psi$ in $(-\epsilon,\epsilon)\times (-\epsilon,\epsilon)$ for some $\epsilon > 0$. Set
$$\varphi(s,t) = \int_{0}^s \psi(\tau,t) d\tau.$$
It is easy to see that $\varphi$ satisfies the equation in \eqref{eq-majorant-pde-second-order-fully-nonlinear-1}
and, for any $|s| <\epsilon$ and $|t| < \epsilon$,
$$\varphi(0,t) = 0\quad\text{and}\quad \varphi(s,0) = \int_0^s \psi(\tau,0)d\tau = 0.$$
Next, by \eqref{eq-majorant-pde-first-order-fully-nonlinear} and $\beta(0) = a(0)=0$, 
we get, for any $|t| < \epsilon$, 
$\partial_{t}\psi(0,t)=0$. Then, $\psi(0,t)$ is a constant, independent of $t$, and hence 
$\psi(0,t)=\psi(0,0)=0$. Thus, 
$$\partial_s\varphi(0,t)=0.$$
Therefore, $\varphi$ is an analytic solution to \eqref{eq-majorant-pde-second-order-fully-nonlinear-1} 
in $(-\epsilon,\epsilon)\times(-\epsilon,\epsilon)$.

By \eqref{eq-estimate-for-R_k,l-fully-nonlinear}, \eqref{eq-def-A_k,l}, 
and $A_{k,1} = \|R_{k,1}\|_*$ for any $k\ge 2$, a simple induction yields, for any $k \ge 2$ and $\ell \ge 1$,
$$\|R_{k,\ell}\|_* \le A_{k,\ell}.$$
Since $\varphi$ is analytic near $0$, we have
$$A_{k,\ell} \le D^{k+\ell},$$
for some constant $D>0$. Therefore, for any $k \ge 2$,
$$\|P_k\|_* \le \sum_{1 \le \ell \le k/2} \|R_{k,\ell}\|_* \le \sum_{1 \le \ell \le k/2} D^{k+\ell} \le \frac{k}{2} D^{2k}.$$
This implies \eqref{eq-exp-growth-P_k-fully-nonlinear} and finishes the proof of the theorem.
\end{proof}

In the rest of this section, we generalize Theorem \ref{thrm-existence-fully-nonlinear} to elliptic systems. 
This is motivated by our study of the Monge-Ampère equation in exterior domains when the dimension $n$ is odd. 
Most of the proofs are the same. We will only provide an outline.

Fix a positive integer $m$. By Lemma \ref{lemma-Laplace-estimate-poly-H*}, for any $Q_k \in \mathscr{P}_k$, 
there exists a unique $P_{k+2} \in \mathscr{H}_{k+2}^\perp$ such that
$$\Delta \big( |y|^m P_{k+2} \big) = |y|^m Q_k,$$
and
$$\|P_{k+2}\|_*
\le \big[(m+2)(m+2k+n)\big]^{-1} \|Q_{k}\|_*.$$
Moreover, if $Q_k
\in |x|^{2\ell}\mathscr{P}_{k-2\ell}$, then
$$\|P_{k+2}\|_*
\le \big[(m+2\ell+2)(m+2k-2\ell+n)\big]^{-1} \|Q_{k}\|_*.$$

Consider
\begin{align}
\label{eq-fully-nonlinear-system}
\begin{split}
\Delta \widetilde{u} &= \widetilde{F}\big(x,\widetilde{u}, \widehat{u}, x \nabla \widetilde{u}, x \nabla \widehat{u}, 
x^2 \nabla^2 \widetilde{u}, x^2 \nabla^2 \widehat{u} \big),\\
\Delta (|x|^m \widehat{u} ) &= |x|^m \widehat{F}\big(x,\widetilde{u}, \widehat{u}, x \nabla \widetilde{u}, x \nabla \widehat{u}, 
x^2 \nabla^2 \widetilde{u}, x^2 \nabla^2 \widehat{u} \big),
\end{split}
\end{align}
where $\widetilde{F}$ and $\widehat{F}$ are analytic functions in $B_1 \times \R^2 \times \R^{2n^2} \times \R^{2n^4}$. 
We will study the existence and uniqueness of purely unharmonic solutions to \eqref{eq-fully-nonlinear-system}.

We first study the structure of the nonlinear functionals in this setting. Let $F$ be an analytic function 
in $B_1 \times \R^2 \times \R^{2n^2} \times \R^{2n^4}$ and $\widetilde u, \widehat u$ be $C^2$-functions 
near $0$ with $\widetilde u (0) = \widehat u(0) = 0$. Then, for $|x|$ sufficiently small,
\begin{equation}
\label{eq-Taylor-series-nonlinear-functional-2}
\begin{aligned}
&F\big(x,\widetilde{u}, \widehat{u}, x \nabla \widetilde{u}, x \nabla \widehat{u}, 
x^2 \nabla^2 \widetilde{u}, x^2 \nabla^2 \widehat{u} \big)\\
&\qquad=\, f_0 + \sum_{i=1}^\infty \sum_{j=0}^\infty  T^{(i,0)}_{ij}(\widetilde u,\cdots, \widetilde u)
+ \sum_{i=1}^\infty \sum_{j=0}^\infty  T^{(0,i)}_{ij}(\widehat u,\cdots, \widehat u)\\
&\qquad\qquad+ \sum_{i=2}^\infty \sum_{j=0}^\infty 
\sum_{\substack{\widetilde \iota + \widehat\iota = i\\ \widetilde\iota, \widehat\iota \ge 1}} 
T^{(\widetilde \iota,\widehat\iota)}_{ij}(\widetilde u,\cdots, \widetilde u, \widehat u,\cdots \widehat u),
\end{aligned}
\end{equation}
for some analytic function $f_0$ near $0$ and multi-linear operators 
$T^{(\widetilde \iota, \widehat \iota)}_{ij} : C^\infty (\R^n)^i \rightarrow C^\infty(\R^n)$ for $i \ge 1$, $j \ge 0$, 
and $\widetilde \iota, \widehat \iota \ge 0$ with $\widetilde \iota + \widehat \iota = i$. 
In the last summation, the expression 
$T^{(\widetilde \iota,\widehat\iota)}_{ij}(\widetilde u,\cdots, \widetilde u, \widehat u,\cdots ,\widehat u)$ 
consists of $\widetilde \iota$-copies of $\widetilde u$ and $\widehat \iota$-copies of $\widehat u$. 
Compare \eqref{eq-Taylor-series-nonlinear-functional-2} with \eqref{eq-Taylor-series-nonlinear-functional-1}. 

The following result is the counterpart of Theorem \ref{thrm-existence-fully-nonlinear} for \eqref{eq-fully-nonlinear-system}.

\begin{theorem}
\label{thrm-existence-fully-nonlinear-system}
Let $\widetilde F,\widehat F$ be analytic functions in $B_1 \times\R^2 \times\R^{2n^2} \times \R^{2n^4}$ 
such that \eqref{eq-Taylor-series-nonlinear-functional-2} holds for $\widetilde F$ and $\widehat F$ near $0$ 
with $(F,f_0,T^{(\widetilde \iota,\widehat\iota)}_{ij})$ replaced by 
$(\widetilde F,\widetilde f_0, \widetilde T^{(\widetilde \iota,\widehat\iota)}_{ij})$ and 
$(\widehat F,\widehat f_0,\widehat T^{(\widetilde \iota,\widehat\iota)}_{ij})$, respectively, 
for some analytic functions $\widetilde f_0,\widehat f_0$ near $0$ and operators 
$\widetilde T^{(\widetilde\iota, \widehat\iota)}_{ij}, \widehat T^{(\widetilde\iota, \widehat\iota)}_{ij}$ of type $(i,j)$ 
that preserve the degree of $|x|^2$.
Then, there exists a unique pair of analytic functions $(\widetilde u_*, \widehat u_*)$ in $B_R$ for some 
$R \in (0,1)$ such that $(\widetilde u, \widehat u) = (|x|^2 \widetilde u_*, |x|^2 \widehat u_*)$ 
is a pair of solutions to \eqref{eq-fully-nonlinear-system} in $B_R$. 
\end{theorem}

We point out that the key assumption in Theorem \ref{thrm-existence-fully-nonlinear-system} 
is the preservation of the degree of $|x|^2$ by $\widetilde T^{(\widetilde\iota, \widehat\iota)}_{ij}$ 
and $\widehat T^{(\widetilde\iota, \widehat\iota)}_{ij}$.

\begin{proof}
The uniqueness can be proved by the same argument as in the proof of Theorem \ref{thrm-existence-fully-nonlinear}. 
We point out that the uniqueness does not require the preservation of the degree of $|x|^2$ 
by $\widetilde T^{(\widetilde\iota, \widehat\iota)}_{ij}$ and $\widehat T^{(\widetilde\iota, \widehat\iota)}_{ij}$. 
We now prove the existence.
By the analyticity of $\widetilde F$ and $\widehat F$, we have, for any $i \ge 1$, $j \ge 0$, 
$\widetilde \iota, \widehat \iota \ge 0$ with $\widetilde \iota + \widehat \iota = i$, $k_1,\cdots, k_i \ge 2$, 
and $P_{k_1} \in \mathscr{P}_{k_1},\cdots, P_{k_i} \in \mathscr{P}_{k_i}$, 
\begin{align*}
&\|\widetilde T^{(\widetilde\iota, \widehat\iota)}_{ij}(P_{k_1},\cdots, P_{k_i})\|_*,\,
\|\widehat T^{(\widetilde\iota, \widehat\iota)}_{ij}(P_{k_1},\cdots, P_{k_i})\|_*\\
&\qquad\leq  A^i B_0 B^j k_1(k_1-1)\|P_{k_1}\|_*\cdots k_i(k_i-1)\|P_{k_i}\|_*,
\end{align*}
for some constants $A,B_0,B>0$ depending only on $n$, $\widetilde F$, and $\widehat F$.  
By an argument similar to that of the proof of Theorem \ref{thrm-existence-fully-nonlinear}, 
it suffices to solve a second order system
\begin{equation}
\label{eq-majorant-PDE-2nd-order-system}
\begin{aligned}
&\partial_{st} \widetilde \varphi = G(s,\partial_{ss} \widetilde \varphi, \partial_{ss}\widehat \varphi),\\
&\partial_{st} \widehat \varphi = G(s,\partial_{ss} \widetilde \varphi, \partial_{ss}\widehat \varphi),\\
&\widetilde \varphi(s,0) = \widetilde \varphi(0,t) 
= \partial_s \widetilde \varphi(0,t) =\widehat \varphi(s,0) = \widehat \varphi(0,t) = \partial_s \widehat \varphi(0,t) = 0,
\end{aligned}
\end{equation}
where $G$ is given by
\begin{align*}
G(s,\partial_{ss} \widetilde \varphi, \partial_{ss}\widehat \varphi) &=\beta (s) + 
\frac{1}{2} B_0 b(s)\big[ a \big( s^2\partial_{ss} \widetilde \varphi \big) +  a\big(s^2 \partial_{ss} \widehat \varphi \big) \big]\\
&\qquad+ \frac{1}{2} B_0 b(s) a\big( s^2\partial_{ss} \widetilde \varphi\big) a\big(s^2\partial_{ss} \widehat \varphi \big),
\end{align*}
for some analytic functions $\beta(s)$ depending only on $\widetilde f_0$ and $\widehat f_0$, $b(s)$ depending only on $B$, 
and $a(\tau)$ depending only on $A$. In addition, $\beta(0) = b(0) = a(0) = 0$. 
Then, we consider the associated first order system
\begin{equation}
\label{eq-majorant-PDE-1st-order-system}
\begin{aligned}
&\partial_{t} \widetilde \psi = G(s,\partial_{s} \widetilde \psi, \partial_{s}\widehat \psi),\\
&\partial_{t} \widehat \psi = G(s,\partial_{s} \widetilde \psi, \partial_{s}\widehat \psi),\\
&\widetilde \psi(s,0) = \widehat \psi(s,0) = 0.
\end{aligned}
\end{equation}
By the Cauchy-Kovalevskaya theorem, \eqref{eq-majorant-PDE-1st-order-system} 
admits a unique pair of solutions $(\widetilde \psi, \widehat \psi)$ 
in $(-\epsilon,\epsilon)\times (-\epsilon,\epsilon)$ for some $\epsilon > 0$. Then,
$$\widetilde \varphi = \int_0^s \widetilde \psi(\tau,t)d\tau
\quad\text{and}\quad \widehat \varphi = \int_0^s \widehat \psi(\tau,t)d\tau$$
are analytic solutions to \eqref{eq-majorant-PDE-2nd-order-system}. This finishes the proof of the theorem.
\end{proof}

\section{Minimal Surface Equation in exterior domains}
\label{sec-MSE}

In this section, we study the minimal surface equation in $\R^n \setminus \overline{B}_1$. 
Let $u\in C^2(\R^n\setminus \overline{B}_1)$ be a solution of 
\begin{equation}\label{eq-MSE} 
\Delta u-\frac{u_iu_j}{1+|\nabla u|^2}u_{ij}=0\quad\text{in }\mathbb R^n\setminus\overline{B}_1.\end{equation}
Han and Wang \cite{Han-Wang2024} proved the following result.

\begin{theorem}
Let $n \ge 2$ and $u \in C^2(\R^n \setminus \overline{B}_1)$ be a solution of \eqref{eq-MSE}, 
with $\nabla u(x)\to 0$ as $x\to\infty$.

$\mathrm{(i)}$ For $n =2$, there exist a constant $a$ and a function $v\in C^\infty(B_1)$ such that, 
for any $x\in \R^2\setminus\overline{B}_1$,  
\begin{equation}\label{eq-MSE-expansion-identity-dim-2}
u(x)=a\log|x|+v\Big(\frac{x}{|x|^2}\Big).\end{equation}

$\mathrm{(ii)}$ For $n \geq 3$, there exist a constant $a$ and a function $v\in C^\infty(B_1)$ such that, 
for any $x\in \mathbb R^n\setminus\overline{B}_1$,  
\begin{equation}\label{eq-MSE-expansion-identity-dim-higher}
u(x)=a+|x|^{2-n}v\Big(\frac{x}{|x|^2}\Big).\end{equation}
\end{theorem}

The function $v$ in \eqref{eq-MSE-expansion-identity-dim-2} and \eqref{eq-MSE-expansion-identity-dim-higher} 
is originally defined in the punctured unit ball $B_1\setminus\{0\}$ and satisfies a nonlinear elliptic equation. 
The proof of the smoothness of $v$ in the entire $B_1$ 
is based on a removable singularity result for $v$. 
We now record the equation satisfied by $v$. Refer to \cite{Han-Wang2024} for details. 

Consider the Kelvin transform 
$$y=\frac{x}{|x|^2}.$$
From now on, we switch the notations for variables from $x\in \mathbb R^n\setminus\overline{B}_1$
to $y\in B_1$.  
The function $v=v(y)$ in \eqref{eq-MSE-expansion-identity-dim-2} and \eqref{eq-MSE-expansion-identity-dim-higher} satisfies 
\begin{equation}\label{eq-MSE-v}
\Delta v+a_{ij}(y,v,\nabla v)\partial_{ij} v=f(y,v,\nabla v)\quad\text{in }B_1,
\end{equation}
where $a_{ij}$ and $f$ are given by the following expressions: 
for $n= 2$, 
\begin{align*}
\begin{split}
a_{ij}(y,v,\nabla v)&=|y|^2\big[a^2-2a\langle y, \nabla v\rangle+|y|^2|\nabla v|^2\big]\delta_{ij}\\
&\qquad -\big[a^2y_iy_j-a|y|^2(y_i\partial_j v+y_j\partial_i v)+|y|^4\partial_i v\partial_j v\big],
\end{split}\end{align*}
and 
\begin{align*}
\begin{split}
f(y,v,\nabla v)&=-a^3+4a^2\langle y, \nabla v\rangle
-2a\langle y, \nabla v\rangle^2\\
&\qquad-3a|y|^2|\nabla v|^2
+2|y|^2\langle y, \nabla v\rangle|\nabla v|^2,
\end{split}\end{align*}
for the constant $a$ as in \eqref{eq-MSE-expansion-identity-dim-2}, and, for $n\ge 3$, 
\begin{align*}
\begin{split}
a_{ij}(y,v,\nabla v)&=|y|^{2n-2}\big[(n-2)^2v^2+2(n-2)v\langle y, \nabla v\rangle+|y|^2|\nabla v|^2\big]\delta_{ij}\\
&\quad -|y|^{2n-4}\big[(n-2)^2y_iy_jv^2+(n-2)|y|^2v(y_i\partial_j v+y_j\partial_i v)+|y|^4\partial_i v\partial_j v\big],
\end{split}\end{align*}
and 
\begin{align*}
\begin{split}
f(y,v,\nabla v)&=|y|^{2n-4}\big[(n-2)^3(n-1)v^3+4(n-2)^2(n-1)v^2\langle y, \nabla v\rangle\\
&\qquad +(n-2)(3n-4)v\langle y, \nabla v\rangle^2+(n-2)(2n-1)v|y|^2|\nabla v|^2\\
&\qquad +2(n-1)|y|^2\langle y, \nabla v\rangle|\nabla v|^2\big].
\end{split}\end{align*}
Note that \eqref{eq-MSE-v} is a quasilinear equation and that 
$a_{ij}$ and $f$ are polynomials in terms of $y$, $v$, and $y\nabla v$. Recall that, for a smooth function $v$, we write 
$$y \nabla v:\,\text{the collection of all }y_i \partial_j v\text{ for }1 \leq i ,j\leq n,$$
and
$$y^2 \nabla^2 v:\,\text{the collection of all }y_i y_j \partial_{k\ell} v\text{ for }1 \leq i ,j,k,\ell\leq n.$$

By a simple rearrangement, we write \eqref{eq-MSE-v} as
\begin{align}\label{eq-nonlinear-Poission-MSE-fully} 
\Delta v=F\quad\text{in }B_1,
\end{align} 
where $F$ 
has the following form: for $n=2$,   
\begin{equation}\label{eq-MSE-dim2-fully}\begin{split}
F&=-|y|^2\big[a^2-2a\langle y, \nabla v\rangle+|y|^2|\nabla v|^2\big]\Delta v\\
&\qquad +\big[a^2y_iy_j-a|y|^2(y_i\partial_j v+y_j\partial_i v)+|y|^4\partial_i v\partial_j v\big]\partial_{ij} v\\
&\qquad-a^3+4a^2\langle y, \nabla v\rangle
-2a\langle y, \nabla v\rangle^2
-3a|y|^2|\nabla v|^2
+2|y|^2\langle y, \nabla v\rangle|\nabla v|^2,
\end{split}\end{equation}
for the constant $a$ in \eqref{eq-MSE-expansion-identity-dim-2}, and, for $n\ge 3$, 
\begin{equation}\label{eq-MSE-dim3-above-fully}\begin{split}
F&=-|y|^{2n-2}\big[(n-2)^2v^2+2(n-2)v\langle y, \nabla v\rangle+|y|^2|\nabla v|^2\big]\Delta v\\
&\qquad +|y|^{2n-4}\big[(n-2)^2y_iy_jv^2+(n-2)|y|^2v(y_i\partial_j v+y_j\partial_i v)+|y|^4\partial_i v\partial_j v\big]\partial_{ij} v\\
&\qquad+|y|^{2n-4}\big[(n-2)^3(n-1)v^3+4(n-2)^2(n-1)v^2\langle y, \nabla v\rangle\\
&\qquad +(n-2)(3n-4)v\langle y, \nabla v\rangle^2+(n-2)(2n-1)v|y|^2|\nabla v|^2\\
&\qquad +2(n-1)|y|^2\langle y, \nabla v\rangle|\nabla v|^2\big].
\end{split}\end{equation}
We note that $F$ is a polynomial of $y$, $v$, $y\nabla v$, and $y^2\nabla^2v$ and, in fact, linear in $\nabla^2v$. 


We first verify that all terms in \eqref{eq-MSE-dim2-fully} and \eqref{eq-MSE-dim3-above-fully} 
form operators preserving the degree of $|y|^2$. 

\begin{lemma}
\label{lemma-operators-in-MSE}
The following operators are of type $(i,j)$ that preserve the degree of $|y|^2$, 
for some $i \geq 1$ and $j \geq 0$:
\begin{itemize}
  \item[(1)] $u\mapsto |y|^2 \Delta u$, $u\mapsto\langle y,\nabla u\rangle$, and $u\mapsto y_iy_j \partial_{ij}u$ are of type $(1,0)$;
  \item[(2)] $(u,v) \mapsto |y|^2 \langle\nabla u, \nabla v\rangle$ 
  and $(u,v)\mapsto |y|^2 y_i \partial_j u \partial_{ij}v$ are of type $(2,0)$;
  \item[(3)] $(u,v,w)\mapsto |y|^4 \partial_i u \partial_j v \partial_{ij} w$ is of type $(3,0)$.
\end{itemize} 
\end{lemma}

The proof is based on \eqref{eq-partial_i(|x|^2l u)} and \eqref{eq-partial_ij(|x|^2l u)}, and is omitted.

We now prove the main result in this section. 

\begin{theorem}\label{thrm-harmonic-fcn-correspondence-MSE} 
Let $h$ be a harmonic function in $B_1\subseteq \R^n$. Then, 
there exists a unique analytic function $w_*$ in $B_R$ for some $R \in (0,1)$ such that 
$v = h + |y|^2w_*$ is a solution to \eqref{eq-nonlinear-Poission-MSE-fully} in $B_R$. 
\end{theorem}

\begin{proof}  
Set $w = |y|^2 w_*$ and consider 
\begin{align}\label{eq-v=H+w-MSE}v = h + w.\end{align}
By substituting \eqref{eq-v=H+w-MSE} in \eqref{eq-nonlinear-Poission-MSE-fully}, we get 
\begin{align}
\label{eq-nonlinear-Poission-MSE-fully-with-harmonic-fcn}\Delta w = F_h(y,w,y\nabla w, y^2 \nabla^2 w),
\end{align}
where $F_h$ is given by
$$F_h(y,w,y\nabla w, y^2\nabla^2 w) = F(y, h+w, y\nabla h + y\nabla w, y^2 \nabla^2 h + y^2 \nabla^2 w).$$
We will prove that  
there exists a unique analytic function $w_*$ in $B_R$ for some $R\in (0,1)$ such that $w=|y|^2 w_*$
is a solution to \eqref{eq-nonlinear-Poission-MSE-fully-with-harmonic-fcn} in $B_R$. 

Note that \eqref{eq-nonlinear-Poission-MSE-fully-with-harmonic-fcn} has the same structure 
as \eqref{eq-fully-nonlinear-equation} with $u$ replaced by $w$ and $x$ replaced by $y$. 
We will verify that $F_h$ satisfies the assumption of Theorem \ref{thrm-existence-fully-nonlinear}.
First, it is obvious that $F_h$ is analytic in $B_1\times \R\times \R^{n^2}\times \R^{n^4}$.
Next, we will prove that, for any $C^2$-function $w$ near $0$ with $w(0) = 0$,
\begin{align}
\label{eq-Taylor-expansion-F_h-MSE}
F_h(y,w,y\nabla w, y^2\nabla^2 w) = f_h(y) +  \sum_{j=0}^\infty T_{1,j}(w) 
+ \sum_{j=0}^\infty T_{2,j}(w,w) + \sum_{j=0}^\infty T_{3,j}(w,w,w),
\end{align}
for some analytic function $f_h$ depending on $h$ and operators $T_{ij}$ of type $(i,j)$ 
that preserve the degree of $|y|^2$, for $i = 1,2,3$ and $j \ge 0$.

For an illustration, we verify this assertion for $|y|^{2n}  \partial_i v \partial_j v \partial_{ij } v$ in $F$ for $n\ge 2$. 
With \eqref{eq-v=H+w-MSE}, we have 
\begin{align*}
|y|^{2n}  \partial_i v \partial_j v \partial_{ij} v&= |y|^{2n}\partial_i h \partial_j h \partial_{ij} h
+2|y|^{2n}  \partial_i h \partial_j w \partial_{ij} h+ |y|^{2n}  \partial_i w \partial_j w \partial_{ij} h\\
&\qquad  +|y|^{2n}\partial_i h \partial_j h \partial_{ij} w
+2|y|^{2n}  \partial_i h \partial_j w \partial_{ij} w+ |y|^{2n}  \partial_i w \partial_j w \partial_{ij} w.
\end{align*}
We first examine the three terms involving $\partial_{ij} w$. 
Let $\{H_m\}_{m\ge 0}$ be the sequence of the homogeneous Taylor polynomials of $h$ at $0$, 
with $H_m \in \mathscr{H}_m$. Then, 
\begin{align*}
|y|^{2n} \partial_i h \partial_j h \partial_{ij} w  &= |y|^{2n-4} \sum_{m=0}^\infty \sum_{p+q = m} |y|^{4} \partial_i H_p \partial_j H_q \partial_{ij} w ,\\
|y|^{2n}  \partial_i h \partial_j w \partial_{ij} w
&= |y|^{2n-4}\sum_{m=0}^\infty |y|^{4} \partial_i H_m \partial_j w \partial_{ij} w,\\
|y|^{2n}  \partial_i w \partial_j w \partial_{ij} w
&= |y|^{2n-4}\big( |y|^{4} \partial_i w \partial_j w \partial_{ij} w\big).
\end{align*}
The corresponding operators are given by 
\begin{align*}
T_{1,m}: &\,w \mapsto |y|^{2n-4} \sum_{p+q = m-2n+4} |y|^4  \partial_i H_p \partial_j H_q \partial_{ij} w,\\
T_{2,m}: &\,(w_1,w_2) \mapsto  |y|^{2n-4}  \big(|y|^4 \partial_i H_{m-2n+4} \partial_j w_1 \partial_{ij} w_2 \big),\\
T_{3,2n-4}: &\,(w_1,w_2,w_3)\mapsto |y|^{2n-4} \big( |y|^4 \partial_i w_1 \partial_j w_2 \partial_{ij} w_3\big).
\end{align*}
By Lemma \ref{lemma-constructing-new-operators} and Lemma \ref{lemma-operators-in-MSE}, 
$T_{1,m}$ is of type $(1,m)$, $T_{2,m}$ is of type $(2,m)$ for $m \ge 2n-4$, and $T_{3,2n-4}$ is of type $(3,2n-4)$, 
and they all preserve the degree of $|y|^2$. 
A similar argument applies to $|y|^{2n}  \partial_i h \partial_j w \partial_{ij} h$ 
and $|y|^{2n}  \partial_i w \partial_j w \partial_{ij} h$. 
Note that $|y|^{2n}\partial_i h \partial_j h \partial_{ij} h$ is an analytic function independent of $w$, 
and hence contributes to $f_h$ in \eqref{eq-Taylor-expansion-F_h-MSE}. 
We can discuss other terms in $F$ and $F_h$ similarly. 
Therefore, by Theorem \ref{thrm-existence-fully-nonlinear}, 
\eqref{eq-nonlinear-Poission-MSE-fully-with-harmonic-fcn} 
admits a unique solution $w$
in $B_R$ for some $R\in (0,1)$, 
with $w=|y|^2w_*$ 
for some analytic function $w_*$ in $B_R$.
\end{proof} 

Theorem \ref{thrm-harmonic-fcn-correspondence-MSE-intro} follows from 
Theorem \ref{thrm-harmonic-fcn-correspondence-MSE}. 

\section{Monge-Amp\`{e}re equation in exterior domains}
\label{sec-MA}

In this section, we study the Monge-Ampère equation in $\R^n \setminus \overline{B}_1$. 
Let $u\in C^2(\mathbb{R}^n\setminus \overline{B}_1)$ be a locally convex solution of 
\begin{equation}\label{eq-MA}\det(\nabla^2 u)=1 \quad\text{in } \R^n\setminus \overline{B}_1.\end{equation}
Han and Wang \cite{Han-Wang2024} proved the following result.

\begin{theorem}
\label{thrm-MA-decomposition}
Let $n \ge 2$ and $u \in C^2(\R^n \setminus \overline{B}_1)$ be a locally convex solution of \eqref{eq-MA}.

$\mathrm{(i)}$ For $n = 2$, there exist some $d\in\mathbb{R}$, $b \in \mathbb{R}^2$, 
a positive definite matrix $A$ with $\det A = 1$, 
and a function $v\in C^\infty(B_r)$ for some $r\in (0,1)$, such that, for any $x\in  \R^2\setminus \overline{B}_R$ for some $R>0$, 
\begin{equation}\label{eq-MA-expansion-identity-dim-2}u(x)=\frac{1}{2}x^TAx+ \langle b, x\rangle+d\log\sqrt {x^TAx}+
v\Big(\frac{\sqrt{A}x}{|\sqrt{A}x|^2}\Big).\end{equation}

$\mathrm{(ii)}$ For $n \geq 3$, there exist some $c\in\mathbb{R}$, $b \in \mathbb{R}^n$, 
a positive definite matrix $A$ with $\det A = 1$, 
and a function $v$ in $B_r$ for some $r\in (0,1)$, with $v\in C^\infty(B_r)$ for $n$ even and $v\in C^{n-1,1}(B_r)$ for $n$ odd, 
such that, for any $x\in  \R^n\setminus \overline{B}_R$ for some $R>0$, 
\begin{equation}\label{eq-MA-expansion-identity-dim-large}u(x)=\frac{1}{2}x^TAx+\langle b, x\rangle+c+
|\sqrt{A}x|^{2-n}v\Big(\frac{\sqrt{A}x}{|\sqrt{A}x|^2}\Big).\end{equation}
Moreover, for $n$ odd, there exist functions $\widetilde{v}, \widehat{v}\in C^\infty(B_r)$ such that, for any $y\in B_r$, 
\begin{equation}\label{eq-MA-decomposition-dim-odd}
v(y)=\widetilde{v}(y)+|y|^{n}\widehat{v}(y).\end{equation}
\end{theorem}

The regularity of $v\in C^{n-1,1}(B_r)$ in the odd dimension is optimal. 
The decomposition \eqref{eq-MA-decomposition-dim-odd} clearly illustrates that $|y|^{n}$ is an obstruction, 
in fact the only obstruction, to the higher regularity of $v$. 

First, by an appropriate affine transform, we can always assume $A=I$, the identity matrix. 
This is done for all dimensions. We now record the equation satisfied by $v$. 
Refer to \cite{Han-Wang2024} for details. 

Consider the Kelvin transform 
$$y=\frac{x}{|x|^2}.$$
As in Section \ref{sec-MSE}, we switch the notations for variables from $x\in \mathbb R^n\setminus\overline{B}_1$
to $y\in B_1$.  
The function $v$ in \eqref{eq-MA-expansion-identity-dim-2} and \eqref{eq-MA-expansion-identity-dim-large} satisfies
\begin{align}
\label{eq-MA-equation-for-v}
\Delta v = -\sum_{k=2}^n |y|^{kn-n-2} \sigma_k\big( M(v) \big),
\end{align}
where $\sigma_k$ is the $k$-th elementary symmetric function and $M(v)$ is an $n\times n$ symmetric matrix given by 
the following expressions: for $n = 2$, \begin{equation}
\label{eq-M(v)-n=2}
\begin{aligned}
M(v)_{ij}&=\big[-2d + 8 \langle y,\nabla v\rangle + 4 y_ky_\ell \partial_{k\ell} v \big]y_i y_j|y|^{-2}\\
&\qquad +\big[ d- 2 \langle y,\nabla v\rangle  \big] \delta_{ij} - 2(y_i \partial_j v + y_j \partial_i v)\\
&\qquad - 2 y_i y_k \partial_{kj} v - 2 y_j y_\ell \partial_{i\ell} v + |y|^2 \partial_{ij} v,
\end{aligned}
\end{equation}
for the constant $d$ in \eqref{eq-MA-expansion-identity-dim-2}, and, for $n \ge 3$, 
\begin{equation}
\label{eq-M(v)-n>=3}
\begin{aligned}
M(v)_{ij}&=\big[ n(n-2)v + 4n \langle y,\nabla v\rangle + 4 y_ky_\ell \partial_{k\ell} v \big]y_i y_j|y|^{-2}\\
&\qquad-\big[ (n-2)v+ 2 \langle y,\nabla v\rangle  \big] \delta_{ij} - n(y_i \partial_j v + y_j \partial_i v)\\
&\qquad - 2 y_i y_k \partial_{kj} v - 2 y_j y_\ell \partial_{i\ell} v + |y|^2 \partial_{ij} v.
\end{aligned}
\end{equation}
In terms of the eigenvalues $\lambda_1,\cdots,\lambda_n$ of a symmetric matrix $M$, 
$\sigma_k(M)$ for $k = 1,\cdots, n$ is given by
$$\sigma_k(M) = \sum_{1 \leq i_1<\cdots < i_k \leq n} \lambda_{i_1}\cdots \lambda_{i_k}.$$
It is well-known that $\sigma_k(M)$ is a polynomial in $\mathrm{Tr}(M)$, $\cdots$, $\mathrm{Tr}(M^k)$. 
For example, for $k=2$, we have 
\begin{align}\label{eq-trace-square}
\sigma_2(M) = \frac{1}{2} [\mathrm{Tr}(M)]^2 - \frac{1}{2}\mathrm{Tr}(M^2).
\end{align}


Note that $ kn-n-2$, as in \eqref{eq-MA-equation-for-v},
is even if $n$ is even or both $n$ and $k$ are odd. In these cases, 
$|y|^{kn-n-2}$ is a polynomial. However, $kn-n-2$ is odd if $n$ is odd and $k$ is even, 
in which case the right-hand side of \eqref{eq-MA-equation-for-v} contains non-smooth expressions. 
As a consequence, different methods are needed for even dimensions and odd dimensions. 
In fact, the study for the even dimensions is similar as that for the minimal surface equation, 
but more complicated for the odd dimensions.

Let $v$ be a smooth function in some open set $\Omega \subseteq \R^n$. 
A straightforward computation based on \eqref{eq-M(v)-n=2} and \eqref{eq-M(v)-n>=3} yields
\begin{align}\label{eq-trace-Laplace}\mathrm{Tr} \big(M(v)\big) = |y|^2 \Delta v.\end{align}
By Lemma \ref{lemma-operators-in-MSE}(1), $v \mapsto \mathrm{Tr}\big(M(v) \big)$ 
is an operator of type $(1,0)$ 
that preserves the degree of $|y|^2$. 
We now prove that higher order traces are operators that preserve the degree of $|y|^2$.
We first consider $n \ge 3$. 

\begin{lemma}
\label{lemma-Monge-Ampere-operators-n>=3}
Let $n \ge 3$ and $M$ be the matrix operator given by \eqref{eq-M(v)-n>=3}. For $2 \leq i \leq n$, define
\begin{align}
\label{eq-T-operator-trace-M(P_i)-n>=3}
T(v_1,\cdots,v_i) = \mathrm{Tr}\big( M(v_1)\cdots M(v_i) \big).
\end{align}
Then, $T$ is an operator of type $(i,0)$ that preserves the degree of $|y|^2$.
\end{lemma}

\begin{proof}
First, we consider a matrix of the form
\begin{align}
\label{eq-M(P)-general-structure}
A = a |y|^{-2} y y^T +b I_n + y c^T + dy^T + |y|^2 e,
\end{align}
where $a,b \in C^\infty (\R^n)$, $c,d$ are $n$-row vectors whose components are smooth functions, 
and $e$ is an $n \times n$ matrix whose components are smooth functions. Then, 
\begin{align}
\label{eq-Trace(A)}
\mathrm{Tr}(A) = a + nb + \langle y,c\rangle + \langle y,d\rangle + |y|^2 \mathrm{Tr}(e)
\end{align}
is a smooth function in $\R^n$.
We say that $A$ in \eqref{eq-M(P)-general-structure} has degree $i\ge 2$ if all functions are homogeneous polynomials, 
$a,b$ have degree $i$, 
the components of $c,d$ have degree $i-1$, and the components of $e$ have degree $i-2$.
Similarly, we say that $A$ has degree $1$ if $a,b$ have degree $1$,
the components of $c,d$ have degree $0$ and $e = 0$; 
and $A$ has degree $0$ if $a,b$ have degree $0$, $c,d = 0$, and $e = 0$.
It is easy to verify that if $A$ has degree $i$, 
then $\mathrm{Tr}(A)$ given by \eqref{eq-Trace(A)} is a homogeneous polynomial of degree $i$.
Now, for a polynomial $P$ in $\R^n$, we can write $M(P)$ as in \eqref{eq-M(P)-general-structure}, with 
\begin{align*}
a &= \big[ n(n-2)P + 4n \langle y,\nabla P\rangle + 4 y_ky_\ell \partial_{k\ell} P \big],\\
b &= -\big[ (n-2)P+ 2 \langle y,\nabla P\rangle  \big],\\
c &= d = -n\nabla P -2 (\nabla^2 P) y,\\
e &= \nabla^2 P.
\end{align*}
It is clear that if $P$ is a homogeneous polynomial of degree $i$, then $M(P)$ given by \eqref{eq-M(v)-n>=3} has degree $i$. 
We now carry out the proof in two steps.

{\it Step 1. Assume that $A_1$ and $A_2$ have the structures given by \eqref{eq-M(P)-general-structure}, 
then $A_1A_2$ also has the structure given by \eqref{eq-M(P)-general-structure}. 
Moreover, if $A_1$ has degree $i$ and $A_2$ has degree $j$, then $A_1A_2$ has degree $i+j$.}
To verify this, we assume 
\begin{align*}
&A_1 = a_1 |y|^{-2} y y^T +b_1 I_n + y c_1^T + d_1y^T + |y|^2 e_1,\\
&A_2 = a_2 |y|^{-2} y y^T +b_2 I_n + y c_2^T + d_2y^T + |y|^2 e_2.
\end{align*}
A straightforward computation yields 
\begin{align*}
A_1 A_2 &= \big[a_1a_2 + a_1b_2 + a_1 \langle y, d_2 \rangle + b_1 a_2 
+ a_2 \langle y,c_1 \rangle + |y|^2 \langle c_1,d_2 \rangle\big] |y|^{-2} y y^T\\
&\qquad+(b_1b_2)I_n
+ y \big[ a_1c_2^T + a_1y^Te_2 + b_1c_2^T + b_2c_1^T + \langle y,c_1\rangle c_2^T + |y|^2 c_1^T e_2\big]\\
&\qquad+(b_1d_2 + a_2d_1 + d_1b_2 + \langle y,d_2\rangle d_1 +a_2e_1y) y^T\\
&\qquad+|y|^2 \big(b_1e_2 + d_1c_2^T + d_1y^Te_2 + b_2e_1 + e_1yc_2^T + e_1d_2y^T + |y|^2e_1e_2\big).
\end{align*}
Hence, the assertion follows. 

{\it Step 2. Let $m$ be a positive integer and $v$ be a smooth function in $\R^n$. Then,
\begin{align}
\label{eq-M(p)-preserves-deg}
M(|y|^m v ) = |y|^m A_m(v),
\end{align}
for some matrix operator $A_m$, and $A_m(v)$ can be written as in \eqref{eq-M(P)-general-structure}.} 
By a straightforward computation, we have 
$$\nabla (|y|^m v) = m |y|^{m-2} v y + |y|^m \nabla v,$$
and
\begin{align*}
\nabla^2(|y|^m v) &= m(m-2) |y|^{m-4} v yy^T + m |y|^{m-2} v I_n\\
&\qquad+ m |y|^{m-2} (y \nabla v^T + \nabla v y^T) + |y|^m \nabla^2 v.
\end{align*}
Then,
\begin{align*}
y \nabla (|y|^m  v )^T &=|y|^m\big[ mv|y|^{-2} y y^T + y \nabla v^T \big],\\
y y^T \nabla^2(|y|^m v ) &=|y|^m \big\{\big[ m(m-1)v +m \langle y,\nabla v\rangle \big] |y|^{-2} yy^T\\
&\qquad\qquad  +my\nabla v^T+yy^T\nabla^2v\big\},
\end{align*}
and
\begin{align*}
|y|^2\nabla^2(|y|^m v ) &= |y|^m \big[ m(m-2) v |y|^{-2} yy^T + m v I_n\\
&\qquad+ m (y \nabla v^T + \nabla v y^T) + |y|^2 \nabla^2 v \big].
\end{align*}
We also have
\begin{align*}
\nabla (|y|^m  v ) y^T &= |y|^m \big[ mv|y|^{-2} y y^T +\nabla v y^T \big],\\
 \nabla^2(|y|^m v ) y y^T &= |y|^m \big\{\big[ {m(m-1)v} +m \langle y,\nabla v\rangle \big] |y|^{-2} yy^T  \\
&\qquad\qquad+ m   \nabla vy^T +  \nabla^2 v y y^T \big\}.
\end{align*}
It is straightforward to check that \eqref{eq-M(p)-preserves-deg} holds and $A_m(v)$ 
can be put in the same form as in \eqref{eq-M(P)-general-structure}.

Combining Step 1 and Step 2, we conclude that the operator $T$ given by \eqref{eq-T-operator-trace-M(P_i)-n>=3} 
is of type $(i,0)$ that preserves the degree of $|y|^2$.
\end{proof}

\begin{remark}
\label{rmk-operator-A-in-M(|x|^m P) = |x|^m A(P)}
By a straightforward computation as in Step 2 in the proof of Lemma \ref{lemma-Monge-Ampere-operators-n>=3}, 
the operator $A_m$ in \eqref{eq-M(p)-preserves-deg} also satisfies
$$A_m(|x|^m v) = |x|^m B_m(v),$$
for some matrix operator $B_m$, and $B_m(v)$ can be written as in \eqref{eq-M(P)-general-structure}.
\end{remark}

Next, we consider $n = 2$. With $M(v)$ in \eqref{eq-M(v)-n=2}, we write
\begin{align}\label{eq-M(v)-n=2-decomposition}
M(v) = -2d |y|^{-2} y y^T + d I_2 + \tilde M(v),
\end{align}
where
\begin{equation}
\label{eq-tilde-M(v)-n=2}
\begin{aligned}
\tilde M(v)&=[ 8 \langle y,\nabla v\rangle + 4 y_ky_\ell \partial_{k\ell} v ]|y|^{-2}yy^T 
-2 \langle y,\nabla v\rangle I_2\\ 
&\qquad -2y( \nabla v^T + y^T\nabla^2 v) 
 - 2(\nabla v + \nabla^2 v y ) y^T + |y|^2 \nabla^2 v.
\end{aligned}
\end{equation}
Hence,
\begin{align*}
\mathrm{Tr}\big( M(v)^2 \big)& = \mathrm{Tr}\big( (-2d|y|^{-2}yy^T + dI_2)^2 \big)\\
&\qquad+ 2\mathrm{Tr}\big( (-2d|y|^{-2}yy^T + dI_2)\cdot \tilde M(v)\big)
+ \mathrm{Tr}\big( \tilde M(v)^2\big).
\end{align*}
As in Lemma \ref{lemma-Monge-Ampere-operators-n>=3}, we have the following result.

\begin{lemma}
\label{lemma-Monge-Ampere-operators-n=2}
Let $n = 2$, $d$ be a constant, and $\tilde M$ be the operator given by \eqref{eq-tilde-M(v)-n=2}. Define
\begin{align*}
T(v) = \mathrm{Tr}\big( (-2d|y|^{-2} yy^T + dI_2)\cdot \tilde M(v) \big),
\end{align*}
and
\begin{align*}
T'(v,w) = \mathrm{Tr}\big( \tilde M(v) \tilde M(w) \big).
\end{align*}
Then, $T$ is an operator of type $(1,0)$ and $T'$ is an operator of type $(2,0)$ 
that preserve the degree of $|y|^2$.
\end{lemma}

The proof is similar as that of Lemma \ref{lemma-Monge-Ampere-operators-n>=3} and is omitted.

We now prove our first main result in this section, concerning the solvability of \eqref{eq-MA-equation-for-v} 
in even dimensions. 

\begin{theorem}\label{thrm-harmonic-fcn-correspondence-MA-even} 
Let $n\ge 2$ be even and $h$ be a harmonic function in $B_1$. Then, 
there exists a unique analytic function $w_*$ in $B_R$ for some $R\in (0,1)$ such that 
$v = h + |y|^2w_*$ is a solution to \eqref{eq-MA-equation-for-v} in $B_R$. 
\end{theorem}

\begin{proof} 
Since $n$ is even, $|y|^{kn-n-2}$, as in \eqref{eq-MA-equation-for-v}, is a homogeneous polynomial of degree $kn-n-2$ 
for $k = 2,\cdots, n$. Set 
\begin{align}\label{eq-MA-expression-F-for-v}
F(y,v,y\nabla v, y^2\nabla^2 v) = - \sum_{k=2}^n |y|^{kn-n-2} \sigma_k\big( M(v) \big).\end{align}
We write \eqref{eq-MA-equation-for-v} as 
\begin{align}\label{eq-MA-equation-F-for-v}
\Delta v=F(y,v,y\nabla v, y^2\nabla^2 v).\end{align}
It is obvious that $F$ is a polynomial in its components in 
$\R^n \times \R \times \R^{n^2} \times \R^{n^4}$. 
Recall that, for a smooth function $v$, we write 
$$y \nabla v:\,\text{the collection of all }y_i \partial_j v\text{ for }1 \leq i ,j\leq n,$$
and
$$y^2 \nabla^2 v:\,\text{the collection of all }y_i y_j \partial_{k\ell} v\text{ for }1 \leq i ,j,k,\ell\leq n.$$

Set $w = |y|^2 w_*$ and consider 
\begin{align}\label{eq-v=H+w-MA}v = h + w.\end{align}
By substituting \eqref{eq-v=H+w-MA} in \eqref{eq-MA-equation-F-for-v}, we get 
\begin{align}
\label{eq-nonlinear-Poission-MA-fully-with-harmonic-fcn}\Delta w = F_h(y,w,y\nabla w, y^2 \nabla^2 w),
\end{align}
where $F_h$ is given by
$$F_h(y,w,y\nabla w, y^2\nabla^2 w) = F(y, h+w, y\nabla h + y\nabla w, y^2 \nabla^2 h + y^2 \nabla^2 w).$$
We will prove that  
there exists a unique analytic function $w_*$ in $B_R$ for some $R\in (0,1)$ such that $w=|y|^2 w_*$
is a solution to \eqref{eq-nonlinear-Poission-MA-fully-with-harmonic-fcn} in $B_R$.

Note that \eqref{eq-nonlinear-Poission-MA-fully-with-harmonic-fcn} has the same structure 
as \eqref{eq-fully-nonlinear-equation} with $u$ replaced by $w$ and $x$ replaced by $y$. 
We will verify that $F_h$ satisfies the assumption of Theorem \ref{thrm-existence-fully-nonlinear}.
First, it is obvious that $F_h$ is analytic in $B_1\times \R\times \R^{n^2}\times \R^{n^4}$.
Next, we will prove that, for any $C^2$-function $w$ near $0$ with $w(0) = 0$,
\begin{align}
\label{eq-Taylor-expansion-F_h-MA}
F_h(y,w,y\nabla w, y^2\nabla^2 w) = f_h(y) +  \sum_{i=1}^n\sum_{j=0}^\infty T_{ij}(w,\cdots, w),
\end{align}
for some analytic function $f_h$ depending on $h$ and operators $T_{ij}$ of type $(i,j)$ 
that preserve the degree of $|y|^2$, for $i = 1,\cdots,n$ and $j \ge 0$.

For an illustration, we verify this assertion for $|y|^{n-2}\sigma_2\big( M(v) \big)$
in $F$, corresponding to $k=2$ in \eqref{eq-MA-expression-F-for-v}. 
We first consider $n\ge 4$. Note that $M(v)$ is linear in $v$ by \eqref{eq-M(v)-n>=3}. 
By \eqref{eq-trace-square} and \eqref{eq-v=H+w-MA}, we have 
\begin{align*}
\sigma_2(M(v)) = \frac{1}{2} [\mathrm{Tr}(M(h))+\mathrm{Tr}(M(w))]^2 
- \frac{1}{2}\mathrm{Tr}[(M(h)+M(w))^2].
\end{align*}
Note that $\mathrm{Tr}\big( M(h) \big) = 0$ by \eqref{eq-trace-Laplace}, since $h$ is harmonic.
A straightforward calculation yields
\begin{align*}
2|y|^{n-2}\sigma_2\big(M(v) \big)&=
|y|^{n-2} \big(|y|^2\Delta w\big)^2 -|y|^{n-2}  \mathrm{Tr}\big(M(h)^2 \big)\\
&\qquad- 2|y|^{n-2}  \mathrm{Tr}\big(M(h) M(w) \big)-|y|^{n-2}  \mathrm{Tr}\big(M(w)^2 \big).
\end{align*}
It is obvious that $|y|^{n-2}  \mathrm{Tr}\big(M(h)^2 \big)$ is an analytic function of $y$, independent of $w$. 
Hence, it contributes to $f_h$ in \eqref{eq-Taylor-expansion-F_h-MA}.
We now prove that all the other terms yield operators that preserve the degree of $|y|^2$. 
Let $\{H_j\}_{j\ge 0}$ be the sequence of the homogeneous Taylor polynomials of $h$ at $0$, with $H_j \in \mathscr{H}_j$. 
Then, the corresponding operators are given by
\begin{align*}
T_{1,j}: \,& w\mapsto |y|^{n-2}  \mathrm{Tr}\big(M(H_{j-n+2}) M(w) \big),\\
T_{2,n-2}: \,&(w_1,w_2) \mapsto |y|^{n-2} \big(|y|^2 \Delta w_1 \big)\big(|y|^2 \Delta w_2 \big) 
-|y|^{n-2} \mathrm{Tr}\big(M(w_1) M(w_2) \big).
\end{align*}
For $n \geq 4$, by Lemma \ref{lemma-constructing-new-operators}, Lemma \ref{lemma-operators-in-MSE}, 
and Lemma \ref{lemma-Monge-Ampere-operators-n>=3}, 
we conclude that $T_{1,j}$ is of type $(1,j)$ for $j \ge n-2$ and $T_{2,n-2}$ is of type $(2,n-2)$, 
and they all preserve the degree of $|y|^2$. 
Next, we consider $n =2$. By \eqref{eq-M(v)-n=2-decomposition}, \eqref{eq-tilde-M(v)-n=2}, and \eqref{eq-v=H+w-MA}, 
we have 
$$M(v)=-2d |y|^{-2} y y^T + d I_2+\tilde M(h)+\tilde M(w)=M(h)+\tilde M(w).$$
Note that $\tilde M(v)$ as given by \eqref{eq-tilde-M(v)-n=2} is linear in $v$. 
Now, with Lemma \ref{lemma-Monge-Ampere-operators-n>=3} replaced 
by Lemma \ref{lemma-Monge-Ampere-operators-n=2}, we get a similar result.  
The same argument applies to other terms in $F$ and $F_h$. 
Therefore, by Theorem \ref{thrm-existence-fully-nonlinear} , 
\eqref{eq-nonlinear-Poission-MA-fully-with-harmonic-fcn} 
admits a unique solution $w$
in $B_R$ for some $R\in (0,1)$, 
with $w=|y|^2w_*$ 
for some analytic function $w_*$ in $B_R$.
\end{proof} 

We next prove our second main result in this section, concerning the solvability of \eqref{eq-MA-equation-for-v} 
in odd dimensions. 

\begin{theorem}\label{thrm-harmonic-fcn-correspondence-MA-odd} 
Let $n\ge 3$ be odd and $h$ be a harmonic function in $B_1$. Then, 
there exist unique analytic functions $\widetilde{w}_*, \widehat{w}_*$ in $B_R$ for some $R\in (0,1)$ such that 
$v = h + |y|^2\widetilde{w}_*+|y|^n \widehat{w}_*$ is a solution to \eqref{eq-MA-equation-for-v} in $B_R$. 
\end{theorem}

\begin{proof} With $F$ as in \eqref{eq-MA-expression-F-for-v}, the equation \eqref{eq-MA-equation-for-v} is 
\eqref{eq-MA-equation-F-for-v}. 
As observed earlier, $kn-n-2$ is odd if $n$ is odd and $k$ is even, 
in which case the right-hand side of \eqref{eq-MA-equation-F-for-v} contains non-smooth factors $|y|^{kn-n-2}$. 
First, we decompose \eqref{eq-MA-equation-F-for-v} into two equations with smooth dependence. 
To this end, we write $v$ as
\begin{align}
\label{eq-MA-decomposition-dim-odd-purely-unharmonic-second-part}
v = \widetilde{v} + |y|^{n-2} \widehat{v}.
\end{align}
We point out that this is different from \eqref{eq-MA-decomposition-dim-odd}. 
Here, a factor of $|y|^2$ in \eqref{eq-MA-decomposition-dim-odd} is absorbed 
by $\widehat{v}$.
Substituting \eqref{eq-MA-decomposition-dim-odd-purely-unharmonic-second-part} 
in \eqref{eq-MA-equation-F-for-v}, we obtain
\begin{align}\label{eq-MA-dim-odd-system}
\begin{split}
\Delta \widetilde{v} &= \widetilde{F}\big(y,\widetilde{v}, \widehat{v}, y \nabla \widetilde{v}, y \nabla \widehat{v}, 
y^2 \nabla^2 \widetilde{v}, y^2 \nabla^2 \widehat{v} \big),\\
\Delta (|y|^{n-2} \widehat{v} )&= |y|^{n-2} \widehat{F}\big(y,\widetilde{v}, \widehat{v}, y \nabla \widetilde{v}, y \nabla \widehat{v}, 
y^2 \nabla^2 \widetilde{v}, y^2 \nabla^2 \widehat{v} \big),
\end{split}
\end{align}
for some polynomials $\widetilde{F}$ and $\widehat{F}$ in $y,\, \widetilde{v},\, \widehat{v},\, y \nabla \widetilde{v},\, 
y \nabla \widehat{v},\, y^2 \nabla^2 \widetilde{v},$ and $y^2 \nabla^2 \widehat{v}$. 
This has the same structure as \eqref{eq-fully-nonlinear-system} with $m = n-2$, an odd number.

For an illustration, we demonstrate how to derive \eqref{eq-MA-dim-odd-system} 
by computing $|y|^{n-2} \sigma_2 \big( M(v) \big)$, 
corresponding to $k=2$ in \eqref{eq-MA-expression-F-for-v}. 
With $m=n-2$, \eqref{eq-M(p)-preserves-deg} implies 
$$M\big( |y|^{n-2} \widehat{v} \big) = |y|^{n-2} A_{n-2} (\widehat{v}),$$
for some matrix operator $A_{n-2}$. For simplicity, we write $A = A_{n-2}$. 
By \eqref{eq-trace-square}, we have
\begin{align*}
\sigma_2 \big( M(v) \big) 
 &= \frac{1}{2} \big[\mathrm{Tr}\big( M(\widetilde{v}) + |y|^{n-2}A(\widehat{v}) \big)\big]^2 
 - \frac{1}{2}\mathrm{Tr}\big[ \big(M(\widetilde{v}) + |y|^{n-2} A(\widehat{v})\big)^2 \big].
\end{align*}
A straightforward calculation yields 
\begin{align*}
2|y|^{n-2} \sigma_2 \big( M(v) \big)  
&=|y|^{n-2} \big(|y|^2 \Delta \widetilde{v}\big)^2 
+|y|^{3n-6}  \big[\mathrm{Tr}\big(A(\widehat{v}) \big)\big]^2\\
&\qquad -|y|^{n-2} \mathrm{Tr} \big(M(\widetilde{v})^2 \big)
 -|y|^{3n-6} \mathrm{Tr}\big( A(\widehat{v})^2 \big)\\
&\qquad+2|y|^{2n-4} \big( |y|^2 \Delta \widetilde{v}\big)\mathrm{Tr}\big(A(\widehat{v}) \big)
-2|y|^{2n-4}\mathrm{Tr}\big( M(\widetilde{v}) A(\widehat{v}) \big),
\end{align*}
where we used $\mathrm{Tr}\big( M(\widetilde{v}) \big) = |y|^2 \Delta \widetilde{v}$ by \eqref{eq-trace-Laplace}. 
The first four terms contain an odd power of $|y|^{n-2}$ and hence contribute to $|y|^{n-2} \widehat{F}$ 
in \eqref{eq-MA-dim-odd-system}, 
while the last two terms contain an even power of $|y|^{n-2}$ and hence contribute to $\widetilde{F}$. 
Note that the operator $A$ has a similar structure as the operator $M$. 
By the computation in the proof of Lemma \ref{lemma-Monge-Ampere-operators-n>=3}, 
$\mathrm{Tr}\big(A(\widehat{v}) \big)$ and $\mathrm{Tr}\big(A(\widehat{v})^2 \big)$ 
are polynomials in $y$, $\widehat{v}$, $y\nabla \widehat{v}$, 
and $y^2 \nabla^2 \widehat{v}$, and $\mathrm{Tr}\big(M(\widetilde v)A(\widehat{v}) \big)$ 
is a polynomial in $y$, $\widetilde v,\widehat{v}$, $y\nabla \widetilde v,y\nabla \widehat{v}$, 
and $y^2 \nabla^2 \widetilde v,y^2 \nabla^2 \widehat{v}$. The same argument applies to $|y|^{kn-n-2} \sigma_k\big( M(v) \big)$
in \eqref{eq-MA-expression-F-for-v}, 
for $3 \leq k \leq n$. Therefore, we conclude 
\eqref{eq-MA-dim-odd-system}.

Now, we rewrite \eqref{eq-MA-decomposition-dim-odd-purely-unharmonic-second-part} as
\begin{align*}
v = h + \widetilde{w} + |y|^{n-2} \widehat{w},
\end{align*}
or 
\begin{align}
\label{eq-MA-decomposition-dim-odd-with-Q}
\widetilde v = h + \widetilde{w}, \quad \widehat{v} = \widehat{w},
\end{align}
where $\widetilde{w}$ and $\widehat{w}$ are given by 
$\widetilde{w}=|y|^2\widetilde{w}_*$ and $\widehat{w}=|y|^2\widehat{w}_*$ 
for some analytic functions $\widetilde{w}_*$ and $\widehat{w}_*$. 
By substituting \eqref{eq-MA-decomposition-dim-odd-with-Q} in \eqref{eq-MA-dim-odd-system}, 
we obtain
\begin{align}
\label{eq-MA-dim-odd-system-with-Q}\begin{split}
\Delta \widetilde{w} &= \widetilde{F}_h\big(y,\widetilde{w}, \widehat{w}, y \nabla \widetilde{w}, 
y \nabla \widehat{w}, y^2 \nabla^2 \widetilde{w}, y^2 \nabla^2 \widehat{w} \big),\\
\Delta (|y|^{n-2} \widehat{w} ) &= |y|^{n-2} \widehat{F}_h\big(y,\widetilde{w}, \widehat{w}, y \nabla \widetilde{w}, 
y \nabla \widehat{w}, y^2 \nabla^2 \widetilde{w}, y^2 \nabla^2 \widehat{w} \big),
\end{split}
\end{align}
for some analytic functions $\widetilde{F}_h$ and $\widehat{F}_h$.
We point out that $\widetilde{F}_h$ and $\widehat{F}_h$ are analytic in $y$ and are polynomials in 
$\widetilde{w}, \widehat{w}$, $y \nabla \widetilde{w}, y \nabla \widehat{w}$, and $y^2 \nabla^2 \widetilde{w}, y^2 \nabla^2 \widehat{w}$.
By Remark \ref{rmk-operator-A-in-M(|x|^m P) = |x|^m A(P)}, we can repeat the argument 
as in the proof of Theorem \ref{thrm-harmonic-fcn-correspondence-MA-even} 
with Theorem \ref{thrm-existence-fully-nonlinear} replaced by Theorem \ref{thrm-existence-fully-nonlinear-system}. 
Therefore, \eqref{eq-MA-dim-odd-system-with-Q} admits a unique pair of solutions $(\widetilde{w}, \widehat{w})$
in $B_R$ for some $R\in (0,1)$, 
with $\widetilde{w}=|y|^2\widetilde{w}_*$ and $\widehat{w}=|y|^2\widehat{w}_*$ 
for some analytic functions $\widetilde{w}_*$ and $\widehat{w}_*$ in $B_R$.
\end{proof} 

Theorem \ref{thrm-harmonic-fcn-correspondence-MA-intro}
follows from Theorem \ref{thrm-harmonic-fcn-correspondence-MA-even}
and Theorem \ref{thrm-harmonic-fcn-correspondence-MA-odd}. 

\section{Optimal regularity for the Monge-Amp\`{e}re equation near infinity in odd dimensions}
\label{sec-optimal-decomposition-MA-odd}

In this section, we assume that $n\ge 3$ is odd
and prove that the decomposition in Theorem \ref{thrm-MA-decomposition}
can be improved from a smooth decomposition to an analytic one.
Note that the smooth decomposition in
\eqref{eq-MA-decomposition-dim-odd} is not unique.
In fact, if $\varphi$ is a nonzero smooth function near $0$ such that
$D^k \varphi(0) = 0$ for any $k \ge 0$, then
$v$ can be written in two different ways as
$$v=\widetilde{v}+|y|^{n}\widehat{v} = 
(\widetilde{v} + \varphi)+|y|^{n}(\widehat{v} - |y|^{-n}\varphi).$$
The proof of our analytic decomposition is based directly 
on \eqref{eq-MA-equation-for-v}, 
without decomposing it into a system as in \eqref{eq-MA-dim-odd-system}.

We first introduce some notation.
Fix $\gamma\in (0,1)$. 
For a nonnegative integer $k$, set
\[
\mathscr P_{\le k}=\bigoplus_{j=0}^k\mathscr P_j,
\]
and set $\mathscr P_{\le k}=\{0\}$ if $k<0$. 
Define
\begin{align}
\label{eq-def-V_k-MA-odd}
\mathscr V_k
=
\mathscr P_{\le k}
+
|y|^n\mathscr P_{\le k-n}.
\end{align}
Since $n$ is odd, the representation in \eqref{eq-def-V_k-MA-odd}
is unique. 

We next introduce a scaling that will be used below. 
If $v$ is a solution of \eqref{eq-MA-equation-for-v} near the origin and
$\epsilon>0$, set
\[
v_\epsilon(y)=v(\epsilon y).
\]
It follows directly from \eqref{eq-M(v)-n>=3} that
\[
M(v_\epsilon)(y)=M(v)(\epsilon y).
\]
Consequently, $v_\epsilon$ satisfies
\begin{align}
\label{eq-MA-scaled-odd}
\Delta v_\epsilon
=
-\sum_{i=2}^n
\epsilon^{(i-1)n}
|y|^{(i-1)n-2}
\sigma_i\big(M(v_\epsilon)\big).
\end{align}

For $2\le i\le n$, let $S_i$ be the unique symmetric $i$-linear functional
such that
\[
S_i(A,\ldots,A)=\sigma_i(A),
\]
for any symmetric matrix $A$.
We also denote by $\Delta^{-1}$ the zero-Dirichlet inverse of the
Laplacian on $B_1$. Namely, for $f\in C^{\gamma}(\overline B_1)$,
$v=\Delta^{-1}f$ is the unique solution of
\begin{align*}
\begin{cases}
\Delta v=f &\text{in }B_1,\\
v=0 &\text{on }\partial B_1.
\end{cases}
\end{align*}

The following lemma lists two properties of the nonlinear
terms that will be used in the proof.

\begin{lemma}
\label{lemma-MA-block-operator-odd}
For $2\le i\le n$, define
\begin{align}
\label{eq-def-T_i-MA-odd}
\mathcal T_i(v_1,\ldots,v_i)
=
\Delta^{-1}\big(
|y|^{(i-1)n-2}
S_i(M(v_1),\ldots,M(v_i))
\big).
\end{align}
Then, for any $v_1,\cdots, v_i \in C^{2,\gamma}(\overline B_1)$,
\begin{align}
\label{eq-estimate-T_i-MA-odd}
|\mathcal T_i(v_1,\ldots,v_i)|_{C^{2,\gamma}(\overline B_1)}
\le
C_i\prod_{j=1}^i
|v_j|_{C^{2,\gamma}(\overline B_1)},
\end{align}
for some positive constant $C_i$ depending only on 
$n$, $i$, and $\gamma$.
Moreover, for any 
$v_1 \in \mathscr V_{k_1},\cdots, v_i \in \mathscr V_{k_i}$,
\begin{align}
\label{eq-T_i-preserves-V_k-MA-odd}
\mathcal T_i(v_1,\ldots,v_i)
\in
\mathscr V_{k_1+\cdots+k_i+(i-1)n}.
\end{align}
\end{lemma}

\begin{proof}
The estimate \eqref{eq-estimate-T_i-MA-odd} follows directly from
Lemma \ref{lemma-Monge-Ampere-operators-n>=3}
and the Schauder estimate for the zero-Dirichlet problem for the Laplacian.

We now prove \eqref{eq-T_i-preserves-V_k-MA-odd}. 
Let $v_1 \in \mathscr V_{k_1}, \cdots, v_i \in \mathscr V_{k_i}$.
By Lemma \ref{lemma-Monge-Ampere-operators-n>=3} and its proof,
it is straightforward to check that
\[
S_i(M(v_1),\ldots,M(v_i))\in \mathscr V_{k_1 + \cdots + k_i}.
\]
On the other hand, for any fixed $m \ge 0$ and $H_j \in \mathscr H_j$,
we have
\[
\Delta^{-1}(|y|^mH_j)
=
\frac{(|y|^{m+2}-1)H_j}
{(m+2)(2j+m+n)}.
\]
By decomposing $P_k \in \mathscr P_k$ as in \eqref{eq-harmonic-decomposition},
we then obtain \eqref{eq-T_i-preserves-V_k-MA-odd}.
\end{proof}

We next prove a uniform estimate. For $B \in \mathscr V_k$,
we can write
\begin{align}
\label{eq-expression-B-components}
B = \sum_{j=0}^k P_j + |y|^n \sum_{j=0}^{k-n} Q_j,
\end{align}
for some $P_j,Q_j \in \mathscr P_j$.
We will control the $L^\infty$-norms of $P_j$ and $Q_j$ by 
$|B|_{L^\infty (\overline B_1)}$.

\begin{lemma}
\label{lemma-MA-separate-block-odd}
Fix $A>1$ and let $B\in\mathscr V_k$ be given by 
\eqref{eq-expression-B-components} with $P_j,Q_j\in\mathscr P_j$. Then, 
\begin{align}
\label{eq-separate-block-estimate-MA-odd}
\sum_{j=0}^k
\rho^j|P_j|_{L^\infty(\mathbb S^{n-1})}
+
\sum_{j=0}^{k-n}
\rho^{n+j}|Q_j|_{L^\infty(\mathbb S^{n-1})}
\le
4A^k|B|_{L^\infty(\overline B_1)},
\end{align}
for some constant $\rho\in(0,1)$ depending only on $A$.
\end{lemma}

\begin{proof}
We first recall the following estimate for complex polynomials in $\C$.
Set
$$\widetilde E_A
=
\left\{
z\in\mathbb C:
|z-1|+|z+1|\le A+A^{-1}
\right\}$$
and
$$
 E_A
= \Big\{
z\in\mathbb C:
|z|+|z-1|
\le
\frac{A+A^{-1}}{2}
\Big\}.
$$
Let $q$ be a complex polynomial of degree at most $k$. Then
\begin{align}
\label{eq-Bernstein--Walsh}
|q|_{L^\infty(\widetilde E_A)}
\le
A^k|q|_{L^\infty([-1,1])}.
\end{align}
This is a corollary of the Bernstein--Walsh inequality.
Refer to \cite{Bernstein1912,Walsh1926} and to 
\cite[Corollary 12.1.6]{RahmanSchmeisser} for a modern treatment.
By a simple reparameterization, 
for any complex polynomial $p$ of degree at most $k$, 
applying \eqref{eq-Bernstein--Walsh} to $q(z)=p\big((z+1)/2\big)$ yields
\begin{align}
\label{eq-Bernstein--Walsh-rescaled}
|p|_{L^\infty( E_A)}
\le
A^k|p|_{L^\infty([0,1])}.
\end{align}
Since $A>1$, $0$ is in the interior of $E_A$. In particular, 
there exists $\rho\in (0,1)$ such that $\{z;\,|z| \le 2\rho\} \subseteq E_A$.
By the Cauchy estimate in complex analysis and
\eqref{eq-Bernstein--Walsh-rescaled}, 
if $p(z) = \sum_{j=0}^k c_j z^j$, we have
\begin{align}
\label{eq-Bernstein--Walsh-polynomial-coeff}
|c_j|
\le
(2\rho)^{-j}A^k|p|_{L^\infty([0,1])}.
\end{align}

Now, fix 
$$B = P+|y|^nQ = \sum_{j=0}^k P_j + |y|^n \sum_{j=0}^{k-n} Q_j \in \mathscr V_k,$$
where $P_j, Q_j \in \mathscr P_j$.
For each $\theta\in\mathbb S^{n-1}$, define
\[
p_\theta(t)
=P(t\theta) + t^n Q(t\theta),
\quad t\in\mathbb R.
\]
Then $p_\theta$ is a polynomial of degree at most $k$, 
$p_\theta(t) = B(t\theta)$ if $t \ge 0$, and
\[
p_\theta(t)=\sum_{j=0}^k \big[P_j(\theta) + Q_{j-n}(\theta) \big] t^j \equiv\sum_{j=0}^k c_j(\theta) t^j,
\]
where we set $Q_\ell =0$ if $\ell < 0$.
We now view $p_\theta$ as a complex polynomial.
By \eqref{eq-Bernstein--Walsh-polynomial-coeff},
\[
|c_j(\theta)|
\le
(2\rho)^{-j}A^k
|p_\theta|_{L^\infty([0,1])}
\le (2\rho)^{-j}A^k|B|_{L^\infty(\overline B_1)}.
\]
By taking the supremum over $\theta \in \mathbb S^{n-1}$ and 
summing over $j$, we get
\begin{align}
\label{eq-cj-summability-MA-odd}
\sum_{j=0}^k
\rho^j
|c_j|_{L^\infty(\mathbb S^{n-1})}
\le
2A^k|B|_{L^\infty(\overline B_1)}.
\end{align}
Finally, since $n$ is odd, $P_j$ and $Q_{j-n}$ have opposite parity;
namely, if $P_j$ is even, then $Q_{j-n}$ is odd, and vice versa. Thus,
\[
P_j(\theta)
=
\frac12
\left[
c_j(\theta)+(-1)^jc_j(-\theta)
\right]
\quad\text{and}\quad Q_{j-n}(\theta)
=
\frac12
\left[
c_j(\theta)-(-1)^jc_j(-\theta)
\right].
\]
It follows that
\[
|P_j|_{L^\infty(\mathbb S^{n-1})}
+
|Q_{j-n}|_{L^\infty(\mathbb S^{n-1})}
\le
2|c_j|_{L^\infty(\mathbb S^{n-1})}.
\]
Together with \eqref{eq-cj-summability-MA-odd}, this proves
\eqref{eq-separate-block-estimate-MA-odd}.
\end{proof}

In the next result, we review some standard estimates for smooth and analytic
functions on $\mathbb S^{n-1}$.
\begin{lemma}
\label{lemma-smooth-analytic-spherical-fcn}
Let $g$ be a smooth function on $\mathbb S^{n-1}$. 
Write $g=\sum_{k=0}^\infty g_k$ in the $L^2$ sense,
where $g_k$ is a spherical harmonic of degree
$k$, and let $H_k\in\mathscr H_k$ be the unique homogeneous 
harmonic polynomial such that $H_k|_{\mathbb S^{n-1}} =g_k$.
Then, we have
\begin{align}
\label{eq-C-2,gamma-smooth-bound}
\sum_{k=0}^\infty |H_k|_{C^{2,\gamma}(\overline B_1)} \le C \|g\|_{H^{n+2}(\mathbb S^{n-1})},
\end{align}
for some positive constant $C$ depending only on $n$ and $\gamma$.
Moreover, if $g$ is analytic, i.e., 
$g$ is the restriction of a real-analytic function 
defined in a neighborhood of $\mathbb S^{n-1}$, then 
\begin{align}
\label{eq-C2-gamma-analytic-bound}
\sum_{k=0}^\infty A^k |H_k|_{C^{2,\gamma}(\overline B_1)} <\infty,
\end{align}
for some positive constant $A>1$ depending on $g$.
\end{lemma}

\begin{proof}
For simplicity, denote by $C$ a positive constant 
(different from line to line) depending only on $n$ and $\gamma$.
Let $\lambda_k=k(k+n-2)$ be the eigenvalue of $-\Delta_{\mathbb S^{n-1}}$ 
corresponding to spherical harmonics of degree $k$.
By the global Schauder estimate for harmonic functions and the Sobolev embedding,
we have, for any $k \ge 0$,
\begin{align}
\label{eq-C^2,gamma-harmonic-global-estimate}
|H_k|_{C^{2,\gamma}(\overline B_1)}
\le
C|g_k|_{C^{2,\gamma}(\mathbb S^{n-1})}
\le C\|g_k\|_{H^{n+1}(\mathbb S^{n-1})}.
\end{align}
Recall that spherical harmonics of different degrees are orthogonal
with respect to $\langle\cdot,\cdot\rangle_{H^m(\mathbb S^{n-1})}$ 
for any $m \ge 1$.
Hence, the series representation $g = \sum_{k=0}^\infty g_k$ holds
in $H^m(\mathbb S^{n-1})$ for any $m \ge 1$.
By the Cauchy inequality,
\begin{align*}
\sum_{k=0}^\infty
\|g_k\|_{H^{n+1}(\mathbb S^{n-1})}
&\le
C\Big(
\sum_{k=0}^\infty(1+k)^{-2}
\Big)^{1/2}
\Big(
\sum_{k=0}^\infty
(1+k)^2
\|g_k\|_{H^{n+1}(\mathbb S^{n-1})}^2
\Big)^{1/2}\\
&\le
C\|g\|_{H^{n+2}(\mathbb S^{n-1})}.
\end{align*}
This proves \eqref{eq-C-2,gamma-smooth-bound}.

Suppose now that $g$ is analytic. Then,
there exist constants $C_0,L>0$ such that, for any $m\ge0$,
\begin{align}
\label{eq-analytic-laplacian-sphere}
\|(-\Delta_{\mathbb S^{n-1}})^m g\|_{L^2(\mathbb S^{n-1})}
\le
C_0L^{2m}(2m)!.
\end{align}
Refer to \cite{Seeley1969}.
Since the spherical harmonics of different degrees are orthogonal
with respect to the $L^2$-inner product on $\mathbb S^{n-1}$, we have
\[
\lambda_k^m\|g_k\|_{L^2(\mathbb S^{n-1})}
\le
\|(-\Delta_{\mathbb S^{n-1}})^m g\|_{L^2(\mathbb S^{n-1})}.
\]
Thus, by \eqref{eq-analytic-laplacian-sphere}, for any $k \ge 1$,
\begin{align}
\label{eq-L^2-upperbound-with-m}
\|g_k\|_{L^2(\mathbb S^{n-1})}
\le
C_0\lambda_k^{-m}L^{2m}(2m)! \le C_0
\Big(\frac{2Lm}{k}\Big)^{2m}.
\end{align}
Choose $\delta>0$ sufficiently small such that $2L\delta<1$. 
Taking $m = [\delta k]$ in \eqref{eq-L^2-upperbound-with-m},
and noting that $m \ge (\delta k)/2$ for $k$ sufficiently large,
we get
\[
\|g_k\|_{L^2(\mathbb S^{n-1})}
\le
C e^{-ck},
\]
for some constants $c$ and $C$.
By enlarging $C$ if necessary, the same estimate holds for all $k\ge0$.
Choosing any $A>1$ sufficiently close to $1$ such that
$Ae^{-c}<1$, \eqref{eq-C^2,gamma-harmonic-global-estimate} then yields
\[
\sum_{k=0}^\infty
A^k|H_k|_{C^{2,\gamma}(\overline B_1)} \le C \sum_{k=0}^\infty A^k (1+k)^{n+1} \|g_k\|_{L^2(\mathbb S^{n-1})} \le C \sum_{k=0}^\infty
(1+k)^{n+1}(Ae^{-c})^k
<\infty.
\]
This proves \eqref{eq-C2-gamma-analytic-bound}.
\end{proof}

We are now ready to prove the optimal decomposition.

\begin{theorem}
\label{thrm-MA-optimal-decomposition-odd}
Let $n\ge3$ be odd and $v$ be a solution of \eqref{eq-MA-equation-for-v} in $B_1$.
Then, there exist a unique harmonic function $h$ and unique analytic
functions $\widetilde w_*$ and $\widehat w_*$ in $B_R$, for some
$R\in (0,1)$, such that
\begin{align}
\label{eq-optimal-decomposition-MA-odd}
v
=
h+|y|^2\widetilde w_*+|y|^n\widehat w_*
\quad\text{in }B_R.
\end{align}
\end{theorem}

\begin{proof}
The uniqueness part is straightforward; we will prove the existence.
The proof consists of three steps.

{\it Step 1. We solve a Dirichlet boundary value problem for the scaled
equation \eqref{eq-MA-scaled-odd}.}
Consider
\begin{align}
\label{eq-scaled-Dirichlet-MA-odd}
\begin{cases}
\displaystyle
\Delta V
=
-\sum_{i=2}^n
\epsilon^{(i-1)n}
|y|^{(i-1)n-2}
\sigma_i\big(M(V)\big)
&\text{in }B_1,\\
V=g&\text{on }\mathbb S^{n-1},
\end{cases}
\end{align}
where $g\in C^\infty(\mathbb S^{n-1})$.
Write $g=\sum_{k=0}^\infty g_k$ as in 
Lemma \ref{lemma-smooth-analytic-spherical-fcn},
where $g_k$ is a spherical harmonic of degree
$k$, and let $H_k\in\mathscr H_k$ be the unique homogeneous 
harmonic polynomial such that $H_k|_{\mathbb S^{n-1}} =g_k$.
Assume 
\begin{align}
\label{eq-boundary-smallness}
\sum_{k=0}^\infty
A^k|H_k|_{C^{2,\gamma}(\overline B_1)}\le\frac K2,
\end{align}
for some constants $K>0$ and $A >1$.
Let $C_2,\cdots,C_n$ be the constants in \eqref{eq-estimate-T_i-MA-odd}.
We claim that, for $\epsilon>0$ sufficiently small with
\begin{align}
\label{eq-smallness-iteration-MA-odd}
\sum_{i=2}^n
C_i(A\epsilon)^{(i-1)n}K^i
\le\frac K2,
\quad
\sum_{i=2}^n
iC_i(A\epsilon)^{(i-1)n}K^{i-1}
\le\frac12,
\end{align}
there exists a solution $V\in C^{2,\gamma}(\overline B_1)$
of \eqref{eq-scaled-Dirichlet-MA-odd}.

To this end, define the weighted space
$
\mathcal X_A
=
\left\{
\mathcal B=(B_k)_{k\ge0};
B_k\in\mathscr V_k,\ 
\|\mathcal B\|_A<\infty
\right\},
$
where
\begin{align}
\label{eq-X_A-norm-MA-odd}
\|\mathcal B\|_A
=
\sum_{k=0}^\infty
A^k|B_k|_{C^{2,\gamma}(\overline B_1)}.
\end{align}
It is clear that $\mathcal X_A$ is a Banach space.
For $\mathcal B = (B_k)_{k\ge 0}\in \mathcal X_A$, 
define $\Phi(\mathcal B)$ by
\begin{align}
\label{eq-iteration-MA-odd}
\Phi(\mathcal B)_k
=
H_k
-
\sum_{i=2}^n
\epsilon^{(i-1)n}
\sum_{\substack{
k_1+\cdots+k_i+(i-1)n=k\\
k_1,\ldots,k_i\ge0
}}
\mathcal T_i(B_{k_1},\ldots,B_{k_i}).
\end{align}
By Lemma \ref{lemma-MA-block-operator-odd},
$\Phi(\mathcal B)_k\in\mathscr V_k$ for $k \ge 0$.
Moreover, by \eqref{eq-estimate-T_i-MA-odd} and a simple rearrangement, 
\begin{align}
\label{eq-Phi-estimate-MA-odd}
\|\Phi(\mathcal B)\|_A
\le
\sum_{k=0}^\infty
A^k|H_k|_{C^{2,\gamma}(\overline B_1)}
+
\sum_{i=2}^n
C_i(A\epsilon)^{(i-1)n}\|\mathcal B\|_A^i.
\end{align}
Indeed, for each $2 \le i \le n$, we have
\begin{align*}
\sum_{k\ge0}A^k
\sum_{\substack{
k_1+\cdots+k_i+(i-1)n=k\\
k_1,\ldots,k_i\ge0
}}
\prod_{j=1}^i
|B_{k_j}|_{C^{2,\gamma}}
=
A^{(i-1)n}
\Big(
\sum_{k\ge0}
A^k|B_k|_{C^{2,\gamma}}
\Big)^i.
\end{align*}
Similarly, if $\|\mathcal B\|_A,\|\widetilde{\mathcal{B}}\|_A\le K$, then the multilinearity
of $S_i$ gives
\begin{align}
\label{eq-Phi-contraction-MA-odd}
\|\Phi(\mathcal B)-\Phi(\widetilde{\mathcal{B}})\|_A
\le
\Big[
\sum_{i=2}^n
iC_i(A\epsilon)^{(i-1)n}K^{i-1}
\Big]
\|\mathcal B-\widetilde{\mathcal{B}}\|_A.
\end{align}
Combining \eqref{eq-boundary-smallness}, 
\eqref{eq-smallness-iteration-MA-odd},
\eqref{eq-Phi-estimate-MA-odd}, and \eqref{eq-Phi-contraction-MA-odd},
we conclude that $\Phi$ is a contraction on the closed ball
$
\{\mathcal B\in \mathcal X_A:\|\mathcal B\|_A\le K\}.
$
By the contraction mapping theorem, there exists a unique
$\mathcal B=(B_k)_{k\ge0}$ in this ball satisfying
$\mathcal B=\Phi(\mathcal B).$
Set
\begin{align}
\label{eq-Dirichlet-solution-V}
V=\sum_{k=0}^\infty B_k.
\end{align}
Since $A\ge1$, the series in \eqref{eq-Dirichlet-solution-V}
converges absolutely in $C^{2,\gamma}(\overline B_1)$.
Moreover, by \eqref{eq-iteration-MA-odd} and the definition of 
$\mathcal T_i$ as in \eqref{eq-def-T_i-MA-odd},
for any $k \ge 0$,
we have 
$$B_k|_{\mathbb S^{n-1}} = \Phi(\mathcal B)_k|_{\mathbb S^{n-1}} = H_k|_{\mathbb S^{n-1}} = g_k.$$
By a simple summation, $V$ satisfies the boundary condition in
\eqref{eq-scaled-Dirichlet-MA-odd}.
On the other hand, by \eqref{eq-estimate-T_i-MA-odd}
and a simple summation of 
\eqref{eq-iteration-MA-odd} over $k$, we obtain
\[
V
=
\sum_{k=0}^\infty H_k
-
\sum_{i=2}^n
\epsilon^{(i-1)n}
\Delta^{-1}\Big(
|y|^{(i-1)n-2}\sigma_i(M(V))
\Big).
\]
Therefore, $V$ solves the equation in \eqref{eq-scaled-Dirichlet-MA-odd}.
We point out that $A = 1$ is sufficient for the contraction mapping argument.
The strict inequality $A>1$ is needed later, when we apply
Lemma \ref{lemma-smooth-analytic-spherical-fcn}.

{\it Step 2. We prove that the solution $V$ of \eqref{eq-scaled-Dirichlet-MA-odd}
constructed in Step 1 admits
an analytic decomposition near $0$ given by
\begin{align}
\label{eq-constructed-P-plus-radial-Q-MA-odd}
V=\widetilde V+|y|^n\widehat V \quad\text{in }B_R,
\end{align}
for some constant $R \in (0,1)$ and analytic functions $\widetilde V$
and $\widehat V$ in $B_R$.}
For each $k$, since $B_k \in \mathscr V_k$, we can write
\[
B_k=P_k+|y|^nQ_k,
\]
where
$P_k\in\mathscr P_{\le k}$ and $Q_k\in\mathscr P_{\le k-n}$,
with $P_k=\sum_{j=0}^k P_{k,j}$ and $Q_k=\sum_{j=0}^{k-n} Q_{k,j}$, 
$P_{k,j}, Q_{k,j} \in \mathscr P_j$.
By Lemma \ref{lemma-MA-separate-block-odd} and a simple summation, we have 
\begin{align*}
\sum_{k=0}^\infty\sum_{j=0}^k
\rho^j
|P_{k,j}|_{L^\infty(\mathbb S^{n-1})}
+
\sum_{k=0}^\infty\sum_{j=0}^{k-n}
\rho^{n+j}
|Q_{k,j}|_{L^\infty(\mathbb S^{n-1})}
\le
C\sum_{k=0}^\infty
A^k|B_k|_{L^\infty(\overline B_1)}
<\infty,
\end{align*}
for some constant $\rho > 0$ depending only on $A$.
By switching the order of summation, we obtain
\[
\sum_{j=0}^\infty\rho^j \sum_{k=j}^\infty
|P_{k,j}|_{L^\infty(\mathbb S^{n-1})}
+
\sum_{j=0}^\infty\rho^{n+j}
\sum_{k=n+j}^\infty|Q_{k,j}|_{L^\infty(\mathbb S^{n-1})}
\le C\sum_{k=0}^\infty
A^k|B_k|_{L^\infty(\overline B_1)}.
\]
Set
$$\widetilde P_j = \sum_{k=j}^\infty
P_{k,j}\quad\text{and}\quad 
\widehat P_j = \sum_{k=n+j}^\infty Q_{k,j}.$$
Then, $\widetilde P_j,\widehat P_j \in \mathscr P_j$ and 
$\{\widetilde P_j\}_{j\ge 0}, \{\widehat P_j\}_{j\ge 0}$ satisfy 
analyticity-type estimates as in \eqref{eq-analyticity-estimate-general}.
By Corollary \ref{cor-constructing-analytic-functions}, the functions
$\widetilde V=\sum_{j=0}^\infty \widetilde P_j$
and $\widehat V=\sum_{j=0}^\infty \widehat P_j$
are analytic in $B_R$ for some $R \in (0,1)$. 
Thus, the desired decomposition
\eqref{eq-constructed-P-plus-radial-Q-MA-odd} holds.

{\it Step 3. We prove the decomposition 
\eqref{eq-optimal-decomposition-MA-odd} for a general solution
$v$ of \eqref{eq-MA-equation-for-v} in $B_1$.}

By Theorem \ref{thrm-MA-decomposition}, we have
\[
v=\widetilde v+|y|^n\widehat v
\]
near the origin, for some smooth functions $\widetilde v$ and $\widehat v$.
For $0<\epsilon\ll1$, set
\[
v_\epsilon(y)=v(\epsilon y)
\quad\text{and}\quad
g_\epsilon=v_\epsilon|_{\mathbb S^{n-1}}.
\]
By \eqref{eq-MA-scaled-odd}, $v_\epsilon$ satisfies
\eqref{eq-scaled-Dirichlet-MA-odd} with the boundary value $g_\epsilon$
on $\mathbb S^{n-1}$.
The smooth decomposition of $v$ gives
\[
g_\epsilon(\theta)
=
\widetilde v(\epsilon\theta)
+
\epsilon^n\widehat v(\epsilon\theta).
\]
Hence, for all sufficiently small $\epsilon>0$,
\begin{align}
\label{eq-uniform-v_epsilon-g_epsilon-bound}
|v_\epsilon|_{C^{2,\gamma}(\overline B_1)}
\le M_0 \quad\text{and}\quad
\|g_\epsilon\|_{H^{n+2}(\mathbb S^{n-1})}
\le M_1,
\end{align}
for some positive constants $M_0,M_1$ independent of $\epsilon$.

Write
\[
g_\epsilon=\sum_{k=0}^\infty g_{\epsilon,k},
\]
where $g_{\epsilon,k}$ is a spherical harmonic of degree $k$, and let
$H_{\epsilon,k}\in\mathscr H_k$ be the homogeneous harmonic extension
of $g_{\epsilon,k}$. By \eqref{eq-uniform-v_epsilon-g_epsilon-bound}
and \eqref{eq-C-2,gamma-smooth-bound}, we have
\begin{align}
\label{eq-uniform-eta1-MA-odd}
\sum_{k=0}^\infty
|H_{\epsilon,k}|_{C^{2,\gamma}(\overline B_1)}
\le M,
\end{align}
for some positive constant $M$ independent of $\epsilon$.

Set
\[
K=\max\{4M,M_0,1\}.
\]
Then, choose $\epsilon>0$ sufficiently small such that
\eqref{eq-smallness-iteration-MA-odd} holds for $A = 2$.
We fix such an $\epsilon$ in the rest of the proof.
Since $\epsilon>0$ is fixed, $g_\epsilon$ is analytic on
$\mathbb S^{n-1}$. Indeed, $v$ is real analytic away from the origin.
By Lemma \ref{lemma-smooth-analytic-spherical-fcn}, there exists
$A_\epsilon>1$ such that \eqref{eq-C2-gamma-analytic-bound} holds with $A$ and
$H_k$ replaced by $A_\epsilon$ and $H_{\epsilon,k}$, respectively.
Then, the function 
$A \mapsto \sum_{k=0}^\infty 
A^k|H_{\epsilon,k}|_{C^{2,\gamma}(\overline B_1)}$
is continuous on $[0,A_\epsilon]$.
Hence, by \eqref{eq-uniform-eta1-MA-odd}, we may choose
$1<A_*<\min\{A_\epsilon,2\}$
such that
\begin{align}
\label{eq-eta-A-epsilon-MA-odd}
\sum_{k=0}^\infty
A_*^k|H_{\epsilon,k}|_{C^{2,\gamma}(\overline B_1)}
\le 2M
\le \frac K2.
\end{align}
Since $A_*<2$, \eqref{eq-smallness-iteration-MA-odd} holds for $A_*$
by our choice of $\epsilon$.
Therefore, all the assumptions in Step 1 are satisfied. We obtain a
solution $V$ of \eqref{eq-scaled-Dirichlet-MA-odd} with
boundary value $g_\epsilon$ such that
$|V|_{C^{2,\gamma}(\overline B_1)}
\le K$.
Moreover, by Step 2, \eqref{eq-constructed-P-plus-radial-Q-MA-odd} holds
for some analytic functions $\widetilde V$ and
$\widehat V$ near the origin.

We next prove that $v_\epsilon=V$.
Set
\[
Z=v_\epsilon-V.
\]
Since $v_\epsilon$ and $V$ have the same boundary value,
$Z=0$ on $\mathbb S^{n-1}$. Subtracting their equations, we obtain
\[
Z
=
-\sum_{i=2}^n
\epsilon^{(i-1)n}
\Delta^{-1}
\left\{
|y|^{(i-1)n-2}
\left[
\sigma_i(M(v_\epsilon))
-
\sigma_i(M(V))
\right]
\right\}.
\]
By \eqref{eq-smallness-iteration-MA-odd} and
a similar estimate to that in \eqref{eq-Phi-contraction-MA-odd}, 
\begin{align*}
|Z|_{C^{2,\gamma}(\overline B_1)}
\le
\left[
\sum_{i=2}^n
iC_i\epsilon^{(i-1)n}K^{i-1}
\right]
|Z|_{C^{2,\gamma}(\overline B_1)}
\le
\frac12
|Z|_{C^{2,\gamma}(\overline B_1)}.
\end{align*}
We point out that
\[
|v_\epsilon|_{C^{2,\gamma}(\overline B_1)}
\le M_0\le K,
\qquad
|V|_{C^{2,\gamma}(\overline B_1)}
\le K.
\]
Therefore, $v_\epsilon=V$. In particular,
$v_\epsilon=\widetilde V+|y|^n\widehat V$ near the origin. 
Rescaling back, the desired decomposition 
\eqref{eq-optimal-decomposition-MA-odd} then follows from
Proposition \ref{prop-decomposition-ana}. 
\end{proof}

\section{Continuous Dependence}
\label{sec-conti-depend}

In this section, we study the dependence of solutions to the minimal surface equation 
and the Monge-Amp\`{e}re equation in exterior domains 
on harmonic functions. 
We should point out that the results in this section do not rely on the existence in the previous two sections.

First, we recall the following estimate for composition of analytic functions. 
It will be used in our derivations of uniform estimates. 
Refer to Proposition 1.4.2 in \cite{KrantzParks2002} for a proof of the $1$-dimensional case. 

\begin{lemma}
\label{lemma-estimate-composition-of-analytic-functions}
Let $U \subseteq \R^{m_1}$ be a neighborhood of $0$ in $\R^{m_1}$ and $V$ be a neighborhood of $0$ in $\R^{m_2}$, 
and $g : U\rightarrow \R^{m_2}$ and $f : V\rightarrow \R$ be analytic functions with $g(0) = 0$. 
Assume, 
for any $\beta \in \Z_+^{m_2}\setminus\{0\}$ and any $\gamma \in \Z_+^{m_1}\setminus\{0\}$, 
$$\frac{1}{\beta!}|\partial^{\beta} f(0)| \leq A_0 A^{|\beta|}\quad\text{and}\quad \frac{1}{\gamma!}|\partial^\gamma g(0)| \leq B_0 B^{|\gamma|},$$
for some positive constants $A_0, A, B_0$, and $B$. Then, for any $\alpha \in \Z_+^{m_1}\setminus\{0\}$,
$$\frac{1}{\alpha!}|\partial^{\alpha} (f\circ g)(0)| \leq C_0C^{|\alpha|},$$
for some positive constants $C_0$ and $C$ depending only on $A_0, A, B_0,B$, and $m_1,m_2$.
\end{lemma}

We start with the minimal surface equation \eqref{eq-MSE} in $\R^n \setminus \overline{B}_1$. 
As in Section \ref{sec-MSE}, we study the equation
\begin{align}\label{eq-nonlinear-Poission-MSE-fully-1} 
\Delta v=F\quad\text{in }B_1,
\end{align}
where $F$ is given by \eqref{eq-MSE-dim2-fully} for $n = 2$ and \eqref{eq-MSE-dim3-above-fully} for $n \ge 3$. 
We point out that $F$ in \eqref{eq-nonlinear-Poission-MSE-fully-1} is an analytic function with respect to 
its components and any smooth solution $v$ to \eqref{eq-nonlinear-Poission-MSE-fully-1} is analytic. 
Hence, if $v$ is a solution to \eqref{eq-nonlinear-Poission-MSE-fully-1}, by Proposition \ref{prop-decomposition-ana}, 
we can decompose $v$ as
$$v = h + |y|^2 w_*\quad\text{in }B_R,$$
for some harmonic function $h$ and analytic function $w_*$ in $B_R$, for some $R\in (0,1)$. 
Conversely, by Theorem \ref{thrm-harmonic-fcn-correspondence-MSE}, for any harmonic function $h$ in $B_1$, 
there exists a unique analytic function $w_*$ in $B_R$ for some $R \in (0,1)$ such that $v = h + |y|^2 w_*$ 
is a solution to \eqref{eq-nonlinear-Poission-MSE-fully-1} in $B_R$. 
In the following, we study the dependence of $w = |y|^2 w_*$ on $h$.
For simplicity, we consider $n \ge 3$.
The discussion for $n = 2$ is similar, and the constants in the estimates
depend additionally on $a$ in \eqref{eq-MSE-expansion-identity-dim-2}.

Recall that $\langle \cdot,\cdot \rangle_{*}$ is the inner product defined in \eqref{eq-H^m-inner-product-with-constant}
and $\|\cdot\|_*$ is the corresponding norm. 

\begin{theorem}
\label{thrm-conti-dependence-MSE}
Let $h_1$ and $h_2$ be harmonic functions in $B_1$, 
$w_{*1}$ and $w_{*2}$ be analytic functions in $B_{R_0}$ for some $R_0 \in (0,1)$ 
such that $v_1 = h_1+ |y|^2 w_{*1}$ and $v_2 = h_2 + |y|^2 w_{*2}$ 
are solutions to \eqref{eq-nonlinear-Poission-MSE-fully-1} in $B_{R_0}$, 
and $\{P_{1,k}\}_{k\ge 2}$ and $\{P_{2,k}\}_{k\ge 2}$ be the sequences 
of the homogeneous Taylor polynomials of $|y|^2 w_{*1}$ and $|y|^2 w_{*2}$ at $0$, respectively, 
with $P_{1,k}, P_{2,k} \in \mathscr{H}_k^\perp$. 
Assume that
\begin{align*}
\|h_1\|_{L^2(B_1)}, \|h_2\|_{L^2(B_1)} \leq M_0,
\end{align*}
for some constant $M_0>0$, and, for any $k \ge 2$,
\begin{align}
\label{eq-upper-bound-purely-unharmonic-part-conti-depend-MSE}
\|P_{1,k}\|_*,\|P_{2,k}\|_* \leq C_*^k,
\end{align}
for some constant $C_*>0$. Then, for any $k\ge 2$,
\begin{align}
\label{eq-estimate-difference-conti-depend-MSE}
\|P_{1,k} - P_{2,k}\|_* \leq  C^k\|h_1 - h_2\|_{L^2(B_1)},
\end{align}
for some constant $C>0$ depending only on $n$, $M_0$, and $C_*$. 
\end{theorem}

\begin{proof}
For simplicity, we write $w_1 = |y|^2 w_{*1}$, $w_2 = |y|^2 w_{*2}$, 
and $u =  w_1 - w_2$. By a simple subtraction, $u$ satisfies a linear equation 
in the form of \eqref{eq-linear-equation-general-special-factors}, i.e., 
\begin{align}\label{eq-linear-equation-conti-depend-MSE}\begin{split} 
\Delta u &= a|y|^2\Delta u+\widetilde{a} y_ky_\ell\partial_{k\ell}u+a^k|y|^2 y_\ell \partial_{k\ell}u
+a^{k\ell}|y|^4 \partial_{k\ell}u\\
&\qquad+ b  \langle y, \nabla u\rangle+ b^k|y|^2\partial_ku + cu + f,
\end{split}\end{align}
for some analytic functions $a, \widetilde{a}, a^k, a^{k\ell}$, $b, b^k$, $c$, and $f$ in $B_{R_0}$.
Moreover, $f$ can be written as
\begin{align}
\label{eq-identity-f-conti-depend-MSE}
f = \widetilde a^{ij} \partial_{ij}(h_1 - h_2) + \widetilde b^i \partial_i (h_1 - h_2) + \widetilde{c} (h_1 - h_2),
\end{align}
for some analytic functions $\widetilde a^{ij}$, $\widetilde b^i$, and $\widetilde c$ in $B_{R_0}$. 
We point out that $a, \widetilde{a}, a^k, a^{k\ell}$, $b, b^k$, $c$, and $\widetilde a^{ij}, \widetilde b^i, \widetilde c$
are analytic functions in $y$, $h_1, y\nabla h_1, y^2\nabla^2 h_1$, 
$h_2, y\nabla h_2, y^2\nabla^2 h_2$, $w_1, y\nabla  w_1, y^2\nabla^2  w_1$, and
$ w_2, y\nabla  w_2, y^2\nabla^2 w_2$, for $k, \ell=1, \cdots,n$ and $i,j = 1,\cdots,n$.
For an illustration, we compute the difference from $vy_iy_j\partial_{ij}v$.
In this case, with $u = w_1 - w_2$, we have
\begin{align*}
&v_1y_iy_j\partial_{ij}v_1 - v_2 y_iy_j\partial_{ij}v_2\\
&\qquad= (w_2 y_iy_j +  h_2 y_i y_j) \partial_{ij} u  + (y_iy_j \partial_{ij} h_1  + y_iy_j \partial_{ij}  w_1) u\\
&\qquad\qquad+ (w_2 y_iy_j  + h_2 y_i y_j) \partial_{ij}(h_1 - h_2) + (y_iy_j \partial_{ij}h_1 
 + y_iy_j\partial_{ij}w_1)(h_1 - h_2).
\end{align*}

Now, assume that $\{H_{1,k}\}_{k\ge 0}$ and $\{H_{2,k}\}_{k\ge 0}$ are the sequences of 
the homogeneous Taylor polynomials of $h_1$ and $h_2$ at $0$, respectively, with $H_{1,k}, H_{2,k} \in \mathscr{H}_k$. 
Then, by the interior estimate for harmonic functions, we have, for any $k \ge 0$,
\begin{align}
\label{eq-upper-bound-harmonic-part-homogeneous-conti-deped-MSE}
\|H_{1,k}\|_*, \|H_{2,k}\|_* \leq M_0 K_0 K^k
\quad\text{and}\quad\|H_{1,k} - H_{2,k}\|_* \leq \|h_1 - h_2\|_{L^2(B_1)} K_0 K^k,
\end{align}
for some positive constants $K_0$ and $K$, depending only on $n$. 
Recall that, by Lemma \ref{lemma-T-preserving}, the estimates in Theorem \ref{thrm-existence-linear-equation} 
are valid for $u$. By Lemma \ref{lemma-estimate-composition-of-analytic-functions}, 
\eqref{eq-upper-bound-purely-unharmonic-part-conti-depend-MSE}, \eqref{eq-linear-equation-conti-depend-MSE}, 
and \eqref{eq-upper-bound-harmonic-part-homogeneous-conti-deped-MSE}, 
the constants $A_0$ and $A$ in \eqref{eq-estimate-T_i P-linear} depend only on $n$, $M_0$ and $C_*$. 
Similarly, by \eqref{eq-identity-f-conti-depend-MSE}, 
we can take $B_0 = K_*\|h_1 - h_2\|_{L^2(B_1)}$,
where $K_*$ is a positive constant depending only on $n$ and $M_0$, and $B$ to be a constant 
depending only on $n$, $M_0$, and $C_*$ in \eqref{eq-estimate-Q_j-linear}.
Therefore, the desired estimate \eqref{eq-estimate-difference-conti-depend-MSE} 
follows from \eqref{eq-final-estimate-linear-k}.
\end{proof}

We now make some remarks. Note that \eqref{eq-estimate-difference-conti-depend-MSE} implies, 
for $r > 0$ sufficiently small,
$$\big||y|^2 w_{*1} - |y|^2 w_{*2} \big|_{L^\infty(B_r)} \leq \frac{1}{2}  \|h_1 - h_2\|_{L^2 (B_1)}.$$
On the other hand, let $h$ be a harmonic function in $B_1$ with $\|h\|_{L^2(B_1)} \leq M_0$ and $F_h$ be the nonlinear functional given by 
$$F_h(y, w,y\nabla  w, y^2\nabla^2  w) 
= F(y,h +  w, y \nabla h + y \nabla  w, y^2 \nabla^2 h + y^2 \nabla^2  w).$$
Write 
\begin{align*}
F_h(y,w,y\nabla w, y^2\nabla^2 w) = f_h(y) +  \sum_{i=1}^3\sum_{j=0}^\infty T_{ij}(w,\cdots,w)
\end{align*}
as in \eqref{eq-Taylor-expansion-F_h-MSE}, for some analytic function $f_h$ and operators $T_{ij}$ of type $(i,j)$ 
that preserve the degree of $|y|^2$, for $i=1,2,3$ and $j \ge 0$. Then, by the interior estimate for harmonic functions, 
\eqref{eq-estimate-for-T_ij-fully-nonlinear} holds for some constants $A,B_0,B>0$ depending only on $n$ and $M_0$. 
In addition, if $\{Q_{h,j}\}_{j\ge 0}$ is the sequence of the homogeneous Taylor polynomial of $f_h$ at $0$, with $Q_{h,j} \in \mathscr{P}_j$, then
$$\|Q_{h,j}\|_* \le C_0C^j,$$
for some positive constants $C_0$ and $C$, depending only on $n$ and $M_0$. Let $ w = |y|^2 w_*$ be the unique analytic function 
in Theorem \ref{thrm-harmonic-fcn-correspondence-MSE} 
such that $v = h+ w$ is a solution to \eqref{eq-nonlinear-Poission-MSE-fully-1}. 
Assume that $\{P_k\}_{k\ge 2}$ is the sequence of the homogeneous Taylor polynomials of $w$ at $0$, with $P_k \in \mathscr{H}_k^\perp$. 
Then, by Theorem \ref{thrm-existence-fully-nonlinear} , we get, for any $k \ge 2$,
$$\|P_k\|_* \leq C_*^k,$$
for some constant $C_*>0$, depending only on $n$ and $M_0$. 
In other words, if $h_1$ and $h_2$ are harmonic functions in $B_1$ with
$$\|h_1\|_{L^2(B_1)}, \|h_2\|_{L^2(B_1)} \leq M_0, $$ 
then there exists a constant $C_*>0$, depending only on $n$ and $M_0$, 
such that \eqref{eq-upper-bound-purely-unharmonic-part-conti-depend-MSE} holds. 
We emphasize again that the proof of Theorem \ref{thrm-conti-dependence-MSE} 
does not rely on the existence results in Section \ref{sec-MSE}.

Next, we study the dependence of solutions to the  Monge-Amp\`{e}re equation in exterior domains on harmonic functions. 
If the dimension $n$ is even, the discussion is similar to that of the minimal surface equation. 
More precisely, consider
$$\Delta v = F,$$
where $F$ is given by \eqref{eq-MA-expression-F-for-v}. 
Since it has the same structure as \eqref{eq-nonlinear-Poission-MSE-fully-1}, 
we can proceed as in the proof of Theorem \ref{thrm-conti-dependence-MSE} to obtain the same conclusion. 
The discussion is more complicated when the dimension $n$ is odd. 
In this case, we need to study a system of the form \eqref{eq-MA-dim-odd-system}, i.e.,
\begin{align*}
\begin{split}
\Delta \widetilde{v} &= \widetilde{F}\big(y,\widetilde{v}, \widehat{v}, y \nabla \widetilde{v}, y \nabla \widehat{v}, 
y^2 \nabla^2 \widetilde{v}, y^2 \nabla^2 \widehat{v} \big),\\
\Delta (|y|^{n-2} \widehat{v} )&= |y|^{n-2} \widehat{F}\big(y,\widetilde{v}, \widehat{v}, y \nabla \widetilde{v}, y \nabla \widehat{v}, 
y^2 \nabla^2 \widetilde{v}, y^2 \nabla^2 \widehat{v} \big).
\end{split}
\end{align*}
In this case, we generalize Theorem \ref{thrm-existence-linear-equation} to elliptic systems 
and employ a similar argument as in the proof of Theorem \ref{thrm-conti-dependence-MSE} 
to obtain the Lipschitz dependence on harmonic functions as in \eqref{eq-estimate-difference-conti-depend-MSE}.
We omit details.


\begin{thebibliography}{99}

\bibitem{Aiolfi&Bustos2022} A. Aiolfi, D. Bustos, J. Ripoll,
{\it On the existence of foliations by solutions to the exterior
{D}irichlet problem for the minimal surface equation.}
Proc. Amer. Math. Soc., 150(2022), 3063-3073.

\bibitem{Alm1966} F. Almgren,
{\it Some interior regularity theorems for minimal surfaces and an extension of Bernstein's theorem.}
Ann. of Math., 84(1966), 277-292.

\bibitem{HarmonicFcnTheory} S. Axler, P. Bourdon, W. Ramey,  
{\it Harmonic Function Theory}, 2nd ed.,
Graduate Texts in Mathematics, Springer, 2001. 

\bibitem{BaoLi2012} J. Bao, H. Li,
{\it On the exterior Dirichlet problem for the Monge–Ampère equation in dimension two.} Nonlinear Anal.: Theory, Methods and Applications,
Vol. 75, Iss. 18(2012), 6448-6455,


\bibitem{BaoLiu2022} J. Bao, Z. Liu, 
{\it Asymptotic expansion at infinity of solutions of special Lagrangian equations.}
J. Geom. Anal., 32(2022), 90.

\bibitem{Bernstein1912}
S. Bernstein,
{\it Sur l'ordre de la meilleure approximation des fonctions continues
par des polyn\^omes de degr\'e donn\'e.}
M\'em. Acad. Roy. Belg. Cl. Sci., \textbf{4} (1912), 1--103.

\bibitem{Bern1015} S. Bernstein,
{\it Sur un th\'eor\`eme de g\'eom\'etrie et ses applications aux \'equations aux d\'eriv\'ees partielles du type elliptique.}
Comm. de la Soc. Math. de Kharkov, 15(1915-1917), 38-45; also German translation:
\"Uber ein geometrisches Theorem und seine Anwendung auf die partiellen Differentialgleichungen vom elliptischen Typus.
Math. Z., 26(1927), 551-558.

\bibitem{BL} L. Bers,
{\it Isolated singularities of minimal surfaces.} 
Ann. of Math., 53(1951), 364-386.

\bibitem{BGG1969} E. Bombieri, E. De Giorgi, E. Giusti,
{\it Minimal cones and the Bernstein theorem.} 
Inventiones Math., 7(1969), 243-269.

\bibitem{Caffarelli1995} L. Caffarelli,
{\it Topics in PDEs: The Monge-Amp\`{e}re equation.} Graduate course. Courant Institute,
New York University, 1995.

\bibitem{CaLi2003} L. Caffarelli, Y. Li,
{\it An extension to a theorem of J\"{o}rgens, Calabi, and Pogorelov.} Comm. Pure Appl. Math.,  56(2003), 549-583.

\bibitem{Calabi1958} E. Calabi,
{\it Improper affine hyperspheres of convex type and a generalization of a theorem by K.
J\"orgens.} Michigan Math. J.,  5(1958), 105-126.

\bibitem{DeGiorgi1965} E. De Giorgi,
{\it Una extensione del teorema di Bernstein.} Ann. Scuola Norm. Sup. Pisa, 19(1965), 79-85.

\bibitem{FeMaMi1999} L. Ferrer, A. Mart\'{i}nez, F. Mil\'{a}n,  
{\it An extension of a theorem by K. J\"orgens and a maximum
principle at infinity for parabolic affine spheres.} Math. Z., 230(1999), 471-486.


\bibitem{Han-Marchenko202?} Q. Han, I. Marchenko,
{\it Solutions of  the special Lagrangian equation near infinity}, preprint.

\bibitem{Han-Wang2024} Q. Han, Z. Wang,
{\it Solutions of the minimal surface equation and of the Monge-Amp\`{e}re equation near infinity}.
J. Reine Angew. Math. (Crelles Journal),
2025(2025), no. 820, 51-73.

\bibitem{Helgason1959}
S. Helgason,
{\it Differential operators on homogeneous spaces},
Acta Math. 102(1959), 239-299.

\bibitem{HongYuan2020} G. Hong, Y. Yuan,
{\it Maximal hypersurfaces over exterior domains.} Comm. Pure Appl. Math., 74(2021), 589-614.

\bibitem{Jorgens1954} K. J\"orgens,
{\it \"Uber die Lösungen der Differentialgleichung $rt -s^2 = 1$.} Math. Ann., 127(1954),
130-134.

\bibitem{KrantzParks2002} S. Krantz, H. Parks,
{\it A Primer of Real Analytic Functions 2nd ed.} Birkh{\"a}user, 2002.

\bibitem{Krust1989} R. Krust,
{\it Remarques sur le problème extérieur de Plateau.} Duke Math. J., 59(1989), 161-173.

\bibitem{KutevTomi1998} N. Kutev, F. Tomi,
{\it Existence and nonexistence for the exterior Dirichlet problem for the minimal surface equation in the plane.} 
Diff. and Integ. Eq.,
11(1998), 917-928.

\bibitem{Kuwert1993} E. Kuwert,
{\it On solutions of the exterior Dirichlet problem for the minimal
 surface equation.} Ann. Inst. Henry Poincar\'e, 10(1993), 445-451.

\bibitem{LLY2020} D. Li, Z. Li, Y. Yuan,
{\it A Bernstein problem for special Lagrangian equations in exterior domains.}
Adv. Math., 361(2020), 106927.

\bibitem{BaoLiu2020} Z. Liu, J. Bao, 
{\it Asymptotic expansion at infinity of solutions of Monge-Amp\`{e}re type equations.}
Nonlinear Anal., 212(2021), Paper No. 112450, 17 pp.
 
\bibitem{Mo} J. Moser,
{\it On Harnack's theorem for elliptic differential equations.} Comm. Pure Appl. Math., 14(1961), 577-591.


\bibitem{Osserman1986} R. Osserman, 
{\it A Survey of Minimal Surfaces.} Dover publ. 2nd ed., New York, 1986.

\bibitem{Pogorelov} A. Pogorelov,
{\it On the improper convex affine hyperspheres.} Geometriae Dedicata, 1(1972),
33-46.

\bibitem{RahmanSchmeisser}
Q. Rahman and G. Schmeisser,
{\it Analytic Theory of Polynomials.}
London Mathematical Society Monographs, New Series, Vol.~26,
Oxford University Press, Oxford, 2002.

\bibitem{Sc} R. Schoen,
{\it Uniqueness, symmetry, and embeddedness of minimal surfaces.}
J. Diff. Geom., 18(1983), 791-809.


\bibitem{Seeley1969}
R. T. Seeley,
{\it Eigenfunction expansions of analytic functions.}
Proc. Amer. Math. Soc., 21(1969), 734--738.

\bibitem{Sim} L. Simon,
{\it Asymptotic behavior of minimal graphs over exterior domains.} 
Ann. Inst. Henri Poincar\'{e}, Anal. Non Lin\'{e}aire, 4(1987), 231-242.

\bibitem{Simon1994} L. Simon,
{\it Asymptotics for exterior solutions of quasilinear elliptic equations. 
Geometry from the Pacific Rim: 
Proceedings of the Pacific Rim Geometry Conference held at National University of Singapore, 
Republic of Singapore, December 12-17, 1994,} edited by A. Jon Berrick, B. Loo and H.-Y. Wang, Berlin, 
De Gruyter, 1997, 343-362. https://doi.org/10.1515/9783110908961.343
                    
\bibitem{Simons1968} J. Simons,
{\it Minimal varieties in Riemannian manifolds.} Ann. of Math., 88(1968), 62-105.

\bibitem{Walsh1926}
J. Walsh,
{\it \"Uber den Grad der Approximation einer analytischen Funktion.}
Sitzungsberichte der Bayerischen Akademie der Wissenschaften,
Mathematisch-Naturwissenschaftliche Abteilung, 1926, 223--229.

\bibitem{WangBao2013} C. Wang, J. Bao,
{\it Necessary and sufficient conditions on existence and convexity of solutions 
for Dirichlet problems of Hessian equations on exterior domains.}
Proc. Amer. Math. Soc., 141(2013), 1289-1296.

\bibitem{Warren2010} M. Warren,
{\it Calibrations associated to Monge-Amp\'{e}re equation.} Trans. American Mathematical
Society. 362 (2010), 3947-3962.




\end{thebibliography}
\end{document}